\documentclass[11pt]{article}
\pdfoutput=1 
\usepackage{amssymb, amsmath, amsthm}
\usepackage[margin=1in]{geometry}

\usepackage{tikz}
\usepackage{booktabs}
\usepackage{multirow} 
\usepackage{tabularx}

\usepackage{natbib}
 \bibpunct[, ]{(}{)}{,}{a}{}{,}%

\usepackage[colorlinks=true,linkcolor=blue,citecolor=blue]{hyperref}
\usepackage{tabularray}
\usepackage[utf8]{inputenc}
\usepackage{tikz}
\usepackage{pgfplots}
\pgfplotsset{compat=newest}
\usepgfplotslibrary{statistics}
\usepackage{tabularray,subcaption}

\usepackage{algorithm}
\usepackage{algorithmic}

\pgfplotsset{
    every axis/.style={
        line width=0.8pt,
        tick style={line width=0.7pt,black},
        grid=both,
        grid style={dotted,gray!50},
        minor grid style={dotted,gray!25},
        minor tick num=1,
        tick align=inside,
        tick label style={font=\small},
        label style={font=\small},
        legend style={
            font=\small,
            draw=none,
            row sep=3pt
        }
    },
    every axis plot/.style={
        line width=1pt,
        mark options={solid}
    }
}

\newtheorem{property}{Property}
\newtheorem{remark}{Remark}[section]

\usepackage[nameinlink]{cleveref}
\crefname{property}{property}{properties}
\Crefname{property}{Property}{Properties}
\newtheorem{theorem}{Theorem}

\usepackage{authblk}
\providecommand{\keywords}[1]{\noindent\textbf{Keywords:} #1}

\begin{document}

\title{Extreme-Scale EV Charging Infrastructure Planning for Last-Mile Delivery Using High-Performance Parallel Computing}


\author[a]{Waquar Kaleem}
\author[b]{Taner Cokyasar}
\author[b]{Jeffrey Larson}
\author[b]{Omer Verbas}
\author[c]{Tanveer Hossain Bhuiyan}
\author[a]{Anirudh Subramanyam}

\affil[a]{Pennsylvania State University, University Park, PA 16802, USA}
\affil[b]{Argonne National Laboratory, 9700 S. Cass Avenue, Lemont, IL, 60439 USA}
\affil[c]{University of Texas at San Antonio, One UTSA Circle, San Antonio, TX 78249, USA}

\maketitle


\begin{abstract}
This paper addresses stochastic charger location and allocation problems under queue congestion for last-mile delivery using electric vehicles (EVs). The objective is to decide where to open charging stations and how many chargers of each type to install, subject to budgetary and waiting-time constraints. We formulate the problem as a mixed-integer non-linear program, where each station–charger pair is modeled as a multiserver queue with stochastic arrivals and service times to capture the notion of waiting in fleet operations. The model is extremely large, with billions of variables and constraints for a typical metropolitan area; even loading the model in solver memory is difficult, let alone solving it. To address this challenge, we develop a Lagrangian-based dual decomposition framework that decomposes the problem by station and leverages parallelization on high-performance computing systems, where the subproblems are solved by using a cutting plane method and their solutions are collected at the master level. We also develop a three-step rounding heuristic to transform the fractional subproblem solutions into feasible integral solutions. Computational experiments on data from the Chicago metropolitan area with hundreds of thousands of households and thousands of candidate stations show that our approach produces high-quality solutions in cases where existing exact methods cannot even load the model in memory. We also analyze various policy scenarios, demonstrating that combining existing depots with newly built stations under multiagency collaboration substantially reduces costs and congestion. These findings offer a scalable and efficient framework for developing sustainable large-scale EV charging networks.
\end{abstract}

\keywords{
Facility Location, Stochastic Charger Location and Allocation, Capacity Allocation, Charging Station, Electric Vehicles, Queuing Systems, Mixed-Integer Quadratic Programming, Large-Scale Optimization}

\section{Introduction}\label{sec:Intro}
The global electric vehicle (EV) market has experienced substantial growth, with the fleet size expanding from 11 million in 2020 to approximately 40 million vehicles in 2023 \citep{IEA2024}. EV sales accounted for nearly 18\% of total car sales in 2023, with projections indicating continued growth to 17 million units in 2024. By 2030, EVs could represent 35\% of global car sales, or even 42\% to 58\% if all manufacturers meet their electrification targets \citep{IEA2024}. To support this increasing demand, government entities and private organizations worldwide must ensure adequate charging infrastructure. The National Renewable Energy Laboratory study \citep{wood20232030} projects that supporting 33 million EVs in the United States by 2030 will require approximately 28 million charging ports of various levels (fast, moderate, and slow). However, depot‐only charging is insufficient for high‐utilization last‐mile fleets: medium‐duty delivery vans and trucks often need opportunistic en‐route fast‐charging infrastructure to sustain continuous operations \citep{nrel,EDF2024ChargingNeeds}.  
Indeed, studies have shown that a lack of adequate EV infrastructure could hamper EV adoption and fleet growth \citep{DolsakPrakash2021}.

Our work is particularly motivated by the ongoing transformation of urban last-mile delivery through electric vehicle adoption \citep{fareye2024}. As e-commerce companies seek sustainable and cost-effective solutions for urban logistics, major carriers are making ambitious commitments to electrification. For example, FedEx has announced plans to convert its entire parcel pickup and delivery fleet to zero-emission electric vehicles by 2040 \citep{fedex2023}. EVs have emerged as an attractive option due to their lower operational costs, reduced maintenance requirements, and significant fuel savings. While this transition promises environmental benefits through reduced carbon emissions and improved air quality, it also presents unique operational challenges. Delivery vehicles require reliable access to charging infrastructure while maintaining delivery schedules and service levels. Class 6 last-mile delivery trucks---larger than cargo vans yet more maneuverable than semi-tractor units---often exceed their planned daily range and must top up at public or semi-public fast chargers during their shifts to avoid missed delivery windows \citep{nrel,EDF2024ChargingNeeds}. The stochastic nature of delivery routes, varying package volumes, and time-sensitive operations make it crucial to strategically locate charging stations—both at depots and en route along delivery corridors—that can handle peak demand periods without causing significant delays in delivery operations. These considerations are especially critical as businesses increasingly scale their electric delivery fleets to meet growing consumer demand for sustainable delivery options.

The EV charging station location problem is a subset of facility location problems that seeks to optimize the locations of charging stations to minimize various costs, including station opening, charger placement, waiting times, and accessibility \citep{kchaou2021charging}. Traditionally, this problem has been modeled assuming 
deterministic demand for electric vehicles and fixed service times \citep{davatgari2024electric}, not necessarily incorporating expected waiting time and queue time in the system. In addition, previous studies have focused mainly on small-scale scenarios (in terms of the number of households, charging stations, and chargers considered), limiting their applicability to the large-scale infrastructure needs projected for the future. 

Our paper presents a new approach to the EV charging station location problem by formulating it as a stochastic location model with congestion and immobile servers (SLCIS) (e.g., see \citet{bermanandkrass}), which we refer to as the stochastic charger location and allocation problem (SCLA). The SLCIS is a facility location model with random service times in which consumers produce streams of random service demands. This results in congestion, as some incoming demands cannot be immediately addressed and must wait in the queue \citep{hale2003location}. By regarding EV charging stations as immobile servers and EV arrivals as stochastic demands, we incorporate queuing theory and formulate the problem as a mixed-integer non-linear program (MINLP). 

We present two exact formulations for the resulting SCLA model. When we scale the models to meet future EV infrastructure needs, however, we find that these exact formulations become prohibitively large and computationally infeasible. In fact, building the model alone can exceed hardware memory limits and make direct solution attempts impractical. To address these scalability challenges, we develop a Lagrangian dual decomposition framework \citep{guignard2003lagrangean}, which decomposes the original problem into station-level subproblems that are then solved in parallel using high-performance computing (HPC). We design a customized cutting-plane method to solve the subproblems in parallel on HPC clusters, which reduces both computation time and memory requirements. In addition, we develop novel rounding heuristics that convert fractional Lagrangian solutions into feasible solutions for the exact model. Our approach yields meaningful solutions for extremely large-scale SCLA problem instances with hundreds of thousands of households and thousands of candidate stations. To the best of our knowledge, such large problems have never been  addressed in the literature. Our specific contributions can be summarized as follows: 

\begin{itemize}
    \item We formulate the SCLA problem as a stochastic location model with congestion and immobile servers (chargers) that captures stochasticity in EV charging demand, charging rates, and queue waiting times, thus constituting a realistic framework for large-scale EV infrastructure planning. 
    
    \item We propose two exact solution methods for the SCLA problem based on mixed-integer nonlinear programming formulations, both of which incorporate cutting-plane methods to handle the nonlinear waiting time constraints.
    
    \item To address the scalability limitations of the exact methods, we design a Lagrangian dual decomposition framework. By decomposing the problem into station-level subproblems and solving them in parallel on HPC clusters, we significantly enhance scalability and computational efficiency for extreme-scale SCLA instances. We also design a three-step rounding heuristic to transform the fractional solutions into feasible integral solutions for the exact formulations. 
    
    \item We implement and demonstrate the proposed framework in the domain of last-mile e-commerce delivery. Through extensive computational experiments and policy analyses, we show that our approach can handle instances far beyond the capabilities of existing exact methods, providing actionable insights for EV infrastructure planning.
    
\end{itemize}

The remainder of this paper is organized as follows. \Cref{sec:Lit_review} reviews the relevant literature, \Cref{sec:problem_def} presents the problem formulation, \Cref{sec:exact_solution_method} and \Cref{sec:lagrengean_based_framework} detail the solution approaches, \Cref{sec:comp_exp} presents findings from our computational experiments, and \Cref{sec:Conclusion} offers concluding remarks. 

\section{Literature Review}\label{sec:Lit_review}
Stochastic location models with congestion and immobile servers  have seen significant developments over the past few decades \citep{bermanandkrass}. This section summarizes some key contributions, focusing on queuing systems, objective functions, solution methods, and problem sizes addressed. For a concise overview of the queuing notation used in this paper, see \citet[Chapter~11.1.3]{stewart2009}. \citet{amiri1997solution} used an \textit{M}/\textit{M}/1 queuing system, with the objective to minimize facility opening, detour, and waiting costs. The author employed Lagrangian relaxation \citep{guignard2003lagrangean} to solve instances with up to 500 user nodes and 40 potential facilities, reporting average CPU times of up to an hour.
In another study, \citet{wang2002algorithms} also utilized an \textit{M}/\textit{M}/1 system with the objective to minimize the expected traveling and waiting costs. They applied greedy heuristics, tabu search \citep{glover1998tabu}, and $\epsilon$-optimal branch-and-bound methods to problems with 459 customer nodes and 84 potential facility sites, reporting solution times of roughly 1,000 seconds. \citet{elhedhli2006service} extended the \textit{M}/\textit{M}/1 model to include server capacity costs and used linearization techniques with piecewise linear approximations to solve instances with 100 households and 20 potential service facilities, reporting solution times of up to ten minutes.

Other studies have developed various reformulations to solve SLCIS models. For example, \citet{Goez_2017} reformulated \textit{M}/\textit{M}/1 systems as a mixed-integer second-order cone optimization problem, solving instances with up to 200 customer nodes and 30 facility locations in roughly half an hour. 
Similarly, \citet{Elhedhli_2018} developed an approach for \textit{M}/\textit{M}/1 systems by reformulating the problem as a mixed-integer quadratic program with fourth-degree polynomial constraints and using Lagrangian relaxation, solving problems with up to 100 customer nodes and 15 facility locations, with service rates as continuous decision variables, in about 2.77 hours.  
In studies by Ahmadi-Javid and collaborators \citep{ahmadi2018location,ahmadi2020linear,ahmadi2022convexification}, the authors used an \textit{M}/\textit{G}/1 model and  developed linear formulations, valid inequalities, and mixed-integer second-order cone programming reformulations for solving SLCIS models with up to 400 customer nodes and 30 facility locations within a time limit of three hours. Likewise, \citet{vidyarthi2014efficient,Vidyarthi_2015} adopted an \textit{M}/\textit{G}/1 queuing system and used linearization with an exact constraint generation algorithm to solve problems with up to 500 customer nodes and 40 facility locations. \citet{Etebari_2019} adopted an \textit{M}/\textit{M}/1/\textit{k} model to maximize system profit, employing column generation and hybrid metaheuristics for problems with up to 105 customer nodes and 42 facility locations, with solution times of roughly two hours. 

\citet{Syam_2008} employed linearization techniques followed by Lagrangian methods \citep{guignard2003lagrangean} for problems with up to 250 districts and 60 facilities, achieving average solution times of about three minutes. 
Similarly, \citet{aboolian2008location} and \citet{Aboolian_2009} used an \textit{M}/\textit{M}/\textit{s} system, proposing descent heuristics, simulated annealing, and an exact iterative method consisting of solving an uncapacitated facility location problem. 
They addressed larger instances with up to 800 customer nodes and 167 potential facility locations, with solution times reaching up to 16 hours. 
\citet{castillo2009social} utilized asymptotic approximations for \textit{M}/\textit{M}/\textit{s} systems, originally provided by \citet{halfin1981heavy} and later extended by \citet{borst2004dimensioning}. The authors applied their model to a case study in Edmonton, Canada, considering 38 facility locations and 222 centroids of household neighborhoods. Their work demonstrated that multiserver facilities could be modeled and solved without increased computational effort compared with single-server facilities.  \citet{aboolian2022efficient} revisited the \textit{M}/\textit{M}/1 model, comparing several approaches including generalized Benders decomposition. They addressed problems with roughly 100 customer nodes and facility locations and 10 service capacity levels, with solution times up to two hours. 

This body of research demonstrates the evolution of solution methods applied to SLCIS, where queuing systems have been modeled by using \textit{M}/\textit{M}/1, M/G/1, and \textit{M}/\textit{M}/\textit{s} models. These models have been applied to diverse settings such as ATM facilities, hospitals, and other service systems. Solution methods have progressed from heuristics to mathematical programming techniques that have enabled  solving  instances with up to a few hundred households and potential locations. 

Numerous studies have explored the optimal placement of electric vehicle charging stations, often relying on various optimization methods \citep{kchaou2021charging}.
For instance, \citet{luo2020electric} used the maximal covering location problem combined with queuing theory to enhance the resource utilization of charging stations in sustainable cities. In contrast, \citet{chen2020optimal} proposed a bilevel mathematical model designed to determine station locations while minimizing construction costs and driver waiting times. In a related effort, \citet{cui2019electric} formulated a mixed-integer nonlinear program to minimize the cost of charging stations, network expansion, voltage regulation, and protection device upgrades, thereby catering to urban EV demands.  \citet{erdougan2022establishing} built on this work by introducing an optimization-based framework that locates fast charging stations along designated corridors. Additional works account for the inherent stochasticity in the system. \citet{davidov2017stochastic} developed a stochastic optimization model that incorporates drivers' uncertain travel patterns, vehicles' driving ranges, and infrastructure service quality. \citet{xi2013simulation} focused on maximizing the usage of existing charging facilities in central Ohio through a simulation-optimization procedure. \citet{jordan2022electric} similarly employed a simulation-optimization approach, augmented by a genetic algorithm, to pinpoint ideal station locations and then examined how these solutions affect overall EV waiting times and charging station utilization.

Building on this existing body of work, the present paper applies SLCIS methodologies to the placement of electric vehicle charging stations and the allocation of chargers in last-mile delivery settings. Our approach adopts an \textit{M}/\textit{M}/\textit{s} system for this application, an area not yet explored in this context. Notably, and in contrast to existing literature, we design new solution methods that can handle extremely large-scale instances, including hundreds of thousands of households and thousands of potential facilities with varying server configurations. Solving such large-scale problems is key to enabling strategic decision-making in last-mile delivery, where charger placement can have critical long-term impacts on EV adoption and sustainability.

\section{Model Description}
\label{sec:problem_def}

Consider a graph $(\mathcal{V}, \mathcal{E})$, where $\mathcal{V}$ represents a set of nodes and $\mathcal{E}$ represents edges between nodes. Let $\mathcal{I} \subset \mathcal{V}$ be a discrete set of households, let $\mathcal{J} \subset V$ be a discrete set of candidate locations for charging stations (\textit{stations} for short), and let $ \mathcal{K}$ be the set of charger types with various power outputs.
The key decision variables are related to location and allocation. Specifically, define $y_j$ to be a binary variable indicating whether a station $j \in \mathcal{J}$ is opened.
Also, define $s_{jk}$ to be an integer variable representing the number of chargers of type $k \in  \mathcal{K}$ that are installed at station $j \in \mathcal{J}$.
Define $x_{ijk}$ to be a binary variable that denotes whether a charger of type $ k\in  \mathcal{K}$ at station $j \in \mathcal{J}$ is assigned to EVs that need recharging (immediately) before they serve household $i \in \mathcal I$. 
Given the need for a long-term strategic solution, we adopt a stochastic model of delivery operations that we describe in the following subsections.

\subsection{Modeling Household Demand}
\label{sec:modeling_household_demand}
We assume that each household $i \in \mathcal I$ demands delivery at a rate of $\gamma_i$ per unit time.
Given some time interval $t$ (e.g., one week), we may thus interpret $\gamma_i t$ as the average number of times household $i \in \mathcal I$ is visited by a delivery vehicle in that time interval. We assume that the demands follow Poisson processes with rates \(\gamma_i\) and that they are independent across households. Over a planning horizon, it is standard to model delivery requests as independent Poisson streams \citep{bertsimas1991stochastic}. In particular, by the thinning property of Poisson processes, if all delivery‐request arrivals form a single Poisson stream and each request is independently assigned to household $i$ with probability $\gamma_i/\sum_j\gamma_j$, then each household’s arrival process is itself Poisson (rate $\gamma_i$) and these substreams are mutually independent.One can also parameterize $\gamma_i$ on temporal features such as the day of the week or season, in an attempt to reflect household consumption patterns. For simplicity, we do not consider such dependencies, since historical data suggest that over a strategic planning horizon spanning several months to years, the average demand rates of households can be assumed roughly constant. Empirical studies of service systems likewise find that, once predictable time‐of‐day and day‐of‐week effects are removed, Poisson arrival models capture the key variability in demand \citep{brown2005statistical}.

Let $\pi_i$ denote the probability of a charging activity occurring immediately before the delivery to household $i \in \mathcal{I}$. The average number of charging activities on the route to household $i$ and the average number of deliveries on the route to household $i$ are denoted by $N_i^C$ and $N_i^D$, respectively. These parameters can be derived from historical delivery data. 
The parameter $N_i^D$ can be obtained directly from conventional vehicle routes. Meanwhile, $N_i^C$ can be approximated by dividing each route into segments according to an assumed EV range.
In the absence of detailed state-of-charge (SOC) information for the EVs, the probability $\pi_i$ can be estimated as $\frac{N_i^C}{N_i^D}$. Let $\lambda_i = \gamma_i\pi_i$ denote the charging rate (also referred to as \textit{arrival rate}) for the vehicle delivering to household $i$. Since $\gamma_i$ follows a Poisson process and $\pi_i$ is estimated from historical data, $\lambda_i$ can also be assumed to follow a Poisson process and to be independent for each household $i$. Moreover, the aggregation of many independent charging demand streams at each station is well approximated by a Poisson process, consistent with the Palm–Khintchine principle \citep{bose2013introduction,albin1984approximating,medhi2002stochastic,whitt2002stochastic}.

Recall that the number of chargers of type $k \in  \mathcal{K}$ installed at station $j \in \mathcal{J}$ is denoted as $s_{jk}$. The charging time of charger type $k$ is assumed to follow an exponential distribution with rate $\mu_k$, similar to prior EV‐charging queueing studies \citep{xie2018long,xie2021integrated,kinay2023charging}. In other words, $\mu_k$ is the average number of EVs that can be charged per unit time. The charging service rate $\mu_k$ depends on the initial and final SOC for the EVs immediately before and after charging, which we assume to be homogeneous across all vehicles, since the actual vehicle routes---and hence, their SOC---are not known at the time of decision-making. 

\subsection{Modeling EV Operations}
We assume that each household $i \in \mathcal I$ is served by an EV that always starts from the depot facility $F_i \in \mathcal{J}$, which is co-located with an existing charging station. This assumption reflects many real-world operations, especially in parcel delivery, where households are preassigned to an existing depot. In practice, last-mile delivery vehicles, for example, Class 6 trucks routinely exceed their planned daily ranges and rely on opportunistic fast-charging during their shifts to avoid missed delivery windows. This mid-route topping-up often occurs at public or semi-public stations along delivery corridors rather than at depots alone \citep{nrel,EDF2024ChargingNeeds}. 
To that end, we introduce a parameter $k_c$, which specifies the number of nearest candidate stations that an EV can consider for recharging immediately before serving household $i \in \mathcal I$. Specifically, we assume that an EV can recharge only at one of these $k_c$ nearest stations. 
This parameter $k_c$ reflects the vehicle’s flexibility, trading off between detouring to nearby but congested stations versus remote but relatively unoccupied ones.  To model this flexibility, let subset $\mathcal{J}_i \subseteq \mathcal{J}$ represent the set of feasible stations based on current SOC that an EV can potentially visit immediately before it serves household $i \in \mathcal{I}$, allowing it to charge before delivering to household $i$. 
Specifically, we define $\mathcal{J}_i$ as the set of $k_c$ stations that are closest to $i$ with respect to (lat, long) coordinates. Conversely, we also define the subset of households $\mathcal{I}_j \subseteq \mathcal{I}$ that could be visited after charging at a particular station $j \in \mathcal{J}$. More precisely, $\mathcal{I}_j$ is the set of households for which $j$ is among their $k_c$ nearest stations. In other words, $\mathcal{I}_j = \{ i \in \mathcal{I} : j \in \mathcal{J}_i \}$. 

\begin{remark}[Neighborhood-level risk pooling]
The parameter $k_c$ implicitly enables a form of neighborhood-level risk pooling because, 
during the optimization, each household $i$ can only be assigned to a station $j \in \mathcal{J}_i$, 
where $\mathcal{J}_i$ contains the $k_c$ nearest stations to $i$. 
This restricts assignments to nearby stations, preventing unrealistically long detours 
and allowing the model to balance charging demand locally across neighboring stations.
\end{remark}

The sets $\mathcal{J}_i$ and $\mathcal{I}_j$ are crucial to limit the feasible options for each household and station, respectively. This approach replaces the commonly adopted \emph{coverage constraints} \citep[Ch. 17, p. 486]{bermanandkrass}, eliminating the need for additional constraints to restrict the stations that could be visited prior to each household. Considering all stations as potential charging locations prior to servicing $i$ would not only explode the problem size but also be unrealistic, as the SOC may not allow visiting any given station in the network. Constructing the sets $\mathcal{J}_i$ and $\mathcal{I}_j$ is explained in the Appendix as \Cref{alg:construct_Ji_Ij}.

Let $T_{ij}$ denote the travel time from household $i \in \mathcal{I}$ to station $j \in \mathcal{J}$. Assuming the triangle inequality, the EV detour travel time to visit the station $j \in \mathcal{J}_i$ is then defined as $T^{\delta}_{ij} := T_{F_{i}j} + T_{ji} - T_{F_{i}i}$ for all $i \in \mathcal{I}$ and $j \in \mathcal{J}_i$. We emphasize that this detour travel time is simply a route-agnostic approximation to the true travel time of the EV that detours immediately before visiting household~$i$. In particular, it depends only on the depot $F_i$ which serves the household and the station~$j$ to which the EV detours. It is nonnegative and equals zero when $j$ lies on a shortest $F_i\!\to\! i$ path. If the predecessor of $i$ on that day’s route were known and equal to $p$, a tour-aware detour would be $\Delta_{pij}:=T_{p j}+T_{j i}-T_{p i}$. However, at the strategic level, the predecessor $p$ is unknown as the routes vary from day to day, so $T^\delta_{ij}$ simply serves as a proxy for charging immediately before $i$. Note that there is no practical constraint that charging must occur before the first stop on the route.

Given the extreme scale of SCLA instances considered in this paper, in our implementation, we first perform a fast nearest-neighbor search using a two-dimensional $k$-d tree on (latitude, longitude) coordinates to select the $k_c$ nearest stations and then compute $T_{ij}^\delta$ only for that subset. We emphasize that $T^{\delta}_{ij}$ serves only as an approximation to the true detour travel time, which is difficult to calculate when routes are unknown. Indeed, during actual operations, a vehicle routing problem will be solved to determine the actual delivery sequences, and each EV will charge as needed before delivering to the next household on the route, allowing a more precise calculation of the detour time \citep{davatgari2024electric}. However, this precise calculation is not possible in our strategic model, since the actual vehicle routes will vary significantly from day to day because of the stochastic nature of household delivery demand, operational uncertainties in travel times, driver availability, and vehicle breakdowns (besides others). We instead adopt a simplified yet realistic model. Crucially, historical delivery and network data can be used to readily inform the parameters of our model, even in the absence of predetermined vehicle routes.

\subsection{Modeling Congestion and Network Closure}\label{sec:model_congestion}
In extended Kendall's notation \citep[Chapter~11.1.3]{stewart2009}, we model the queue at each charger type $k$ at station $j$ as an $M/M/s_{jk}/\infty/\infty$/FCFS system, where the first $M$ denotes a Poisson arrival process and the second $M$ denotes exponentially distributed service times. This effectively means that the processes of EV arrivals for recharging and the subsequent service provided by chargers are memoryless and independent. Recall that the term $s_{jk}$ represents the decision variable representing the number of chargers of type $k$ at station $j$.
Also, the first $\infty$ represents an infinite-length (i.e., uncapacitated) queue, and the second $\infty$ indicates an infinite population (i.e., very large EV demand at each station-charger pair).  $FCFS$ indicates a first-come-first-served queuing policy. Moreover, \textit{Burke’s theorem} \citep{burke1956output} states that, for any stable $M/M/1$ queue, the departure process is Poisson with the same rate as its arrival process.  This result extends to $M/M/s_{jk}$ queues \citep{bose2013introduction}, implying that each charger node in our model—a stable $M/M/s_{jk}$ queue—emits a Poisson departure stream.  Consequently, when the charging stations are linked together, they form an open multi‐server Jackson network \citep{jackson1957networks} whose steady‐state joint distribution admits a product‐form solution.  In such a network, the arrival process into every downstream charger‐station $(j,k)$ remains Poisson and independent of upstream activity.  This closure property allows us to model the EV-charging infrastructure at a strategic level when no information about individual routes is available.

Our assumptions imply that EVs form a single queue for each charger type $k$ at station $j$, reflecting that service rates and hardware characteristics differ across types. 
If certain charger types at a charging station are technically identical, their capacities can be combined and modeled as a single effective type without affecting our modeling approach. We also highlight that our focus in this paper is on strategic allocation decisions, while any finer within-station pooling or switching across heterogeneous chargers correspond to an operational-level load-balancing mechanism. Modeling such a level of operational detail in a strategic context is beyond the scope of this paper. 

The expected waiting time, which we denote as $\mathbb{W}$, includes the queuing and charging times at charger type $k$ at station $j$. 
For practical and system efficiency purposes, we impose a prespecified upper bound of $EW$ on $\mathbb W$, reflecting a maximum acceptable limit on the expected waiting time. Let $\bar{\lambda}_{jk}$ denote the total (decision-dependent) EV demand rate at charger $k$ of station $j$.
Then, the utilization $\rho_{jk}$ is the fraction of time that the particular station-charger pair is busy, and it is simply the ratio of the total demand rate to the total charging service rate at station $j$ for charger type $k$:
\[\rho_{jk} = \frac{\bar{\lambda}_{jk}}{s_{jk} \mu_k} = \frac{\sum_{i\in \mathcal{I}_j} \lambda_{i} x_{ijk}}{s_{jk} \mu_k},\]
where we have exploited the fact that the total demand rate $\bar{\lambda}_{jk}$ at station-charger pair $(j, k)$ can be expressed by using the allocation variables as $\bar{\lambda}_{jk} = \sum_{i\in \mathcal{I}_j} \lambda_{i} x_{ijk}$.
For the queue to be stable, we need to ensure $\rho_{jk} < 1$. By introducing a small safety margin, $\epsilon$, we can linearize the constraint $\rho_{jk} < 1$ as follows: $\mu_k s_{jk} (1 - \epsilon) \geq \sum_{i\in \mathcal{I}_j} \lambda_{i} x_{ijk}$. 

The expected waiting time $\mathbb W$ depends on the probability $\mathbb{P}$ that all $s_{jk}$ chargers are busy.
Both $\mathbb W$ and $\mathbb P$ are functions of the utilization ratio $\rho_{jk}$, the number of chargers (servers) $s_{jk}$, and the charging service rate $\mu_k$.
They both admit closed-form expressions, shown below, which can be derived based on standard textbook arguments \citep[Chapter~11.4.1]{stewart2009}.

\begin{align}
    \mathbb{P} \left(\rho_{jk}, s_{jk}\right) &=\frac{(\rho_{jk} s_{jk})^{s_{jk}}}{(1-\rho_{jk}) s_{jk}! (T_1 + T_2)}, \label{eq:probablity_busy} \\
    \mathbb{W} \left(\rho_{jk}, s_{jk}\right) &= \frac{\mathbb{P}\left(\rho_{jk},s_{jk}\right)}{\mu_k s_{jk}(1-\rho_{jk})}
    + \; \frac{1}{\mu_k},
    \label{eq:expected_waiting_time}
\end{align}
where we have defined
\[
    T_1 = \frac{(\rho_{jk} s_{jk})^{s_{jk}}}{(1-\rho_{jk}) s_{jk}!}, \;
    T_2 = \sum_{r=0}^{{s_{jk}}-1}\frac{(\rho_{jk} s_{jk})^r}{r!}.
\]

Finally, at the household or cluster level, delivery demand counts may exhibit 
\textit{under-dispersion} relative to a Poisson process, especially in dense urban areas where many delivery points (e.g., apartment complexes or multi-tenant buildings) receive at most one delivery per day from a carrier. While we assume independent Poisson requests at the household level, we show in Appendix that even if demand counts were under-dispersed and followed a general renewal process, the resulting expected waiting times at each charger would be lower than those given by the $M/M/s_{jk}$ system. Metropolitan areas such as Chicago (see Figure~1a) have heterogeneous household densities—dense urban clusters and sparse suburban regions—so using a conservative $M/M/s_{jk}$ model provides a better approximation across spatial regimes. Hence, the $M/M/s_{jk}$ assumption provided conservative congestion estimates (see the Appendix for details).

\subsection{Modeling Costs and Budget Constraints}
We let $C^{\phi}_j$ be the fixed cost of opening a station $j \in \mathcal{J}$ per unit time, and we let $C^{\xi}_k$ be the cost of charging the EV using a particular charger type $k \in  \mathcal{K}$ per unit time. The detour travel time cost $C^{\delta}$ represents the cost incurred from extending the route because of the need to charge the EV.
We let $C^{\tau}$ denote the charging service cost per unit time.

Additionally, we specify budgets for locating stations and allocating chargers, which we denote as $B^\phi$ and $B^\xi$ per time unit, respectively.
We also let $ \overline{Y} $ represent the maximum number of stations that can be activated, and we let $ \overline{S}_{jk} $ represent the maximum number of type $k$ chargers that can be allocated to station $j$.

To express the objective function more compactly, we define the indexing set 
$ \mathcal{M} := \{ (i,j,k) \mid j \in \mathcal{J},\, i \in \mathcal{I}_j,\, k \in \mathcal{K}\}$.
Then, we can define 
$f_\phi(y) := \sum_{j \in \mathcal{J}} C^\phi_j y_j$ as the station opening cost, 
$f_\xi(s) := \sum_{j \in \mathcal{J}}\sum_{k \in \mathcal{K}} C^\xi_k s_{jk}$ as the charger installation cost, 
$f_\delta(x) := \sum_{(i,j,k)\in \mathcal{M}} \lambda_i C^\delta T^\delta_{ij} x_{ijk}$ as the detour cost, and 
$f_\tau(x,W) := \sum_{(i,j,k)\in \mathcal{M}} \lambda_i C^\tau W_{jk} x_{ijk}$ as the congestion (waiting) cost. Therefore, the total objective can be expressed as
$f_{\phi}(y) + f_{\xi}(s) + f_{\delta}(x) + f_{\tau}(x,W)$.

\subsection{Final Optimization Model}
\label{sec:final_optimization_model}
Given the aforementioned assumptions and definitions, the SCLA problem can now be formulated as the following mixed-integer nonlinear program, which we denote as $\mathbb{G}$.
\Cref{tab:notations_1} and \Cref{tab:notations_2} summarize the key notation that we use in this model.

\begin{align}
    \mathbb{G}: \mathop{\mathrm{minimize}}_{x,y,W,s} \quad & f_{\phi}(y) + f_{\xi}(s) + f_{\delta}(x) + f_{\tau}(x,W) \label{eq:objective_function} \\
    \mathrm{subject\ to} \quad
    & x_{ijk} \leq y_j \quad \forall  j \in \mathcal{J},i \in \mathcal{I}_j, k \in  \mathcal{K}\label{eq:charge_from_only_open_stations} \\
    & \sum_{j \in \mathcal{J}, k \in  \mathcal{K}} x_{ijk} = 1 \quad \forall i \in \mathcal{I}
    \label{eq:one_assignment} \\
    & \mu_k s_{jk} (1-\epsilon) \geq \sum_{i \in \mathcal{I}_j} \lambda_i x_{ijk} \quad  \forall j \in \mathcal{J},  \forall k \in  \mathcal{K}  \label{eq:rho}\\
    & \mathbb{W} \left( \frac{\sum_{i\in \mathcal{I}_j} \lambda_{i} x_{ijk}}{s_{jk} \mu_k}, s_{jk} \right) \leq W_{jk} \quad
     \forall j \in \mathcal{J}, k \in  \mathcal{K} \label{eq:W_def} \\
    & W_{jk} \leq EW + \frac{1}{\mu_k} \quad \forall j \in \mathcal{J}, k \in  \mathcal{K} \label{eq:W_bound} \\
    & f_{\phi}(y) \leq B^\phi\label{eq:budget_station} \\
    & f_{\xi}(s) \leq B^\xi\label{eq:budget_charger} \\
    & \sum_{j \in \mathcal{J}} y_j \leq \overline{Y}\label{eq:station_limit} \\
    & s_{jk} \leq \overline{S}_{jk} \quad \forall j \in \mathcal{J}, k \in  \mathcal{K}\label{eq:charger_limit} 
    \\
    & s_{jk} \leq \overline{S}_{jk} \sum_{i \in \mathcal{I}_j} x_{ijk} \quad  \forall j \in \mathcal{J}, k \in \mathcal{K} \label{eq:charger_allocation} \\
    & y_j \leq \sum_{k \in \mathcal{K}} s_{jk} \quad \forall j \in \mathcal{J} \label{eq:station_open_charger_allocation}\\
    & x_{ijk}, y_j \in \{0, 1\} \quad \forall j \in \mathcal{J}, i \in \mathcal{I}_j, k \in  \mathcal{K} \label{eq:non-neg1}\\
    & s_{jk} \in \mathbb{Z}_{\geq 0} \quad \forall j \in \mathcal{J}, k \in  \mathcal{K}\label{eq:non-neg2}\\
    & W_{jk} \in \mathbb{R}_{\geq 0} \quad \forall j \in \mathcal{J}, k \in  \mathcal{K}\label{eq:non-neg3}
\end{align}

The objective function (\ref{eq:objective_function}) in $\mathbb{G}$ minimizes the total costs for station opening, charger, detour, and congestion (waiting) per time unit. Note that the congestion cost term, $f_\tau(x,W)$, features bilinear products between decision variables, which makes the objective quadratic. Note that $W_{jk}$ is a dummy variable representing the expected wait time to receive charging service from a charger of type $k$ at station $j$.

Constraints (\ref{eq:charge_from_only_open_stations}) ensure that charger $k$ at station $j$ cannot be visited by the vehicle serving $i$ unless station $j$ is opened. Constraints (\ref{eq:one_assignment}) enforce that vehicles serving each household must be charged by exactly one charger of some type $k$ at some station $j$. Constraints (\ref{eq:rho})--(\ref{eq:W_bound}) are related to queue congestion and are based on the equations derived in \Cref{sec:model_congestion}. Specifically, Constraints (\ref{eq:rho}) ensure that the queue stability condition $\rho_{jk} < 1$ is satisfied, that is, the service rate per charger type is strictly greater than its total charging demand rate.
 Additionally,  Constraints (\ref{eq:W_def}) and (\ref{eq:W_bound}) set a lower bound for the service-inclusive waiting time, based on the function $\mathbb{W}$, and an upper bound with $EW$, the maximum acceptable queue waiting time.

Constraints (\ref{eq:budget_station})--(\ref{eq:budget_charger}) ensure that the costs of station locations and charger allocations do not exceed their respective budgets. Constraints (\ref{eq:station_limit}) and (\ref{eq:charger_limit}) limit the number of active stations and the allocation of type $k$ chargers at each station, respectively. Constraints (\ref{eq:charger_allocation}) ensure that a charger of type $k$ can  be allocated at station $j$ only if at least one household $i$ is assigned to it. Constraints (\ref{eq:station_open_charger_allocation}) ensure that a station $j$ can  be opened only if it has at least one charger assigned of any type.  Constraints (\ref{eq:non-neg1}), (\ref{eq:non-neg2}), and (\ref{eq:non-neg3}) define the variable domains. 

\begin{remark}[Probabilistic interpretation of fractional assignments]
The optimization problem $\mathbb{G}$ presented above considers $x_{ijk}$ as binary decision variables. 
However, the relaxation $x_{ijk} \in [0,1]$ is also examined within the decomposition and rounding framework 
(see~\Cref{sec:rounding_heuristics}). 
In which case, $x_{ijk}$ can be interpreted as the long-run probability that household $i$ is served by charger type $k$ at station $j$. 
\end{remark}

\begin{table*}[!ht]
  \caption{Sets and Parameters used in model $\mathbb{G}$. \label{tab:notations_1}}
  \centering\small
  \begin{tabularx}{\textwidth}{@{}l@{\quad}X@{}}
    \toprule
    \textbf{Set}       & \textbf{Definition} \\
    \midrule
    $\mathcal{V, E}$   & vertex set and edge set in the graph representation of the road network.\\
    $\mathcal{I}$      & set of households\\
    $\mathcal{I}_j$    & subset of households that could be visited followed by charging at station $j\in\mathcal{J}$\\
    $\mathcal{J}$      & set of \textit{candidate} charging stations\\
    $\mathcal{J}_i$    & subset of charging stations that could be visited before serving household $i\in\mathcal{I}$\\
    $\mathcal{K}$      & set of charger types\\
    \midrule
    \textbf{Parameter} & \textbf{Definition} \\
    \midrule
    $\epsilon$          & an infinitesimal number\\
    $\lambda_i$         & charging demand rate (arrival rate) for the vehicle delivering to household $i\in\mathcal{I}$, $\lambda_i = \gamma_i\pi_i$\\
    $\gamma_i$          & delivery rate to household $i\in\mathcal{I}$\\
    $\pi_i$             & probability of a charging activity occurring before the delivery to household $i\in\mathcal{I}$, estimated as\\
                        & $ \pi_i = \frac{N^C_i}{N^D_i}$, where $N^C_i$ and $N^D_i$ are counts of charging vs.\ delivery events on the route to $i$\\
    $\mu_k$             & service rate of charger type $k\in\mathcal{K}$\\
    $T_{ij}$            & travel time from $i\in\mathcal{I}$ to $j\in\mathcal{J}$\\
    $T^\delta_{ij}$     & detour travel time when station $j\in\mathcal{J}_i$ is visited immediately before serving household $i\in\mathcal{I}$\\
    $C^\delta$          & cost of detour per time unit\\
    $C^\tau$            & cost of wait time and service time per time unit\\
    $C^\phi_j$          & fixed cost of charging station $j\in\mathcal{J}$ per time unit\\
    $C^\xi_k$           & cost of installing a type $k$ charger per time unit\\
    \bottomrule
  \end{tabularx}
\end{table*}

\begin{table*}[!ht]
  \caption{Variables and Cost Functions used in model $\mathbb{G}$. \label{tab:notations_2}}
  \centering\small
  \begin{tabularx}{\textwidth}{@{}l@{\quad}X@{}}
    \toprule
    \textbf{Variable}      & \textbf{Definition} \\
    \midrule
    $\rho_{jk}$ & utilization rate of charger type $k$ at station $j$, 
                   $\rho_{jk}\in\mathbb{R}_{\ge0},\;\rho_{jk}<1$, 
                   $ \rho_{jk}
                     = \frac{\sum_{i\in\mathcal{I}_j}\lambda_i\,x_{ijk}}{\mu_k\,s_{jk}}$\\
    $s_{jk}$    & number of type $k$ chargers allocated to station $j\in\mathcal{J}$, $s_{jk}\in\mathbb{Z}_{\ge0}$\\
    $W_{jk}$    & expected waiting (including charging) time for charger type $k\in\mathcal{K}$ at station $j\in\mathcal{J}$, 
                   $W_{jk}\in\mathbb{R}_{\ge0}$\\
    $x_{ijk}$   & $
                   \begin{cases}
                    1 & 
                    \shortstack[l]{if charger type $k\in\mathcal{K}$ at station $j\in\mathcal{J}_i$ is visited to charge \\
                    the vehicle preceding household $i\in\mathcal{I}$},\\
                     0 & \text{otherwise}
                   \end{cases}$\\
    $y_j$       & $
                   \begin{cases}
                     1 & \text{if station $j\in\mathcal{J}$ is active},\\
                     0 & \text{otherwise}
                   \end{cases}$\\
    \midrule
    \textbf{Cost Functions} & \textbf{Definition} \\
    \midrule
    $f_\phi(y)$   & station opening cost: $ f_\phi(y)=\sum_{j\in\mathcal{J}} C^\phi_j\,y_j$\\
    $f_\xi(s)$    & charger installation cost: $ f_\xi(s)=\sum_{j\in\mathcal{J}}\sum_{k\in\mathcal{K}} C^\xi_k\,s_{jk}$\\
    $f_\delta(x)$ & detour cost: $ f_\delta(x) =\sum_{(i,j,k)\in\mathcal{M}} \lambda_i\,C^\delta\,T^\delta_{ij}\,x_{ijk}$\\
    $f_\tau(x,W)$ & congestion (waiting) cost: $ f_\tau(x,W) =\sum_{(i,j,k)\in\mathcal{M}} \lambda_i\,C^\tau\,W_{jk}\,x_{ijk}$\\
    \bottomrule
  \end{tabularx}
\end{table*}

\subsection{Differences with Existing Models}
Compared with the model from \citet{aboolian2008location,aboolian2022efficient,berman2007multiple}, our proposed model $\mathbb{G}$ is more general since we do not constrain each household $i$ to be connected to its closest open station $j$ (i.e., the one with the shortest travel time). Instead, we introduce the sets $\mathcal{J}_i$ and $\mathcal{I}_j$ to represent the feasible stations for each household and the feasible households for each station, respectively. Recall that the set $\mathcal{J}_i$ contains the $k_c$ nearest stations to household $i$ based on travel times $T_{ij}$, while $\mathcal{I}_j$ includes the households for which station $j$ is among their $k_c$ nearest stations.
The studies by \citet{aboolian2008location,aboolian2022efficient} assume that (dropping the charger index $k$ without loss of generality) $x_{ij}=1$ for the household $i$ that is closest to a selected station $j$. With such an approach, the problem complexity is significantly reduced because the optimal solution of all variables can be derived once the open stations are known. In our proposed model $\mathbb G$, even if we were to fix the open stations (i.e., by fixing the values of the $y_j$ variables), we cannot infer the optimal solution to $x_{ijk}$ for that specific subset of $\mathcal{J}$. The reason is that the model selects the optimal values of the $x_{ijk}$ (i.e., allocation) decisions not only based on the travel time or the detour travel time between $i$ and $j$ but also considering the queue congestion and the particular charger $k$ used. Moreover, the use of $\mathcal{J}_i$ and $\mathcal{I}_j$ sets in model $\mathbb{G}$ limits the feasible options for each household and station, respectively, replacing the commonly adopted coverage constraints as explained earlier. For these reasons, the exact solutions from \citet{aboolian2008location,aboolian2022efficient} may not be feasible in our proposed model. 

\section{Exact Solution Method}
\label{sec:exact_solution_method}
In this section, we present an exact solution method for the SCLA model. The method consists of several components, including: a valid inequality (\Cref{{sec:valid_inequality}}), a binary representation of the number of chargers $s_{jk}$ (\Cref{{sec:binarizing_num_chargers}}), linearization of the quadratic objective function (\Cref{{sec:linearizing_the_objective}}), and most notably, a cutting plane method for handling the nonlinear waiting time constraints (\Cref{{sec:cutting_plane}}), two exact formulations that capture all of these features (\Cref{{sec:exact_formulations}}). We conclude this section by presenting computational challenges of the exact formulation for extreme-scale problems  (\Cref{{computational_challenges}}).

\subsection{A Valid Inequality}
\label{sec:valid_inequality}
Recall that for each station $j\in\mathcal{J}$, $\mathcal{I}_j$ is the set of households that can potentially be assigned to $j$. The household $i\in \mathcal{I}_j$ cannot be assigned to more than one charger type $k\in\mathcal{K}$ at station $j$. This constraint can be expressed by the inequality $ \sum_{k \in \mathcal{K}} x_{ijk} \leq 1 $ for all $i\in\mathcal{I}_j$. Imposing these inequalities for each station $j$ strengthens formulation $\mathbb G$ without cutting off any of its integer feasible solutions. 
The following theorem formally proves this claim.

\begin{theorem}
For any $j\in\mathcal{J}$, the inequality $ \sum_{k \in \mathcal{K}} x_{ijk} \leq 1 $ for all $i\in\mathcal{I}_j$ is valid for model $\mathbb G$.
\end{theorem}

\begin{proof}
Assume for contradiction that there exist $j\in\mathcal{J}$ and $i\in\mathcal{I}_j$ such that $ \sum_{k \in \mathcal{K}} x_{ijk} > 1 $. The set partitioning constraint, $ \sum_{j'\in\mathcal{J}, k\in\mathcal{K}} x_{ij'k}=1 $, ensures that household $i$ is assigned to exactly one station-charger pair in the entire network. Now, since $ \sum_{k \in \mathcal{K}} x_{ijk}>1 $ at station $j$, then $i$ must be assigned to multiple chargers at the same station, violating the uniqueness imposed by the set partitioning constraint. This contradiction implies that no feasible integral solution can violate the inequality. 
\end{proof}

\subsection{Binary Representation of the Number of Chargers}
\label{sec:binarizing_num_chargers}
To handle the discrete nature of charger allocation, we introduce a binary representation that allows us to precisely model the number of chargers at each station. 
To that end, define the set $\mathcal{S} = \{0, 1, 2, \ldots, \max_{j \in  \mathcal{J},\; k \in  \mathcal{K}} \overline{S}_{jk})\}$, which represents all possible discrete numbers of chargers that can be allocated at station $j$ for each charger type $k$.  
We then introduce the binary variable $z_{cjk}$ to indicate whether $c \in \mathcal S$ number of type $k$ chargers are allocated to station $j$. 
In the following, we present model $\mathbb{G}$ with the incorporation of the valid inequality from the previous section and the binarization of charger allocations.


\begin{align}
    \mathop{\mathrm{min}}_{x,y,W,s,z} \quad & f_{\phi}(y) + f_{\xi}(s) + f_{\delta}(x) + f_{\tau}(x,W) \notag \\
    \mathrm{s.\; t.} \quad
    & (\ref{eq:charge_from_only_open_stations})-(\ref{eq:non-neg3}) \notag \\
    &\sum_{k\in\mathcal{K}}x_{ijk} \leq 1 \quad \forall j \in \mathcal{J}, i \in \mathcal{I}_j \label{eq:exact_valid_ineq} \\
    & s_{jk} = \sum_{c\in\mathcal{S}} c \cdot z_{cjk} \quad \forall j \in \mathcal{J}, k\in \mathcal{K} \label{eq:exact_s_to_z} \\
    & \sum_{c\in\mathcal{S}} z_{cjk} \leq y_j \quad \forall j \in \mathcal{J}, k\in \mathcal{K} \label{eq:exact_z_link_y} \\
    & z_{cjk} \in \{0,1\} \quad  \forall c\in\mathcal{S}, j \in \mathcal{J},  k\in \mathcal{K}
    \label{eq:exact_z_cjk_domain}
\end{align}

The constraints (\ref{eq:exact_valid_ineq}) ensure that given a station $j$, each household is assigned to at most one charger type. Constraints (\ref{eq:exact_s_to_z}) link the number of chargers $s_{jk}$ to the binary variables $z_{cjk}$, ensuring that the correct number of chargers are allocated. Constraints (\ref{eq:exact_z_link_y}) ensure that chargers are  allocated only if the station is open. 

\subsection{Linearizing the Objective Function}
\label{sec:linearizing_the_objective}
The bilinear term $W_{jk}x_{ijk}$ in the objective function term, $f_{\tau}(x,W)$, 
makes the solution of the model computationally challenging. To address the bilinear terms, we explore the use of McCormick inequalities \citep{mccormick1976computability}. 
Consider linearizing $q_{ijk} = W_{jk}x_{ijk}$ over the domain $\mathcal{B} = \left[0, EW + 1/\mu_k\right] \times \{0, 1\}$. Recall that $W_{jk} \in \left[0, EW + 1/\mu_k\right]$ is the expected waiting time including charging time and $x_{ijk} \in \{0, 1\}$ is a binary variable. 
It is well known that the convex hull of the set $\{(W_{jk}, x_{ijk}, q_{ijk}) \in (\mathcal{B} \times \mathbb{R}) \; |\;  q_{ijk} = W_{jk}  x_{ijk}\}$ can be described using the following inequalities:

\begin{equation}
\begin{aligned}
& q_{ijk} \geq 0, \\
& q_{ijk} \geq W_{jk} - \left( EW + \frac{1}{\mu_k} \right)(1 - x_{ijk}), \\
& q_{ijk} \leq \left( EW + \frac{1}{\mu_k} \right) x_{ijk}, \\
& q_{ijk} \leq W_{jk}.
\label{exact_conv_so}
\end{aligned}
\end{equation}

This set of inequalities (\ref{exact_conv_so}), also known as the McCormick inequalities, defines the convex and concave envelopes of $W_{jk} x_{ijk}$ on the domain $\mathcal{B}$. 
Specifically, the first inequality ensures non-negativity of $q_{ijk}$. The second inequality ensures $q_{ijk} \geq W_{jk}$ when $x_{ijk} = 1$, whereas the third inequality ensures $q_{ijk} = 0$ when $x_{ijk} = 0$. The fourth inequality ensures $q_{ijk}$ never exceeds $ W_{jk}$. We note that these inequalities provide an exact linearization of the bilinear term, since $x_{ijk}$ is binary.
Therefore, we can equivalently rewrite the objective by replacing $W_{jk}x_{ijk}$ with $q_{ijk}$, resulting in: 

\begin{align}
    \mathop{\mathrm{min}}_{x,y,W,s,z,q} & f_{\phi}(y) + f_{\xi}(s) 
    + f_{\delta}(x)  + \sum_{(i,j,k)\in \mathcal{M}} \lambda_i C^{\tau} q_{ijk} \label{eq:objective_milp} \\
    \text{s.t.} \quad
    & (\ref{eq:charge_from_only_open_stations})-(\ref{exact_conv_so})  \notag 
\end{align}
where it is understood that the McCormick inequalities (\ref{exact_conv_so}) are written for all $(i,j,k)\in\mathcal{M}$.


\subsection{Cutting-Plane Approach for Nonlinear Constraints}\label{sec:cutting_plane}
We now address the nonlinearity arising from the expected waiting time function $\mathbb{W}$. We employ a cutting-plane approach that iteratively approximates $\mathbb{W}$ by introducing valid linear inequalities. The expected waiting time function, \(\mathbb{W}(\rho_{jk}, s_{jk})\), which models an \(M/M/s_{jk}\) queue, exhibits fundamental characteristics derived from the Erlang delay formula. Prior results \citep{grassmann1983convexity,lee1983note} establish that $\mathbb{W}$ is strictly increasing and strictly convex with respect to its first argument, the utilization ratio $\rho_{jk}$, whenever the stability condition $0<\rho_{jk}<1$ is satisfied.

\begin{property}\label{W_convex_in_rho}
For all $s_{jk} \in \mathbb Z_{>0}$, the function $\mathbb{W}(\cdot, s_{jk})$ is strictly increasing and strictly convex over the domain $(0, 1)$. 
\end{property}

We can exploit \Cref{W_convex_in_rho} to design a Kelley-type cutting plane method \citep{kelley1960cutting,cokyasar2023additive} to handle the nonlinear waiting time constraints. These constraints iteratively enforce lower bounds on $\mathbb{W}$ by adding supporting hyperplanes as cuts. To see this, temporarily fix $s_{jk}$ and observe that \Cref{W_convex_in_rho} also applies to the function $\mathbb{W}^\nu(\rho_{jk})$ defined as follows:

\begin{equation}
\label{W_nu}
    \mathbb{W}^\nu(\rho_{jk}) = \frac{\mathbb{P}(\rho_{jk}, s_{jk})}{1-\rho_{jk}}
\end{equation}

For any fixed candidate value $0 < \Tilde{\rho}_{jk} < 1$, define \[B_{jk}=\frac{\partial \mathbb{W}^\nu}{\partial \rho_{jk}}(\Tilde{\rho}_{jk}), \;\; A_{jk}=\mathbb{W}^\nu(\Tilde{\rho}_{jk}) - B_{jk} \Tilde{\rho}_{jk}.\] 
Then, \Cref{W_convex_in_rho} implies that the line, $A_{jk}+B_{jk}\rho_{jk}$, supports the graph of the function, $\mathbb{W}^\nu(\rho_{jk})$, at $\rho_{jk}=\Tilde{\rho}_{jk}$.
This implies $W_{jk} \geq (A_{jk}+B_{jk}\rho_{jk})/(\mu_k s_{jk}) + 1/\mu_k$.
However, note that this constraint is nonlinear in $s_{jk}$.
To address this nonlinearity, we introduce this constraint only conditionally on the current values of $\rho_{jk}$ and $s_{jk}$. Enforcing the constraints conditionally  will ensure a valid lower bound on $\mathbb{W}$ for all possible values of the decision variables, thus ensuring that they can be added as global cuts within a branch-and-bound search process. The conditional requirement can be imposed by exploiting the binary representation of $s_{jk}$ that we introduced previously. Also, since $\rho_{jk}$ is a function of only $s_{jk}$ and $x_{ijk}$, the conditional requirements can thus effectively be enforced using the values of only the binary variables $z_{cjk}$ and $x_{ijk}$, as shown below.

\begin{equation}
\label{W_cut_exact}
\begin{aligned}
    W_{jk} \geq &\ \frac{A_{jk}z_{c jk}}{\mu_k c} 
     + \frac{B_{jk} \sum_{i\in \mathcal{I}_j} \lambda_i (z_{c jk} + x_{ijk} - 1)}{\mu_k^2 c^2} + \frac{z_{c jk}}{\mu_k} \quad \forall c \in \mathcal{S} \setminus \{0\}, j \in\mathcal{J}, k \in\mathcal{K}
\end{aligned}
\end{equation}

The constraints in equation (\ref{W_cut_exact}) leverage the binary variable $z_{cjk}$ to establish the relationship between waiting time $W_{jk}$ and the number of allocated chargers. Recall that $z_{c jk}$ indicates whether $c$ chargers of type $k$ are assigned to station $j$. When $z_{c jk} = 1$, the constraint becomes active for the specific number of chargers $c$, enforcing a lower bound on the waiting time based on the allocated chargers. On the other hand, when $z_{cjk} = 0$, the constraint becomes slack, allowing the solver to effectively disregard that particular constraint.

\begin{remark}[Extension to the $M/G/s_{jk}$ Case]
The nonlinear waiting-time constraints in our formulation stem from the $M/M/s_{jk}$ queuing model, but the decomposition, linearization, and cutting-plane framework are largely queuing model-agnostic. 
Our same approach applies to alternative service-time distributions such as $M/G/s_{jk}$, requiring only a scaled cut instantiation based on the service-time coefficient of variation. For completeness, we add a detailed appendix presenting the relationship between the $M/M/s_{jk}$ and $M/G/s_{jk}$ systems and the corresponding cut instantiation.
\end{remark}

In actual implementation of the cutting planes method, we remove the constraint (\ref{eq:W_def}) from the original model and instead add the conditional constraints (\ref{W_cut_exact}) as lazy constraints every time an integer solution is obtained during the branch-and-bound process.
Notably, we do not need to add the constraint (\ref{W_cut_exact}) for all possible values of $c $, $j$ and $k$ simultaneously. Instead, as the solver progresses and comes up with integer solutions for a specific configuration of chargers for a given $j$ and $k$, we apply this cut only for the current $c$ at the current node. This approach ensures accurate bounding of $W_{jk}$ based on the actual number of allocated chargers. The binary nature of $z_{c jk}$ along with constraint (\ref{eq:exact_z_link_y}) guarantees that only one constraint is active for each combination of $j$ and $k$, corresponding to the chosen number of chargers in the current integer solution. Note that when $s_{jk}=0$, the constraint is undefined and never added. 

To measure convergence during the solution process, let $\mathbb{C}$ denote the objective value of the best-known integral feasible solution found so far, and let $\mathbb{\overline{C}}$ be the current best-known upper bound on the optimal objective. Define the solution gap as $Q := 1 - (\mathbb{C}/\mathbb{\overline{C}})$, which quantifies how close the current solution is to optimality.

The branch-and-bound procedure starts off by solving the SCLA problem formulation without constraint (\ref{eq:W_def}). 
Each time we encounter a feasible MIP node satisfying all other constraints including integrality constraints, we enter a solver callback, that  uses the solution vectors $\mathbf{s}$, $\mathbf{x}$, and $\mathbf{y}$ to calculate (i) $\Tilde{\rho}_{jk}= \sum_{i\in\mathcal{I}_j}\lambda_i x_{ijk}/(\mu_k s_{jk})$; (ii) $\overline{W}_{jk}$, the upper bound (value) of $W_{jk}$ via $\Tilde{\rho}_{jk}$; (iii) $\mathbb{\overline{C}}$, the upper bound (value) of $\mathbb{C}$ via $\overline{W}_{jk}$; and (iv) the solution gap $Q:=1-(\mathbb{C}/\mathbb{\overline{C}})$. While $Q>Q^\tau$, where $Q^\tau$ is an acceptable gap threshold, we use $s_{jk}$ as calculated from the binary variables and introduce (\ref{eq:W_def}) to $j$ and $k$ only when $s_{jk}>0$ and $W_{jk}<\overline{W}_{jk}$. Therefore, constraints (\ref{eq:W_def}) are  introduced only for specific charger types at stations where (i) at least one charger is allocated and (ii) the value of $W_{jk}$ is underestimated. When the condition $Q \leq Q^\tau$ is met or a predefined solution time limit is reached, the solver terminates and reports the best solution found.

Our previous work \citep{cokyasar2023additive} adopted the $M/G/k$ discipline and established the convexity of the expected waiting-time approximation with respect to utilization rate which is an essential property for maintaining global optimality under a cutting-plane approach. Extending the present framework to an $M/G/k$ formulation is conceptually straightforward, as also noted by the reviewer. Nevertheless, the $M/M/s$ queuing model used in this paper remains a standard and widely used approximation in the location science literature \citep{Aboolian_2009,aboolian2008location,castillo2009social,jung2014stochastic,kinay2023charging,xie2021integrated,xie2018long}. Regardless of the queuing discipline, we highlight in the manuscript that incorporating congestion effects at charging stations is a key element in the strategic charger allocation decision, which constitutes an important contribution of our paper.

To further clarify, we review Cosmetatos’ approximation \citep{cosmetatos1976some} for the $M/G/k$ queue in \citep[Eq.~(6.101)]{bolch2006queueing}, we observe that the waiting time for this principle is expressed as a weighted sum of the expected waiting times for $M/M/s$ and $M/D/s$ systems, that is, using Bolch’s notation:
\[
\overline{W}_{M/G/m} \approx c_B^2\,\overline{W}_{M/M/m} + (1-c_B^2)\,\overline{W}_{M/D/m}
\]
where $c_B^2$ denotes the squared coefficient of variation of the service time. Because large-scale empirical data on charging-time variability is still scarce (especially for commercial last-mile fleets), we assume $c_B^2 = 1$ corresponding to exponential service times. This aligns with the convention in the location-science and charging-infrastructure literature. 

\subsection{Summary of Exact Formulations}
\label{sec:exact_formulations}

After replacing constraint (\ref{eq:W_def}) with constraint (\ref{W_cut_exact}), we obtain two exact formulations for the SCLA problem. The mixed-integer quadratic program (MIQP) is denoted as $\text{E}_{\text{MIQP}}$ and shown below.

\begin{equation*}
\begin{aligned}
    \mathop{\mathrm{min}}_{x,y,W,s,z} \quad & f_{\phi}(y) + f_{\xi}(s) + f_{\delta}(x) + f_{\tau}(x,W) \\
    \mathrm{s.t.} \quad
    & (\ref{eq:charge_from_only_open_stations}) - (\ref{eq:rho}), \ (\ref{eq:W_bound}) - (\ref{eq:exact_z_cjk_domain}), \ (\ref{W_cut_exact}) 
\end{aligned}
\end{equation*}

The mixed-integer linear program (MILP) is denoted $\text{E}_{\text{MCC}}$ and shown below.
\begin{equation*}
\begin{aligned}
    \mathop{\mathrm{min}}_{x,y,W,s,z,q} \quad & (\ref{eq:objective_milp}) \\
    \mathrm{s.t.} \quad
    & (\ref{eq:charge_from_only_open_stations}) - (\ref{eq:rho}), \ (\ref{eq:W_bound}) - (\ref{exact_conv_so}), \ (\ref{W_cut_exact}) 
\end{aligned}
\end{equation*}

\subsection{Computational Challenges} 
\label{computational_challenges}
Our primary aim in this paper is to solve extremely large-scale instances of the SCLA problem, encompassing an immense number of households, candidate stations, and charger types. Even with a carefully formulated MIQP or MILP model, the dimensionality of the problem escalates rapidly. As $|\mathcal{I}|$, $|\mathcal{J}|$, and $|\mathcal{K}|$ grow large, the indexing set $\mathcal{M} := \{(i,j,k) \mid j \in \mathcal{J}, i \in \mathcal{I}_j, k \in \mathcal{K}\}$ can contain an enormous number of elements. Each such $(i,j,k)$ triplet requires a binary assignment variable $x_{ijk}\in\{0,1\}$, and every station-charger pair $(j,k)$ may demand multiple binary configuration variables $z_{cjk}\in\{0,1\}$ indexed by $c \in \mathcal{S}$. The result is a model with a massive number of binary variables, supplemented by integer charger allocation variables $s_{jk}$ and continuous waiting time variables $W_{jk}$, among others. When McCormick linearizations are applied in the $\text{E}_{\text{MCC}}$ formulation, each bilinear term $W_{jk}x_{ijk}$ is replaced by a continuous variable $q_{ijk}$ and multiple additional constraints. This further inflates the number of variables and constraints, pushing the limits of memory and computational time. Coupled with the complex nonlinear constraints related to the queuing-based expected waiting times, the resulting MIQP or MILP is extraordinarily challenging for commercial solvers. Constructing, loading, and solving such a gigantic formulation directly is impractical; and obtaining high-quality solutions within a reasonable timeframe is effectively impossible.
%
These computational challenges necessitate an alternative approach that exploits problem structure, reduces overhead, and can be parallelized. In the following sections, we propose a decomposition framework based on Lagrangian relaxation, designed to handle extreme-scale SCLA instances efficiently and reliably. 

\section{A Lagrangian-Based Decomposition Framework}
\label{sec:lagrengean_based_framework}
Our proposed solution approach for the SCLA problem integrates multiple components. We begin by applying Lagrangian relaxation (\Cref{sec:LR}) to introduce dual variables associated with key constraints. Building on this, we employ Lagrangian dual decomposition (\Cref{sec:LD}) to decompose the relaxed problem into station-level subproblems, each of which can be solved independently and in parallel. To solve each station-level subproblem, we employ a partial relaxation strategy combined with cutting-plane methods (\Cref{sec:semi_relaxation_cutting_plane}). We iteratively update the Lagrange multipliers via a subgradient method (\Cref{sec:sub_graident_method}). We then apply a three step rounding heuristic (\Cref{sec:rounding_heuristics}) to convert the fractional assignments from the subproblems into a feasible integral solution for the entire SCLA problem. Together, these techniques produce an integrated, scalable, and effective framework capable of handling large-scale SCLA problem instances. 

\subsection{Lagrangian Relaxation}
\label{sec:LR}
Observe that the set partitioning constraint (\ref{eq:one_assignment}) and the budget constraints (\ref{eq:budget_station})--(\ref{eq:station_limit}) are complicating constraints in formulation $\mathbb{G}$: the latter completely decomposes into smaller and independent subproblems, in the absence of these constraints. To exploit this structure, we introduce a Lagrange multiplier $\zeta_i$ for each equality constraint (\ref{eq:one_assignment}), corresponding to household $i$. To ensure non-negativity of these multipliers, we first modify the equality constraint (\ref{eq:one_assignment}) to a set covering constraint: $\sum_{j \in \mathcal{J}, k \in  \mathcal{K}} x_{ijk} \geq 1 \quad \forall i \in \mathcal{I}$. This is without loss of generality, since doing so  does not change the optimal solution of model $\mathbb G$. However, $\zeta_i$  can now only take non-negative values ($\zeta_i \geq 0$).  Similarly, we also introduce multipliers $\beta_{\phi}$, $\beta_{\xi}$, and $\nu$ corresponding to the station budget constraint (\ref{eq:budget_station}), the charger budget constraint (\ref{eq:budget_charger}), and the limit on the number of active stations (\ref{eq:station_limit}). The resulting Lagrangian relaxation, $L(\zeta, \beta_{\phi}, \beta_{\xi}, \nu)$, of model $\mathbb{G}$ is given by: 

\begin{align}
     \mathop{\mathrm{min}}_{x, y, W, s, z} \quad & f_{\phi}(y) + f_{\xi}(s) + f_{\delta}(x) + f_{\tau}(x,W) +  \sum_{i \in \mathcal{I}} \zeta_i \left(1 - \sum_{j \in \mathcal{J}, k \in  \mathcal{K}} x_{ijk}\right) \notag \\
    & + \beta_{\phi} \bigl( f_{\phi}(y) - B^\phi \bigr) 
      + \beta_{\xi} \bigl( f_{\xi}(s) - B^\xi \bigr) \notag \\ 
    & + \nu \left(  \sum_{j \in \mathcal{J}} y_j - \overline{Y}\right) \label{eq:lag_objective_function}  \\
    \mathrm{s. \ t.} \quad
    & (\ref{eq:charge_from_only_open_stations}), (\ref{eq:rho})-(\ref{eq:W_bound}), (\ref{eq:charger_limit})-(\ref{eq:exact_z_cjk_domain}). \notag
\end{align}


\subsection{Lagrangian Dual and Decomposition}
\label{sec:LD}
The Lagrangian dual $D$ of model $\mathbb{G}$, 
seeks the tightest lower bound on the optimal value of $\mathbb{G}$. Formally, if $L(\zeta,\beta_{\phi},\beta_{\xi},\nu)$ denotes the optimal objective value of the Lagrangian relaxation corresponding to fixed multipliers $(\zeta,\beta_{\phi},\beta_{\xi},\nu)$, 
then the dual problem $D$ can be defined as:
\[
D = \max_{\zeta,\beta_{\phi},\beta_{\xi},\nu\geq 0} L(\zeta,\beta_{\phi},\beta_{\xi},\nu)
\]
To demonstrate the decomposition arising in the inner minimization, we use the extensive form notation from \Cref{tab:notations_1} and \Cref{tab:notations_2}. 
In this form, the Lagrangian separates by station index $j$, leading to independent station-level subproblems $SP_j(\zeta,\beta_{\phi},\beta_{\xi},\nu)$. Formally, we have:
\begin{align}
    L(\zeta,\beta_{\phi}, \beta_{\xi}, \nu) = \quad & \sum_{j \in \mathcal{J}} SP_j(\zeta, \beta_{\phi}, \beta_{\xi}, \nu) + \sum_{i \in \mathcal{I}} \zeta_i - \beta_{\phi} B^\phi - \beta_{\xi} B^\xi - \nu \overline{Y}
    \label{eq:Lagrangean_dual_sub_problem} 
\end{align}

Because the Lagrangian function is separable by station index $j$, the minimization problem splits into independent subproblems $SP_j(\zeta,\beta_{\phi},\beta_{\xi},\nu)$, each involving only the variables and parameters associated with that station. This decomposition allows for the concurrent solution of all $SP_j$ and efficient handling of extremely large-scale instances. For each fixed station $j \in \mathcal{J} $, the subproblem ($SP_j$) can be stated as follows:

\begin{align}
    \min_{x, y, W, s, z} &  \big[(1 + \beta_\phi)C^{\phi}_j + \nu \big] y_j +  \left. \sum_{k \in \mathcal{K}} \big[(1 + \beta_\xi)C^{\xi}_k \big] s_{jk} \notag  \right.\\ + & \sum_{i \in \mathcal{I}} \sum_{k \in \mathcal{K}} (C^\delta T^{\delta}_{ij} \lambda_i + C^\tau W_{jk} \lambda_i - \zeta_i) x_{ijk}  
    \label{eq:Lagrangean_sub-problem_objective}\\
    \mathrm{s.\ t.} \quad
    & x_{ijk} \leq y_j \quad \forall i \in \mathcal{I}_j, k \in  \mathcal{K}\label{eq:lag_sub-problem_assign_only_if_open} \\
    &\sum_{k\in\mathcal{K}}x_{ijk} \leq 1 \quad \forall i \in \mathcal{I}_j \label{eq:lag__valid_ineq} \\
    & \mu_k s_{jk} (1 - \epsilon) \geq \sum_{i \in \mathcal{I}_j} \lambda_i x_{ijk} \quad \forall k \in  \mathcal{K}\label{eq:lag_sub-problem_rho}\\
    & \mathbb{W} \left( \frac{\sum_{i\in \mathcal{I}_j} \lambda_{i} x_{ijk}}{s_{jk} \mu_k}, s_{jk} \right) \leq W_{jk} \quad \forall k \in  \mathcal{K}  \label{eq:lag_sub-problem_W_def} \\
    & W_{jk} \leq EW + \frac{1}{\mu_k} \quad \forall k \in  \mathcal{K} \label{eq:lag_sub-problem_W_def_2} \\
    & s_{jk} \leq \overline{S}_{jk} \quad \forall k \in  \mathcal{K}\label{eq:lag_sub-problem_number_of_chargers} \\
    & s_{j,k} \leq \overline{S}_{jk} \sum_{i \in \mathcal{I}_j} x_{ijk} \quad \forall k \in \mathcal{K}\label{eq:lag_charger_allocation} \\
    & y_j \leq \sum_{k \in \mathcal{K}} s_{jk}\label{eq:lag_station_open_charger_allocation}\\
    & s_{jk} = \sum_{c\in\mathcal{S}} c \cdot z_{cjk} \quad k\in \mathcal{K} \label{eq:s_to_z} \\
    & \sum_{c\in\mathcal{S}} z_{cjk} \leq y_j \quad \forall k\in \mathcal{K} \label{eq:z_link_y} \\
    & z_{cjk} \in \{0,1\} \quad  \forall c\in\mathcal{S}, k\in \mathcal{K} \label{eq:z_cjk_domain}\\
    & x_{ijk}, y_j \in \{0, 1\} \quad \forall i \in \mathcal{I}_j, k \in  \mathcal{K} \label{eq:lag_sub-problem_nonneg1}\\
    & s_{jk} \in \mathbb{Z}_{\geq 0}, \quad W_{jk} \in \mathbb{R}_{\geq 0} \qquad \forall k \in  \mathcal{K}.\label{eq:lag_sub-problem_nonneg2}
\end{align}
Each station-level subproblem $SP_j$ focuses on station opening, charger installation, detour costs, and waiting times for station $j$ only, independent of all other stations. 
In the next section we present two solution approaches for solving $SP_j$.

\subsection{Partial Linear Relaxation with Cutting Planes}
\label{sec:semi_relaxation_cutting_plane}

Observe that subproblem $SP_j$ is a simpler version of the full-scale problem, restricted to station $j$ and without any complicating constraints that couple these stations.
Therefore, it is possible directly apply the solution method from \Cref{sec:exact_solution_method} to solve $SP_j$.
However, computational experiments on extremely large-scale instances demonstrate significant hurdles. Even when focusing on a single station $j$, the cardinality $|\mathcal{I}_j|$ of the household set and the associated binary decision variables $x_{ijk}$ can be exceptionally large. As a result, the decomposed subproblems also become computationally intractable, exceeding practical time limits.

To address this challenge, we present modified variants of the exact methods, $\text{E}_{\text{MIQP}}$ and $\text{E}_{\text{MCC}}$, for the Lagrangian subproblems $SP_j$. These variants, denoted as $\text{L}_{\text{MIQP}}$ and $\text{L}_{\text{MCC}}$, are similar to their exact counterparts. The only modification is the relaxation of integrality conditions on the binary variables $x_{ijk}$. As a consequence, they only constitute relaxations of $SP_j$. Although these relaxations simplify the subproblem and reduce computational complexity, it is important to note that they continue to provide rigorous lower bounds on the optimal value of the original SCLA problem. However, at the same time, they may also yield fractional solutions that are not integer-valued. To regain integer feasibility, we develop specialized rounding heuristics that transform the fractional solutions into integral ones.


\subsubsection{$\text{L}_{\text{MIQP}}$ Formulation for Solving $SP_j$.}
We relax $x_{ijk} \in [0,1]$. Note that the other integer variables, namely $y_j$, $z_{cjk}$, and $s_{jk}$ retain their original discrete domains. The objective and constraints remain identical to the problem $SP_j$ defined in equations (\ref{eq:Lagrangean_sub-problem_objective})--(\ref{eq:lag_sub-problem_nonneg2}) except we replace $x_{ijk}\in\{0,1\}$ with $0 \leq x_{ijk} \leq 1$. The resulting formulation can be summarized as follows.

\begin{align}
    \min_{x,y,W,s,z} \;&(\ref{eq:Lagrangean_sub-problem_objective}) \notag \\[6pt]
    \text{s.t.}\quad
    &(\ref{eq:lag_sub-problem_assign_only_if_open})-(\ref{eq:z_cjk_domain}), (\ref{eq:lag_sub-problem_nonneg2})\notag\\
    & 0 \leq x_{ijk} \leq 1,\; y_j \in \{0,1\} .\label{eq:semi_relaxed_nonneg1}
\end{align}

\subsubsection{$\text{L}_{\text{MCC}}$ Formulation for Solving $SP_j$.}

Similar to $\text{L}_{\text{MIQP}}$, the $\text{L}_{\text{MCC}}$ formulation relaxes $x_{ijk} \in [0, 1]$ in the original $\text{E}_{\text{MCC}}$ model, whereas $y_j$, $z_{cjk}$, and $s_{jk}$ continue to remain discrete. Similarly, the bilinear term $W_{jk}x_{ijk}$ is linearized via an auxiliary variable $q_{ijk}$ and appropriate McCormick-type inequalities. The objective and constraints reflect those of the problem $SP_j$ but with $x_{ijk}\in[0,1]$ and $q_{ijk}$ introduced to represent the convex hull of the product $W_{jk}x_{ijk}$. 

\begin{align}
    \min_{x, y, W, s, z, q} & \big[(1 + \beta_\phi)C^{\phi}_j + \nu \big] y_j +  \sum_{k \in \mathcal{K}} \big[(1 + \beta_\xi)C^{\xi}_k \big] s_{jk} \notag \\
    +& \sum_{i \in \mathcal{I}} \sum_{k \in \mathcal{K}} ((C^\delta T^{\delta}_{ij} \lambda_i - \zeta_i) x_{ijk} + C^\tau \lambda_{i} q_{ijk})
    \label{eq:sp_j_bilinear_objective_mc}\\[6pt]
    \text{s.t.}\quad
    & (\ref{exact_conv_so}), (\ref{eq:lag_sub-problem_assign_only_if_open})-(\ref{eq:z_cjk_domain}), (\ref{eq:lag_sub-problem_nonneg2}), (\ref{eq:semi_relaxed_nonneg1}) \notag 
\end{align}


In both \(\text{L}_{\text{MIQP}}\) and \(\text{L}_{\text{MCC}}\) subproblem formulations, the nonlinear waiting-time constraint \eqref{eq:lag_sub-problem_W_def} is replaced by the  cutting-plane constraint \eqref{W_cut_exact}, applied for each fixed station \(j\).


\subsection{Subgradient Method for Updating Lagrangian Multipliers}
\label{sec:sub_graident_method}
The subproblems $SP_j$ are solved for fixed values of the Lagrange multipliers $\zeta, \beta_{\phi}, \beta_{\xi}, \nu$. Aggregating the subproblem solutions allows us to calculate the Lagrangian objective value, $L(\zeta, \beta_{\phi}, \beta_{\xi}, \nu)$, shown in (\ref{eq:Lagrangean_dual_sub_problem}). The Lagrangian dual, which provides a lower bound to the original problem $\mathbb{G}$, can further be improved by updating the Lagrange multipliers. To that end, we employ the classical subgradient method \citep{boyd2003subgradient,barahona2005near,barahona2000volume} to update the multipliers, which consist of the following steps.

\textbf{Step 0:} Let $\eta = (\zeta, \beta_{\phi}, \beta_{\xi}, \nu)$ denote the tuple of Lagrange multipliers. Initialize $\eta^0 = 0$. Solve subproblem $SP_j$ for all station locations $j \in \mathcal{J}$, independently and in parallel, to obtain subproblem-optimal solutions $x_{ijk}^0, y_j^0, s_{jk}^0$, $W_{jk}^0$ for each $SP_j$. Aggregate these solutions and evaluate the Lagrangian \(\mathcal{L}(\eta^0)\). Run a rounding heuristic (Section \ref{sec:rounding_heuristics}) on the aggregated solution to obtain a feasible solution with objective value \(z^0\), and set the initial upper bound \(UB = z^0\). Set $t=1$, $\lambda_\eta^1 = 0.1$, and $\eta' = \eta^0$.

\textbf{Step 1:} Compute the subgradients for all Lagrange multipliers at iteration $t$ as follows:

\begin{align}
    & g_{\zeta_i}^t = 1 - \sum_{j \in \mathcal{J},\, k \in \mathcal{K}} x_{ijk}^{t-1} \qquad \forall i \in \mathcal{I} \\
    & g_{\beta_\phi}^t = \sum_{j \in \mathcal{J}} C^\phi_j y_j^{t-1} - B^\phi \\
    & g_{\beta_\xi}^t = \sum_{j \in \mathcal{J},\, k \in \mathcal{K}} C^\xi_k s_{jk}^{t-1} - B^\xi \\
    & g_{\nu}^t = \sum_{j \in \mathcal{J}} y_j^{t-1} - \overline{Y}.
\end{align}

Compute step sizes $s_{\eta}^t$ for each component $\eta \in \left\{ \{\zeta_i\}_{i \in \mathcal I}, \beta_{\phi}, \beta_{\xi}, \nu \right\}$ as follows:
\[
s_{\eta}^t = \lambda_{\eta}^t \frac{\Delta}{||g_{\eta}^t||^2}, \text{ where } \Delta = UB - \mathcal{L}(\eta'),
\]
$UB$ is any valid upper bound on the optimal objective value of the original formulation $\mathbb G$, $\eta'$ represents the best multipliers found so far, and 
$\lambda_{\eta}^t \in (0,2)$ is a step length parameter for each $\eta$.
Using these step sizes, we update the multipliers as follows:
\[
\eta^t = \max(0, \eta^{t-1} + s_{\eta}^t g_{\eta}^t).
\]
This process ensures that each Lagrange multiplier is updated independently, taking into account its specific subgradient and step size.

\textbf{Step 2:} Evaluate $\mathcal{L}(\eta^t)$ by solving $SP_j$ for all station locations $j \in \mathcal{J}$, independently and in parallel, 
to obtain subproblem-optimal solutions $x_{ijk}^t, y_j^t, s_{jk}^t$, $W_{jk}^t$ for each $SP_j$. Run a rounding heuristic (Section \ref{sec:rounding_heuristics}) on the new aggregated solution to obtain a feasible solution with objective value \(z^t\). If \(z^t < UB\), update the upper bound by setting $UB \gets z^t.$

\textbf{Step 3:} For each $\eta$, update $\lambda_\eta^t$ as follows.
\begin{itemize}
    \item If \(\mathcal{L}(\eta^t) \leq \mathcal{L}(\eta')\), then no improvement has been achieved: label the iteration as \textit{red}, indicating that the current step size for $\eta$ is likely too large. Reduce \(\lambda_{\eta}^t = \frac{9}{10}\lambda_{\eta}^{t-1}\) to encourage progress in subsequent iterations.

    \item If \(\mathcal{L}(\eta^t) > \mathcal{L}(\eta')\), perform a directional test by forming direction vectors from current \((g_{\eta}^t)\) and aggregated \((\bar{g}_{\eta}^t)\) subgradients, i.e.\ let \(v_{\eta}^t = g_{\eta}^t - \bar{g}_{\eta}^t\). Then compute the inner product \(d_{\eta} = \langle v_{\eta}^t,\, g_{\eta}^t\rangle\):
    \begin{itemize}
        \item If \(d_{\eta} \geq 0\), label the iteration as \textit{green}, suggesting that even larger steps, \(\lambda_{\eta}^t = 1.1 \lambda_{\eta}^{t-1}\), could yield further improvement.
        \item If \(d_{\eta} < 0\), label the iteration as \textit{yellow}, indicating that the current step size is roughly appropriate, and we keep \(\lambda_{\eta}^t = \lambda_{\eta}^{t-1}\) unchanged.
    \end{itemize}

\end{itemize}
Since each \(\lambda_{\eta}^t\) is updated independently based on the multiplier’s own classification, this process systematically refines step sizes, guiding the subgradient method toward improved dual bounds and enhanced solution quality over time.

\textbf{Step 4:} If $\mathcal{L}(\eta^t) > \mathcal{L}(\eta')$, update $\eta' \gets \eta^t$. Go to \textbf{Step 1}.

\subsection{Rounding Heuristic}
\label{sec:rounding_heuristics}
The Lagrangian-based decomposition provides lower bounds on the original problem \(\mathbb{G}\), but may not produce primal feasible solutions. To address this issue, we develop a rounding heuristic that uses the station-level solutions from the Lagrangian subproblems \(SP_j\) to construct a feasible primal solution obtain an upper bound on \(\mathbb{G}\). We emphasize that the heuristic is not an independent procedure; rather, it is explicitly driven by and dependent on the subproblem solutions \((x_{ijk}^t, y_j^t, s_{jk}^t, W_{jk}^t)\) obtained after solving all subproblems \(SP_j\) for every \(j \in \mathcal{J}\), in Step~2 of the procedure described in the previous subsection.

The rounding procedure (\Cref{alg:primal_feasible}) consists of three main steps, each realized through a set of algorithms that refine the partially relaxed solution into a fully integral one. The procedure begins with the station opening decisions \(y_j^t\) and fractional household assignments \(x_{ijk}^t\) at iteration \(t\), and produces a final solution \((x', y', s', W')\) that satisfies all integrality and feasibility constraints.

\subsubsection*{Step 1: Adaptive Station Selection for Household Coverage}
At the outset, we start with a candidate set of open stations \(\{ j \in \mathcal{J} \mid y_j^t = 1 \}\) derived from the Lagrangian solutions. This set may initially be empty (if no station was selected by the subproblem solutions), partially cover the households (if some are covered and some are not), or already fully cover all households. Our adaptive station selection heuristic (\Cref{alg:adaptive_station_opening}) is designed to handle all three scenarios while respecting the natural constraint that no more than \(|\mathcal{J}|\) stations can be opened in total.

At any given iteration, let \(\mathcal{U}^c \subseteq \mathcal{I}\) denote the set of uncovered households. For each unopened station \(j \in \mathcal{J}\), we compute a weight \(w_j = |\mathcal{I}_j \cap \mathcal{U}^c|\), representing how many uncovered households could be covered by station \(j\). Stations that cannot cover any uncovered households are ignored. Among the candidate stations, the algorithm selects one to open based on probabilities proportional to the weights \(w_j\). This probabilistic selection approach reduces the risk of biases and prevents repeatedly opening the same or suboptimal stations.

By construction, \Cref{alg:adaptive_station_opening} must terminate with \(\mathcal{U}^c = \emptyset\). Specifically, if at some point no candidate stations remained to cover \(\mathcal{U}^c\), it would imply that each uncovered household had no accessible station in its predefined neighborhood \(\mathcal{J}_i\). However, \Cref{alg:construct_Ji_Ij} ensures \(\mathcal{J}_i \neq \emptyset\) for every household \(i\), making such a scenario impossible. Therefore, with each station opening strictly reducing \(|\mathcal{U}^c|\) by at least one, and given that we cannot open more than \(|\mathcal{J}|\) stations, we must achieve full coverage of all households within at most \(|\mathcal{I}|\) iterations. Hence, the algorithm ends with \(\mathcal{U}^c = \emptyset\), ensuring that every household is  covered by at least some open station. While this coverage alone does not guarantee feasible assignments under all constraints (e.g., charger capacities and waiting times), it establishes a fully covered station configuration for subsequent assignments.

\subsubsection*{Step 2a: Deterministic Household Assignment to Stations and Chargers}
In Step 2a our goal is to transform the fractional assignment results \(x_{ijk}^t\) from the subproblem solutions \(SP_j\) into integral assignments \(x_{ijk}' \in \{0,1\}\) that  satisfy all the constraints. \Cref{alg:household_iteration} manages this task. In algorithm, the input flag \texttt{probabilistic} is set to \texttt{false}, that lets the algorithm to enforce a deterministic rounding scheme in which each household is assigned to exactly one feasible station–charger pair.

We prioritize households by sorting them in descending order of their demand rates \(\lambda_i\), thus addressing the most demanding households first. For each household \(i\), we consider those stations \(j\) that are open (\(y_j'=1\)) and accessible to \(i\) (\(j \in \mathcal{J}_i\)). We then generate a candidate list \(\mathcal{C}_i\) of all station-charger pairs \((j,k)\) that are both open and capable of serving household \(i\). To select a suitable \((j,k)\) from \(\mathcal{C}_i\), we impose a heuristic sorting that considers the following:  
\begin{enumerate}
    \item The fractional assignments \(x_{ijk}^t\) from the subproblem solutions, leveraging prior knowledge about the likelihood of a successful integral assignment
    \item Charger usage counts \(n_{jk}\) and installation costs \(C^\xi_k\), guiding us toward cost-effective and less congested chargers
    \item Current waiting times \(W_{jk}\), ensuring that we assign households to station-charger pairs that can meet the waiting time constraints
\end{enumerate}

This ordered candidate list is processed by the \textsc{CheckAndAssign} subroutine (\Cref{alg:check_and_assign}), which attempts to fix \(x_{ijk}'=1\) while adjusting the number of chargers \(s_{jk}'\) and updating \(W_{jk}'\) to maintain feasibility. If an appropriate \((j,k)\) is found, household \(i\) is successfully assigned. Otherwise, if no candidate pair meets all constraints, the household is added to the set \(\mathcal{U}^a\) of unassigned households for further handling in subsequent steps. 

\subsubsection*{Step 2b: Probabilistic Household Assignment to Stations and Chargers}
Alternatively, when the input flag \texttt{probabilistic} is set to \texttt{true} in \Cref{alg:household_iteration}, Step 2b becomes a probabilistic assignment mode that constructs a fractional household station–charger allocation consistent with the relaxed Lagrangian solution. 
Here, the fractional variables \(x_{ijk}^t \in [0,1]\) are interpreted as probabilities of assigning household \(i\) to charger type \(k\) at station \(j\). 
For each household \(i\), we define the feasible set 
\(\mathcal{C}_i = \{(j,k) : j \in \mathcal{J}_i,\, y_j' = 1,\, x_{ijk}^t > 0\}\) 
and normalize:
\[
p_{ijk} = \frac{x_{ijk}^t}{\sum_{(j',k') \in \mathcal{C}_i} x_{ij'k'}^t},
\qquad (j,k)\in\mathcal{C}_i
\]
Each \(p_{ijk}\) represents the probability that household \(i\) is served by charger type \(k\) at station \(j\). The corresponding effective arrival rate into each station–charger pair becomes 
\(\bar{\lambda}_{jk} = \sum_{i \in \mathcal{I}_j} \lambda_i p_{ijk}\), 
ensuring that all queuing and capacity constraints defined by \(\mathbb{W}(\cdot)\) remain valid.

\subsubsection*{Step 3: Adjusting for Overload and Ensuring Final Feasibility}

If after attempting the assignment, some households remain unassigned (\(\mathcal{U}^a \neq \emptyset\)), we attempt to open additional stations to reduce overload and accommodate these households. \Cref{alg:additional_stations_overload} chooses additional stations, if possible, to support the uncovered demand. After updating the station configuration, we rerun the household assignment iteration until no further improvements are possible. This iterative refinement ensures that we eventually reach a stable configuration where all households are assigned.

Once every household \(i\) is assigned, we verify that the number of chargers \(s_{jk}'\) and expected waiting times \(W_{jk}'\) are consistent with the constraints. The computations for \(s_{jk}'\) and \(W_{jk}'\) follow directly from the final assignments and the queuing system constraints integrated into \textsc{CheckAndAssign}. We produce a feasible solution \((x_{ijk}', y_j', s_{jk}', W_{jk}')\) for \(\mathbb{G}\).

\subsection{Combining Cutting-Planes Approach, Subgradient Method, and Rounding Heuristic}

At any iteration \( t \), given Lagrangian multipliers \(\eta^t = (\zeta^t, \beta_{\phi}^t, \beta_{\xi}^t, \nu^t)\), we solve all subproblems \( SP_j \) for each \( j \in \mathcal{J} \) in parallel in Step~2. This produces solutions \((x_{ijk}^t, y_j^t, s_{jk}^t, W_{jk}^t)\) and a Lagrangian value \(\mathcal{L}(\eta^t)\), providing a lower bound on \(\mathbb{G}\). We then apply a rounding heuristic to these fractional solutions, yielding an integral feasible solution with objective value \(UB^t\). Depending on the value of the input flag \texttt{probabilistic}, the rounding phase may generate either a deterministic assignment policy (Step~2a, with $x_{ijk}'\!\in\!\{0,1\}$) or a probabilistic one (Step~2b, with fractional $x_{ijk}'\!=\!p_{ijk}$ interpreted as long-run allocation probabilities). This allows the same decomposition framework to output a fully integral plan or a probabilistic sharing policy consistent with the relaxed solution.
The primal-dual gap at iteration \(t\) is $\Delta^t = 1 - \frac{\mathcal{L}(\eta^t)}{UB^t}$. If \(\Delta^t > \delta\), where \(\delta > 0\) is a specified convergence criterion, we update the Lagrangian multipliers \(\eta^t\) using the subgradient method and, if necessary, refine the upper bound \(UB^t\). We then re-solve the subproblems \( SP_j \) in parallel with these updated multipliers for potentially reducing the gap. This iterative process continues until \(\Delta^t \leq \delta\). In practice, convergence may be determined either by reaching the specified tolerance \(\delta\) or by imposing computational limits, such as a maximum number of iterations or a time budget. By adapting the multipliers at each iteration and repeatedly solving the subproblems with refined multipliers, we ensure that our method eventually produces a high-quality solution, making it suitable for extreme-scale SCLA instances. 
We provide additional details about our parallel computing implementation in \Cref{sec:large_scale_implementation}.

\section{Computational Experiments}\label{sec:comp_exp}

All experiments were conducted on Intel~Xeon~CPU~6248R. We use Gurobi~10.0.2 \citep{Gurobi2024} as the MIQP and MIP solver, and we use mpi4py \citep{dalcin2008mpi} for parallel processing.

\subsection{Design of Experiments}
\label{sec:DOE}
To test our proposed frameworks' ability to solve extreme-scale SCLA problem instances, we conduct computational experiments using data from the Chicago metropolitan area. We utilize e-commerce daily demand information and road network data, both sourced from POLARIS (Planning and Operations Language for Agent-based Regional Integrated Simulation).  POLARIS is an advanced simulation tool developed by researchers at Argonne National Laboratory \citep{auld2016polaris}. 

We defined six scenarios based on two key factors: charging locations and policy frameworks. For charging locations, planners can install chargers at  traffic analysis zone (TAZ) sites,  existing depots, or  a combination of both. Regarding policy frameworks, scenarios use either single-agency or multiagency operations. Under single-agency operation, each household belongs to exactly one agency and may use only  that agency’s stations (depot or TAZ). Under multiagency collaboration, however, all depots are pooled across agencies, and TAZ-based stations carry no agency affiliation; thus every household can access any open station, whether it is at a depot or a TAZ.

\begin{figure}
    \centering
    \begin{subfigure}[t]{0.48\textwidth}
        \centering
        \includegraphics[width=\textwidth]{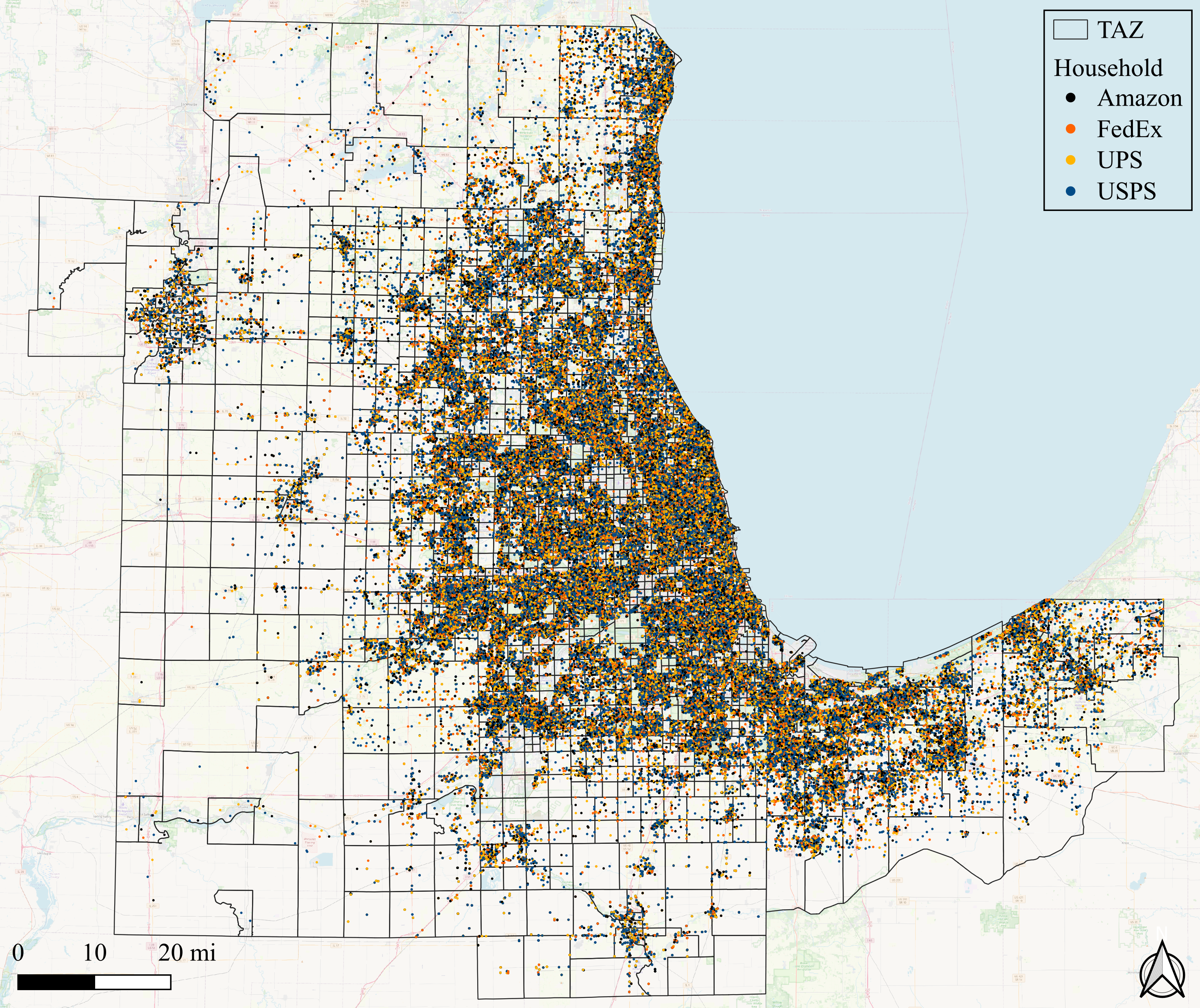}
        \caption{Households.}\label{fig:household_distribution}
    \end{subfigure}
    \quad
    \begin{subfigure}[t]{0.48\textwidth}
        \centering
        \includegraphics[width=\textwidth]{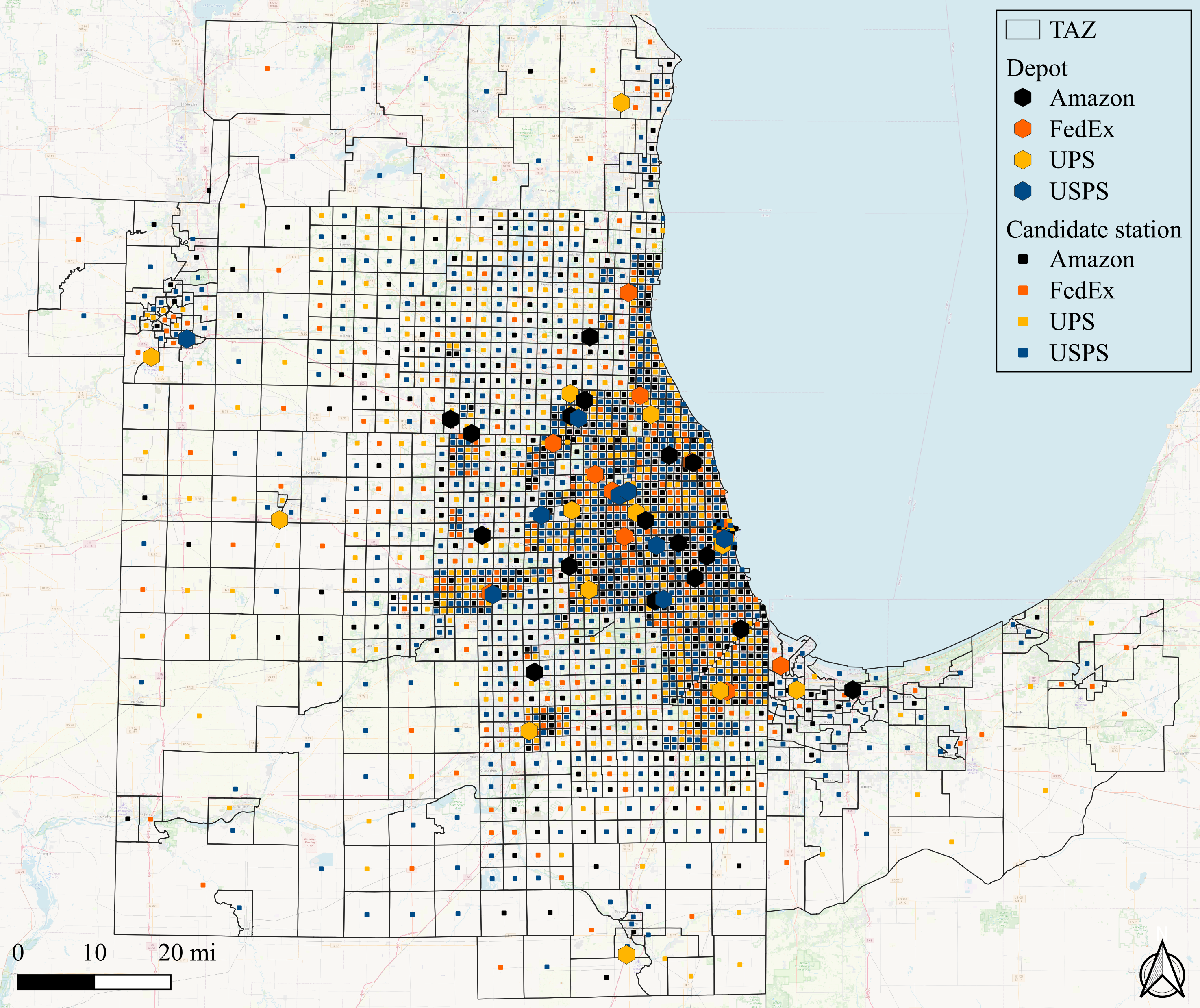}
        \caption{Candidate charging stations and depots.}\label{fig:depot_candidate_station_distribution}
    \end{subfigure}

    \caption{Illustration of the spatial distribution of households, potential charging stations, and existing charging depots across the Chicago metropolitan area.}
    \label{fig:spatial}
\end{figure}

\Cref{fig:spatial} shows the spatial distribution of households and candidate charging stations across the Chicago metropolitan area. The map \Cref{fig:spatial} (a) highlights potential charging stations (squares) located at centroids of TAZs. Existing depots for four major agencies---Amazon, FedEx, UPS, and USPS---are depicted with larger symbols of distinct shapes and colors (e.g., black circles for Amazon, orange hexagons for FedEx). These depots can also be used as candidate station sites, although their station-opening cost \(C^\phi\) is zero. 

The map \Cref{fig:spatial} (b) overlays the same TAZ grid but now includes the potential station locations via depots and TAZs, again color-coded by agency affiliation. Each agency represents a single subscription policy in one scenario (restricting households to use that agency’s stations) or participates in a multiagency collaboration in other scenarios. The dataset consists of 449,367 households, 1,958 potential station sites (TAZ centroids), and 53 existing depots, making it the largest known instance for a queue-based charger location-allocation study. 

The problem parameters used are based on the literature \citep{RivianEDV2024,UPSDriverSalary2024,davatgari2021location,ElectrifyAmerica2024,LightningEMotors2022,SmithCastellano2015,Williams2020}. To compute the detour and waiting costs per minute, we utilized financial and operational data from FedEx \citep{FedExRevenue2024,FedExOperationalDays2024,FedExFleetSize2024,FedExDailyPackagesDelivered2024} and UPS \citep{UPSDriverSalary2024}. 
The detour cost was derived by looking at FedEx's revenue metrics. We calculate the revenue per vehicle per day by distributing FedEx's annual revenue over its operating days and fleet size, assuming that any detour could potentially cut into this daily revenue, thereby translating into a cost. For the waiting cost, we based our calculation on the hourly wage rate of UPS drivers, treating time spent waiting as a direct labor cost. By estimating the cost per minute from the hourly wage, we capture the expenses incurred during periods when drivers are idle. Alongside these costs, $\mu_k$ quantifies the number of full charging cycles that each type of charger can complete in one minute. It is calculated by dividing the total charging capacity needed to fully charge a vehicle from a specific initial SOC by the charger’s power output. The result is then converted from the total time required to charge (in seconds) to how many cycles can be completed per minute; the parameters are summarized in \Cref{tab:general_parameters}. 

For simplicity, we assume each vehicle completes exactly one charging stop per day, making \(1/\nu\) the fraction of deliveries that need a charge, where \(\nu\) is the average number of deliveries per vehicle per day derived from historical data. Hence, for household \(i\), \(N_i^C = N_i^D / \nu\), leading to \(\pi_i = N_i^C / N_i^D = 1/\nu\). Finally, \(\lambda_i = \gamma_i\,\pi_i\) is the overall charging rate for that household. 

Our analysis incorporates three charger types of varying power outputs along the same lines in the literature \citep{liu2017locating,yilmaz2012review}). In order to standardize the objective function, all costs are converted to USD per day. This conversion assumes lifespans of 10 years for chargers and 40 years for facilities, based on estimates by \citet{bennett2022estimating}. \Cref{tab:charger_specific_parameters} presents a summary of the charger-specific parameters.

\begin{table*}
\centering
\small
\caption{General parametric values}
\label{tab:general_parameters}
\begin{tabular}{llllllll}
\hline
$C^\phi$ & $C^\delta$ & $C^\tau$ & $B^\phi$ & $B^\xi$ & $\overline{Y}$ & $EW$ & Lifetime \\
(USD) & (USD/min) & (USD/min) & (USD/day) & (USD/day) & (units) & (min) & (years)\\
\hline
1,000,000 & 1.09 & 0.70 & 150,000 & 500,000 & 1000 & 30 & 40 \\
\hline
\end{tabular}
\end{table*}

\begin{table*}
\centering
\small
\caption{Charger-specific parametric values}
\label{tab:charger_specific_parameters}
\begin{tabular}{llllll}
\hline
Type & $C^\xi$ (USD) & Power (kWh) & $\mu_k$ (charging/hour) & $\overline{S}_{jk}$ (units) & Lifetime (years) \\
\hline
Basic    & 73,000  & 50  & 0.53 & 20 & 10 \\
Moderate & 157,000 & 180 & 1.90 & 20 & 10 \\
Fast     & 228,000 & 360 & 3.81 & 20 & 10 \\
\hline
\end{tabular}
\end{table*}

\subsection{Computational Performance of Exact Cutting-Planes Method and Lagrangian Dual Decomposition Approach}
\label{sec:exact_vs_lagrangean}

To effectively solve the SCLA problem, we first develop an understanding of our  solution methods by analyzing the computational performance of our exact cutting planes approach and comparing it with the Lagrangian dual decomposition method. This analysis utilizes cardinality of $\mathcal{J}$, $k_c$, and $\mathcal{I}$ as key determinants of problem size. A testbed of instances was generated by utilizing e-commerce daily demand information and road network data of the Chicago metropolitan area serving e-commerce deliveries from POLARIS (\cite{auld2016polaris}). As a baseline, the parametric design provided in \Cref{sec:DOE} was utilized. 

For the computational study in this section we focus on the central section of Cook County around the Chicago Loop area, which includes \(\lvert\mathcal{J}\rvert = 20\) potential charging locations. We assume a multiagency policy, meaning any open station whether at a TAZ or an existing depot can serve any household. Within this region, we sample \(\lvert \mathcal{I} \rvert\) households. Details of each scenario appear in \Cref{tab:exact_vs_lagrangean}.

We compare four different approaches to solve five instances of a specific scenario with a time limit of 1 hour. The first method, denoted as $\text{E}_\text{MIQP}$, directly solves the exact MIQP formulation. It provides  valid upper and lower bounds denoted $\text{UB}^{\text{E}_\text{MIQP}}$ and $\text{LB}^{\text{E}_\text{MIQP}}$, respectively. 
The second method, denoted as $\text{E}_\text{MCC}$, relaxes the MIQP objective by replacing bilinear $W_{jk} x_{ijk}$ terms with McCormick envelopes. This approach yields a valid lower bound $\text{LB}^{\text{E}_\text{MCC}}$ but requires recalculating the true upper bound $\text{UB}^{\text{E}_\text{MCC}}$ with respect to the original MIQP objective (equation \ref{eq:objective_function}). The third method, denoted as $\text{L}_\text{BIL}$, employs Lagrangian dual decomposition where the subproblems $SP_j$ are solved by using the bilinear formulation. It produces a valid lower bound $\text{LB}^{\text{L}_\text{BIL}}$ from the Lagrangian dual problem and obtains a valid upper bound $\text{UB}^{\text{L}_\text{BIL}}$ via primal heuristics. The fourth method, denoted as $\text{L}_\text{MCC}$, is similar to $\text{L}_\text{BIL}$ but solves the subproblems $SP_j$ by relaxing the bilinear term $W_{jk} x_{ijk}$ with McCormick envelopes. This method provides a lower bound $\text{LB}^{\text{L}_\text{MCC}}$ and an upper bound $\text{UB}^{\text{L}_\text{MCC}}$ again via primal heuristic.
For each method we calculate the model gap,  denoted as $\text{Gap}^m$, where $m \in \{\text{E}_\text{MIQP}, \text{E}_\text{MCC}, \text{L}_\text{BIL}, \text{L}_\text{MCC}\}$. The model gap is calculated as $\text{Gap}^m = \frac{|\text{UB}^m - \text{LB}^m|}{\text{UB}^m} \times 100\%$. For the first method, $\text{E}_\text{MIQP}$, the model gap $\text{Gap}^{\text{E}_\text{MIQP}}$ is straightforward to compute since both the upper and lower bounds $\text{UB}^{\text{E}_\text{MIQP}}$ and $\text{LB}^{\text{E}_\text{MIQP}}$ are directly obtained from the solver. In the second method $\text{E}_\text{MCC}$, the upper bound $\text{UB}^{\text{E}_\text{MCC}}$ is recalculated with the original MIQP objective (Equation \ref{eq:objective_function}) to accurately calculate $\text{Gap}^{\text{E}_\text{MCC}}$. For both Lagrangian methods $\text{L}_\text{BIL}$ and $\text{L}_\text{MCC}$ the valid upper and lower bounds $\text{UB}^m$ and $\text{LB}^m$ allow for direct calculation of the primal dual gaps $\text{Gap}^{\text{L}_\text{BIL}}$ and $\text{Gap}^{\text{L}_\text{MCC}}$. By comparing the model gaps across these four methods, we can analyze their computational performance and effectiveness in solving five instances of SCLA of a particular scenario. Additionally we introduce a secondary gap metric $\text{Gap}^t_m$ to assess each method against the tightest lower bound achieved across all approaches. This bound $LB^*$ is the maximum of all lower bounds: $LB^* = \max \left\{ \text{LB}^{E_{\text{MIQP}}}, \text{LB}^{E_{\text{MCC}}}, \text{LB}^{L_{\text{BIL}}} \text{LB}^{L_{\text{MCC}}} \right\}$. Using $LB^*$, we calculate $\text{Gap}^t_m$ by comparing each method's upper bound and calculate $\text{Gap}^t_m = \frac{|\text{UB}^m - LB^*|}{\text{UB}^m} \times 100\%$. This metric provides a measure of how each method's upper bound approaches the best-known lower bounds. By making use of both $\text{Gap}^m$ and $\text{Gap}^t_m$, we gain an understanding of the computational performance and efficiency of each method in solving the SCLA problem. We average the $\text{Gap}^m$ and $\text{Gap}^t_m$ over five instances of a particular scenario and report the results in \Cref{tab:exact_vs_lagrangean}.

\begin{table*}[!htb]
\centering
\caption{ Comparison of MIQP and Lagrangian solution methods 
}
\label{tab:exact_vs_lagrangean}
\begin{center}
\scriptsize
\begin{tabular}{@{}cc|cc|cc|cc|cc@{}}
\toprule
\multicolumn{2}{c|}{Scenario} & \multicolumn{2}{c|}{$\text{E}_\text{MIQP}$} & \multicolumn{2}{c|}{$\text{E}_\text{MCC}$} & \multicolumn{2}{c|}{$\text{L}_\text{BIL}$} & \multicolumn{2}{c}{$\text{L}_\text{MCC}$} \\
\midrule
$|\mathcal{I}|$ & $k_c$ & $\text{Gap}^m$ (\%) & $\text{Gap}^t$ (\%) & $\text{Gap}^m$ (\%) & $\text{Gap}^t$ (\%) & $\text{Gap}^m$ (\%) & $\text{Gap}^t$ (\%) & $\text{Gap}^m$ (\%) & $\text{Gap}^t$ (\%) \\
\midrule
50 & 1  & $<$1 & $<$1 & $<$1 & $<$1 & $<$1 & $<$1  & $<$1 & $<$1 \\
100 & 1  & $<$1 & $<$1 & $<$1 & $<$1 & $<$1 & $<$1 & 4.08 & $<$1 \\
200 & 1  & 17.67 & $<$1 & 3.91 & $<$1 & $<$1 & $<$1 & 12.49  & $<$1 \\
400 & 1  & 37.80 & $<$1 & 24.75 & $<$1  & $<$1 & $<$1 & 27.69 & $<$1 \\
800 & 1  & 55.94 & 1.03 & 53.44 & $<$1 & $<$1  & $<$1 &  46.79 & $<$1 \\
1600 & 1  & 74.84 & 1.12 & 68.28 & $<$1& $<$1 & $<$1 & 65.47 & $<$1 \\
3200 & 1  & 85.80 & $<$1 & 82.27 & 1.26 & $<$1 & $<$1 & 79.82  & $<$1 \\
\midrule
50 & 2  & $<$1 & $<$1 & $<$1 & $<$1 & 1.73 & 1.42 & 6.19 & $<$1  \\
100 & 2  & 3.79 & $<$1 & 1.36 & $<$1 & 5.38 & 5.06 & 10.20 & $<$1 \\
200 & 2  & 18.81 & 1.21 & 16.71 & 1.24  & 14.26 & 14.26 & 18.59 & 1.23 \\
400 & 2  & 41.44 & 1.51 & 39.00 & 4.54  & 16.46 & 16.46 & 34.11 & 1.59 \\
800 & 2  & 62.49 & 2.25 & 55.53 & 3.5 & 20.49 & 20.49 & 52.39 & 3.06 \\
1600 & 2  & 73.28 & 2.17 & 71.19 & 2.45  & 30.22 & 30.22 & 69.38 & 2.46 \\
3200 & 2  & 85.60 & 25.86 & 88.21 & 40.66 & 50.12 & 50.12 & 82.62 & 21.84 \\
\midrule
50 & 4  & $<$1 & $<$1 & $<$1 & $<$1 & 5.35 & 4.25 & 10.00 & 1.75  \\
100 & 4  & 9.40 & $<$1 & 6.65 & $<$1 & 4.37  & 4.24 & 17.33 & $<$1 \\
200 & 4  & 31.93 & 1.41 & 31.19 & 1.49 & 8.17 & 8.17 & 34.53 & 3.13 \\
400 & 4  & 58.88 & 13.08 & 49.26 & 3.41 & 10.59 & 10.59 & 49.40 & 4.94 \\
800 & 4  & 66.93 & 3.66 & 63.84 & 4.65 & 14.96 & 14.96 & 68.14 & 10.41 \\
1600 & 4  & 81.12 & 7.72 & 78.21 & 9.23 & 26.91 & 26.91 & 80.13 & 9.80 \\
3200 & 4 & 90.73 & 33.12 & 91.08 & 38.65 & 57.31 & 57.31 & 89.96 & 23.94 \\
\midrule
50 & 8  & 1.24 & $<$1 & $<$1 & $<$1 & 7.11 & 6.12 & 18.60 & 4.61 \\
100 & 8  & 22.70 & 3.69 & 20.77 & 3.83 & 14.68 & 14.68 & 30.71 & 6.30  \\
200 & 8  & 44.25 & 8.68 & 37.87 & 8.78 & 22.56 & 22.56 & 48.60 & 16.03 \\
400 & 8  & 58.22 & 8.16 & 52.06 & 8.60 & 24.36 & 24.36 & 61.89 & 16.69 \\
800 & 8  & 74.92 & 7.30 & 68.79 & 7.32 & 29.25 & 29.25 & 73.66 & 15.69 \\
1600 & 8  & 85.39 & 23.14 & 83.23 & 26.39 & 61.19 & 61.19 & 85.30 & 27.81 \\
3200 & 8  & 93.10 & 48.01 & 92.89 & 51.43 & 86.47 & 86.70 & 92.77 & 46.47 \\
\bottomrule
\end{tabular}
\end{center}
\end{table*}

In \Cref{tab:exact_vs_lagrangean_size} we compare the problem sizes in terms of integer, binary, and continuous variables for all four methods $E_{\text{MIQP}}$, $E_{\text{MCC}}$, $L_{\text{BIL}}$, and $L_{\text{MCC}}$ across various scenarios characterized by $|\mathcal{I}|$ and $k_{\text{c}}$. To assess the worst-case computational burden within the Lagrangian framework, we identify $SP_j^{\max}$ as the subproblem with the largest number of continuous variables among all subproblems $SP_j$ for each scenario. We detail this comparison in \Cref{tab:exact_vs_lagrangean_size}. In \Cref{tab:times_best_sol} we present the number of instances (out of five) where each method achieved the tightest upper bound (within 1\%) within one hour. In the subsequent subsections, we discuss insights from the computational experiments detailed in \Cref{sec:exact_vs_lagrangean}.

\begin{table*}[!htb]
\centering
\caption{Problem size comparison of MIQP and Lagrangian relaxation methods}
\label{tab:exact_vs_lagrangean_size}
\scriptsize
\begin{tabular}{@{}cc|cc|cc|c|cc|cc@{}}
\toprule
\multicolumn{2}{c|}{Scenario} & \multicolumn{2}{c|}{$\text{E}_\text{MIQP}$} & \multicolumn{2}{c|}{$\text{E}_\text{MCC}$} & \multirow{2}{*}{Master (Cont.)} & \multicolumn{2}{c|}{$\text{L}_\text{BIL}$ $SP_{j}^{max}$} & \multicolumn{2}{c}{$\text{L}_\text{MCC}$ $SP_{j}^{max}$} \\
\cline{1-6} \cline{8-11}
$|\mathcal{I}|$ & $k_c$ & Int. (Bi) & Cont. & Int. (Bi) & Cont. &  & Int. (Bi) & Cont. & Int. (Bi) & Cont. \\
\midrule
50 & 1  & 1490 (1430) & 60 & 1490 (1430) & 210 & 54 & 67 (64) & 39 & 67 (64) & 75 \\
100 & 1  & 1640 (1580) & 60 & 1640 (1580) & 360 & 104 & 67 (64) & 69 & 67 (64) & 135 \\
200 & 1  & 1940 (1880) & 60 & 1940 (1880) & 660 & 204 & 67 (64) & 120 & 67 (64) & 237 \\
400 & 1  & 2540 (2480) & 60 & 2540 (2480) & 1260 & 404 & 67 (64) & 243 & 67 (64) & 483 \\
800 & 1  & 3740 (3680) & 60 & 3740 (3680) & 2460 & 804 & 67 (64) & 477 & 67 (64) & 951 \\
1600 & 1  & 6140 (6080) & 60 & 6140 (6080) & 4860 & 1604 & 67 (64) & 891 & 67 (64) & 1779 \\
3200 & 1  & 10940 (10880) & 60 & 10940 (10880) & 9660 & 3204 & 67 (64) & 1707 & 67 (64) & 3411 \\
\midrule
50 & 2  & 1556 (1496) & 60 & 1556 (1496) & 276 & 54 & 67 (64) & 45 & 67 (64) & 87 \\
100 & 2  & 1808 (1748) & 60 & 1808 (1748) & 528 & 104 & 67 (64) & 72 & 67 (64) & 141 \\
200 & 2  & 2249 (2189) & 60 & 2249 (2189) & 969 & 204 & 67 (64) & 150 & 67 (64) & 297 \\
400 & 2  & 3197 (3137) & 60 & 3197 (3137) & 1917 & 404 & 67 (64) & 294 & 67 (64) & 585 \\
800 & 2  & 5081 (5021) & 60 & 5081 (5021) & 3801 & 804 & 67 (64) & 615 & 67 (64) & 1227 \\
1600 & 2  & 8807 (8747) & 60 & 8807 (8747) & 7527 & 1604 & 67 (64) & 1197 & 67 (64) & 2391 \\
3200 & 2  & 16226 (16166) & 60 & 16226 (16166) & 14946 & 3204 & 67 (64) & 2325 & 67 (64) & 4647 \\
\midrule
50 & 4  & 1757 (1697) & 60 & 1757 (1697) & 477 & 54 & 67 (64) & 69 & 67 (64) & 135 \\
100 & 4  & 2204 (2144) & 60 & 2204 (2144) & 924 & 104 & 67 (64) & 129 & 67 (64) & 255 \\
200 & 4  & 3068 (3008) & 60 &3068 (3008) & 1788 & 204 & 67 (64) & 234 & 67 (64) & 465 \\
400 & 4  & 4736 (4676) & 60 & 4736 (4676) & 3456 & 404 & 67 (64) & 420 & 67 (64) & 837 \\
800 & 4  & 8090 (8030) & 60 & 8090 (8030) & 6810 & 804 & 67 (64) & 822 & 67 (64) & 1641 \\
1600 & 4  & 14918 (14858) & 60 & 14918 (14858) & 13638 & 1604 & 67 (64) & 1665 & 67 (64) & 3327 \\
3200 & 4  & 28517 (28457) & 60 & 28517 (28457) & 27237 & 3204 & 67 (64) & 3327 & 67 (64) & 6651 \\
\midrule
50 & 8  & 2033 (1973) & 60 & 2033 (1973) & 753 & 54 & 67 (64) & 87 & 67 (64) & 171 \\
100 & 8  & 2747 (2687) & 60 & 2747 (2687) & 753 & 104 & 67 (64) & 171 & 67 (64) & 339 \\
200 & 8  & 4049 (3989) & 60 & 4049 (3989) & 2769 & 204 & 67 (64) & 336 & 67 (64) & 669 \\
400 & 8  & 6704 (6644) & 60 & 6704 (6644) & 5424 & 404 & 67 (64) & 648 & 67 (64) & 1293 \\
800 & 8 & 12263 (12203) & 60 & 12263 (12203) & 10983 & 804 & 67 (64) & 1308 & 67 (64) & 2613 \\
1600 & 8  & 23207 (23147) & 60 & 23207 (23147) & 21927 & 1604 & 67 (64) & 2496 & 67 (64) & 4989 \\
3200 & 8  & 44453 (44393) & 60 & 44453 (44393) & 43173 & 3204 & 67 (64) & 4896 & 67 (64) & 9789 \\
\bottomrule
\end{tabular}%
\end{table*}

\begin{table}[!htb]
\centering
\caption{ Number of instances each method reached the tightest upper bound within 1 hour}
\label{tab:times_best_sol}
\scriptsize
\begin{tabular}{@{}cccccc@{}}
\toprule
$|\mathcal{I}|$ & $k_c$ & $\text{E}_\text{MIQP}$ & $\text{E}_\text{MCC}$ & $\text{L}_\text{BIL}$ & $\text{L}_\text{MCC}$ \\
\midrule
50 & 1  & 5/5 & 5/5 & 5/5 & 5/5 \\
100 & 1  & 5/5 & 5/5 & 5/5 & 5/5 \\
200 & 1  & 4/5 & 5/5 & 5/5 & 5/5 \\
400 & 1  & 5/5 & 3/5 & 5/5 & 5/5 \\
800 & 1  & 3/5 & 5/5 & 5/5 & 5/5 \\
1600 & 1  & 3/5 & 5/5 & 5/5 & 5/5 \\
3200 & 1  & 4/5 & 3/5 & 5/5 & 5/5 \\
\midrule
50 & 2  & 5/5 & 5/5 & 3/5 & 4/5 \\
100 & 2  & 5/5 & 5/5 & 1/5 & 5/5 \\
200 & 2  & 5/5 & 5/5 & 0/5 & 5/5 \\
400 & 2  & 5/5 & 4/5 & 0/5 & 5/5 \\
800 & 2  & 5/5 & 4/5 & 0/5 & 2/5 \\
1600 & 2  & 5/5 & 4/5 & 0/5 & 4/5 \\
3200 & 2  & 0/5 & 0/5 & 0/5 & 5/5 \\
\midrule
50 & 4  & 5/5 & 5/5 & 3/5 & 4/5 \\
100 & 4  & 5/5 & 5/5 & 3/5 & 4/5 \\
200 & 4  & 5/5 & 5/5 & 0/5 & 2/5 \\
400 & 4  & 4/5 & 5/5 & 1/5 & 2/5 \\
800 & 4  & 4/5 & 4/5 & 0/5 & 0/5 \\
1600 & 4  & 3/5 & 4/5 & 0/5 & 0/5 \\
3200 & 4  & 0/5 & 0/5 & 0/5 & 5/5 \\
\midrule
50 & 8  & 5/5 & 5/5 & 3/5 & 3/5 \\
100 & 8  & 5/5 & 5/5 & 0/5 & 2/5 \\
200 & 8  & 5/5 & 4/5 & 1/5 & 1/5 \\
400 & 8  & 5/5 & 4/5 & 0/5 & 0/5 \\
800 & 8  & 5/5 & 5/5 & 0/5 & 0/5 \\
1600 & 8  & 3/5 & 3/5 & 0/5 & 0/5 \\
3200 & 8  & 1/5 & 1/5 & 0/5 & 4/5 \\
\midrule
Total &  & 114/140 & 113/140 & 50/140 & 93/140 \\
\bottomrule
\end{tabular}
\end{table}

\subsubsection{Effect of $k_c$ and $|\mathcal{I}|$ on Model Gap ($\text{Gap}^m$)}

For the exact methods ($\text{E}_\text{MIQP}$ and $\text{E}_\text{MCC}$), the model gap $\text{Gap}^m$ remains low (less than 1\%) for small problem sizes (e.g., $|\mathcal{I}| = 50$ and $k_c = 1$). As $|\mathcal{I}|$ increases, however, $\text{Gap}^m$ increases significantly. For instance, at $|\mathcal{I}| = 200$ and $k_c = 1$, $\text{E}_\text{MIQP}$ has a $\text{Gap}^m$ of 17.67\%, while $\text{E}_\text{MCC}$ improves on this with a gap of 3.91\%. When $|\mathcal{I}| = 3200$, the gaps become substantial, exceeding 80\% for both exact methods, indicating that they struggle to find tight lower bounds within the time limit for larger problems. Increasing $k_c$ exacerbates the computational difficulty for the exact methods. With more charging stations considered per household, the problem's connectivity increases, leading to a higher number of variables and constraints (see \Cref{tab:exact_vs_lagrangean_size}). For example, at $|\mathcal{I}| = 200$ and $k_c = 8$, $\text{E}_\text{MIQP}$ exhibits a $\text{Gap}^m$ of 44.25\%, and $\text{E}_\text{MCC}$ has a gap of 37.87\%. At $|\mathcal{I}| = 3200$ and $k_c = 8$, the model gaps reach over 92\% for both methods. The Lagrangian methods ($\text{L}_\text{BIL}$ and $\text{L}_\text{MCC}$) display better scalability with respect to $|\mathcal{I}|$ and $k_c$. $\text{L}_\text{BIL}$ maintains low model gaps for small $k_c$ values across all $|\mathcal{I}|$. Even at $|\mathcal{I}| = 3200$ and $k_c = 1$, $\text{L}_\text{BIL}$ achieves a $\text{Gap}^m$ of less than 1\%. As $k_c$ increases, however, the model gaps for Lagrangian methods also increase but remain lower than those of the exact methods. For high $k_c$ values, $\text{L}_\text{MCC}$ sometimes outperforms $\text{L}_\text{BIL}$ in reaching tighter upper bounds.

\subsubsection{Gap Relative to Tightest Lower Bound ($\text{Gap}^t$)}
For the exact methods ($\text{E}_\text{MIQP}$ and $\text{E}_\text{MCC}$), $\text{Gap}^t$ is generally significantly lower than their model gaps $\text{Gap}^m$, especially in larger problem instances. The reason  is that while the exact methods may struggle to find tight lower bounds themselves, they benefit from the tighter lower bounds obtained by the Lagrangian methods, particularly $\text{L}_\text{MCC}$. Examining \Cref{tab:exact_vs_lagrangean}, at $|\mathcal{I}| = 200$ and $k_c = 1$, we see that $\text{E}_\text{MIQP}$ has a model gap $\text{Gap}^m$ of 17.67\%, but its gap relative to the tightest lower bound $\text{Gap}^t$ is less than 1\%. This significant reduction indicates that another method ($\text{L}_\text{BIL}$ in this case) has provided a much tighter lower bound, effectively improving the perceived quality of $\text{E}_\text{MIQP}$'s solution when compared with the best-known solution. As the problem size increases, the contribution of the Lagrangian methods to the tightest lower bound becomes greater. Specifically, $\text{L}_\text{MCC}$ demonstrates better performance in providing tight lower bounds for larger problem instances with higher $k_c$ values. For example, at $|\mathcal{I}| = 3200$ and $k_c = 8$, $\text{L}_\text{MCC}$ achieves a $\text{Gap}^m$ of 92.77\% and a $\text{Gap}^t$ of 46.47\%, which is significantly better than the exact methods. In contrast, $\text{E}_\text{MIQP}$ has a $\text{Gap}^m$ of 93.10\% and a $\text{Gap}^t$ of 48.01\%, indicating that $\text{L}_\text{MCC}$ has provided a tighter lower bound that benefits all methods when computing $\text{Gap}^t$. Furthermore, \Cref{tab:times_best_sol} shows that $\text{L}_\text{MCC}$ reaches the tightest upper bound in 4 out of 5 instances at $|\mathcal{I}| = 3200$ and $k_c = 8$, outperforming the exact methods, which  achieve this only in 1 out of 5 instances. This demonstrates $\text{L}_\text{MCC}$'s effectiveness  not only in providing tight lower bounds but also in finding high-quality feasible solutions in large problems.

\subsubsection{Problem Size Comparison} 
From \Cref{tab:exact_vs_lagrangean_size} we can see that for the exact methods ($\text{E}_\text{MIQP}$ and $\text{E}_\text{MCC}$), the number of integer variables escalates dramatically as both the number of households $|\mathcal{I}|$ and the number of nearest charging stations $k_c$ increase. For instance, at $|\mathcal{I}| = 3200$ and $k_c = 8$, the exact methods require over 44,000 integer variables and approximately 43,000 continuous variables. Constructing such large MIQP models is not only computationally intensive but often infeasible in practice because of memory limitations and processing constraints of optimization solvers such as Gurobi. In many cases, especially at this scale, the solver cannot even build the full model, rendering the exact methods impractical for large-scale problems. In contrast, the Lagrangian methods ($\text{L}_\text{BIL}$ and $\text{L}_\text{MCC}$) decompose the original problem into smaller, more manageable subproblems $SP_j$, each associated with a charging station $j \in \mathcal{J}$. The master problem remains relatively small, containing only continuous variables (e.g., 3,204 variables when $|\mathcal{I}| = 3200$). Each subproblem $SP_j$ involves a fixed number of integer variables (67 in all scenarios) and a manageable number of continuous variables that scale linearly with $|\mathcal{I}|$ and $k_c$. However, since the number of households associated with each charging station varies, the size of each subproblem $SP_j$ can differ. Specifically, at $|\mathcal{I}| = 3200$ and $k_c = 8$, the largest subproblem $SP_j^{\max}$ in $\text{L}_\text{MCC}$ has only 67 integer variables and approximately 9,789 continuous variables. In \Cref{tab:exact_vs_lagrangean_size} we report the problem size for $SP_j^{\max}$ to illustrate the maximum computational effort required for any single subproblem within each scenario. This decomposition significantly reduces the computational burden, enabling the Lagrangian methods to handle extremely large-scale instances effectively where the exact methods cannot even construct the model. The sharp contrast in problem sizes highlights the scalability advantage of the Lagrangian methods. While the exact methods' problem size grows exponentially with $|\mathcal{I}|$ and $k_c$, making them unsuitable for large problems, the Lagrangian methods maintain a linear growth in subproblem sizes. By focusing on the maximum subproblem size $SP_j^{\max}$, we demonstrate that even the most demanding subproblem remains computationally tractable. Thus, the Lagrangian methods to exploit problem structure and computational resources more efficiently,  being able to solve instances that are beyond the capability of the exact methods.

\subsubsection{Ability to Reach the Tightest Upper Bound}

From \Cref{tab:times_best_sol}, we can see that the exact methods achieve the tightest upper bound consistently in smaller problem instances. However, their performance deteriorates as both $|\mathcal{I}|$ and $k_c$ increase. For example, at $|\mathcal{I}| = 3200$ and $k_c = 8$, $\text{E}_\text{MIQP}$ reaches the tightest upper bound in only one out of five instances. In contrast, the Lagrangian methods demonstrate superior performance in larger instances. Specifically, $\text{L}_\text{MCC}$ reaches the tightest upper bound in four out of five instances at $|\mathcal{I}| = 3200$ and $k_c = 8$, outperforming all other methods. This result shows the effectiveness of the Lagrangian dual decomposition approach in scenarios that are both large scale and highly connected. While the exact methods are well suited for small-scale problems with a low $k_c$, providing high-quality solutions within reasonable computational times, they suffer from a lack of scalability due to the exponential increase in problem size with larger $|\mathcal{I}|$ and $k_c$. On the other hand, the Lagrangian methods, which leverage problem decomposition, offer a scalable alternative capable of handling larger problem instances effectively. They maintain acceptable model gaps and are likely to reach an acceptable upper bound within the given time limit, especially in scenarios with a high $k_c$. The choice between bilinear subproblems and McCormick relaxation within the Lagrangian framework depends on the specific scenario, where $\text{L}_\text{BIL}$ tends to perform better for small to moderate $k_c$ values, while $\text{L}_\text{MCC}$ may be more advantageous for higher $k_c$ values in terms of reaching the tightest upper bounds.

\subsection{Charging Location Strategies and Policy Frameworks}
In this section we analyze the impact of different charging location strategies and policy frameworks on the objective ($\ref{eq:objective_function}$). By evaluating various scenarios, we aim to provide managerial insights into optimizing infrastructure investments and operational policies. The combination of charging location strategies and management policies enables us to evaluate their combined impact on the total objective of SCLA (objective $\ref{eq:objective_function}$).  

\begin{table*}[ht]
\centering
\scriptsize
\caption{Percentage change in the objective for various scenarios using $\text{L}_{\text{MCC}}$ method}
\label{tab:scenario-costs}

\begin{tabular}{c cc cc cccc}
\toprule
Scenario &
\multicolumn{2}{c}{Charging Location} &
\multicolumn{2}{c}{Policy} &
\multicolumn{4}{c}{\% $\uparrow$ in $\mathbb{C}$} \\
\cmidrule(lr){2-3} \cmidrule(lr){4-5} \cmidrule(lr){6-9}
 & TAZs & Depots & Single & Multi & Mean & Std Dev & Min & Max \\
\midrule
1 & $\checkmark$ & $\times$ & $\checkmark$ & $\times$ & 91.24 & 7.82 & 82.79 & 100.74 \\
2 & $\checkmark$ & $\times$ & $\times$ & $\checkmark$ & 44.82 & 10.08 & 32.61 & 62.46 \\
3 & $\times$ & $\checkmark$ & $\checkmark$ & $\times$ & 42.33 & 11.14 & 27.39 & 60.93 \\
4 & $\times$ & $\checkmark$ & $\times$ & $\checkmark$ & 32.23 & 3.95 & 25.34 & 35.93 \\
5 & $\checkmark$ & $\checkmark$ & $\checkmark$ & $\times$ & 14.51 & 3.97 & 8.07 & 19.18 \\
6 & $\checkmark$ & $\checkmark$ & $\times$ & $\checkmark$ & - & - & - & - \\
\bottomrule
\end{tabular}
\end{table*}

\Cref{tab:scenario-costs} demonstrates that strategic decisions on charging infrastructure placement at new locations (TAZs) versus existing facilities (depots) and operational policies significantly impact total costs compared with the baseline configuration (Scenario 6). For each scenario we generate five random instances with 100 households each ($|\mathcal{I}|$ = 100), and the results presented are averaged over these five instances. To account for detour times across all households, we scale the detour cost parameter \(C^\delta\) by a fixed multiplier. The transition from Scenario 6 (optimal baseline) to Scenario 5 involves switching from multiagency to single-agency policy while keeping combined TAZs and depot locations, which increases total costs by 14.51\%. The switch from the baseline to Scenario 4 restricts charging stations to only existing depot locations under a multiagency policy and increases costs by 32.23\%. Scenario 3 applies a single-agency policy to depot-only locations, which increases costs by 42.33\%. Similarly, Scenario 2 constrains locations to only new potential TAZ locations while maintaining a multiagency policy, which drives costs up by 44.82\%. Scenario 1 represents the most restrictive case, combining TAZ-only locations with single-agency policy, which produces the highest cost increase of 91.24\% from the baseline.

This systematic pattern of cost increases across scenarios offers key insights for infrastructure planning and policy development. The cost increase from the baseline Scenario 6 to Scenario 1 (91.24\%) demonstrates how limiting charging infrastructure to only new locations and restricting operational flexibility through single-agency policy lead to substantial cost increases. Policymakers should recognize that the shift from multiagency policy alone (Scenario 5) produces the smallest cost increase, whereas restricting charging stations to either only new locations or only existing facilities creates larger cost impacts. We believe that for long-term sustainability policymakers should prioritize implementation plans that consider both new potential locations (TAZs) and existing facilities (depots) for charging infrastructure development rather than restricting installations to only existing facilities. Our analysis demonstrates that combining charging station placement at new strategic locations with existing depot facilities while enabling multiagency collaborative operations leads to the most cost-effective outcomes.

\subsection{Large-Scale Case Study}\label{sec:large_scale_implementation}
Building on the insights from the preceding section, which demonstrated the scalability and effectiveness of Lagrangian methods combined with primal heuristics in solving large-scale SCLA problems, we present a comprehensive large-scale case study. The dataset includes 449,367 households, 1,958 potential charging stations, and 53 existing charging depots. This dataset represents the largest ever used in studies of charger location allocation with queue congestion. We assume a multiagency policy, meaning any open station whether at a TAZ or an existing depot can serve any household. To tackle this extensive problem, we apply the Lagrangian dual decomposition approach, employing the McCormick relaxation subproblem formulation (\(\text{L}_{\text{MCC}}\)).

\subsubsection{Multipartition Parallelization Strategy}

In our large-scale SCLA implementation, we faced the challenge of solving $2011$ station-level subproblems $\{SP_j : j \in \mathcal{J}\}$ across a heterogeneous HPC environment. While an ideal solution would allocate one CPU core per subproblem for simultaneous execution, our available HPC infrastructure consisted of three distinct partitions with more limited resources. The Short Partition ($P_S$) provided $400$ CPU cores with a one-hour runtime limit, while the Group Partition ($P_G$) provided $180$ CPU cores without time constraints, and the Open Partition ($P_O$) provided $100$ CPU cores with a $48$-hour time limit. Communication between these partitions is not possible in our HPC environment. This resource configuration requires strategic distribution of the subproblems across partitions. An example allocation may look like assigning approximately $800$ subproblems to $P_G$, $600$ to $P_O$, and the remaining $611$ to $P_S$. The exact numbers may vary, but each partition faces more subproblems than CPU cores; and so, to address this situation, we employ a hierarchical master-worker architecture and a partitioning and scheduling strategy.

\paragraph{Global master (Lagrangian dual controller)}: The \textit{global master} runs on one of the cores in the group partition \(P_G\). It initializes the multipliers \((\zeta,\beta_\phi,\beta_\xi,\nu)\) and randomly partitions \(\mathcal{J}\) into three subsets \(\mathcal{J}_S,\mathcal{J}_G,\mathcal{J}_O\), assigning each subset to \(P_S, P_G,\) and \(P_O\), respectively. Each subproblem \(SP_j\) is defined over a station \(j\), its corresponding set \(\mathcal{I}_j\) of households, and multiple charger types, resulting in varying problem sizes and complexities across the subproblems.

\paragraph{Partition masters and workload distribution}:  
Within each partition, a \textit{partition master} receives its assigned subset of subproblems $\mathcal{J}_p$ ($p \in \{S,G,O\}$). Since $|\mathcal{J}_p|$ exceeds the available CPU cores in each partition, the partition master must distribute multiple subproblems to each worker. For load balancing, we first sort the subproblems in $\mathcal{J}_p$ by descending order of $|\mathcal{I}_j|$ (number of households assigned to location $j$), using this as a proxy for computational complexity. The sorted subproblems are then distributed once using round-robin assignment: if there are $w$ workers, worker $i$ receives its complete batch of problems $\{j_{i}, j_{i+w}, j_{i+2w}, ...\}$.  

\paragraph{Execution and completion of subproblems}:  
Each worker then processes its assigned batch independently, without further communication with the \textit{partition master}. The \textit{partition master} simply waits for all workers to complete their assigned problems before declaring completion. If any subproblem assigned \(P_S\) fails to complete, its index is recorded, and the \textit{global master} subsequently resubmits it to the group partition \(P_G\). This mechanism ensures that every subproblem is ultimately solved. 

\paragraph{Iterative updates by the global master}:  
After all subproblems in \(\mathcal{J}_S,\mathcal{J}_G,\mathcal{J}_O\) are completed, the global master collects the results. At this stage, it performs subgradient updates, considering a subgradient calculation time \(T^{\text{subgrad}}\), and executes primal heuristics with a heuristic execution time \(T^{\text{heuristic}}\). These adjustments to the multipliers may prompt another iteration of subproblem assignments or lead to termination once convergence criteria are met. 
Our multipartition parallelization strategy leverages random partitioning of \(\mathcal{J}\), hierarchical control, sorted-round-robin scheduling, and a mechanism for reassigning failed tasks. Even with far fewer cores than subproblems and differing time limits and complexities, this approach allows us to efficiently and systematically solve all subproblems over time. This enables the \textit{global master} to advance through the Lagrangian dual iterations. We present the results below.

\begin{table*}[ht]
\centering
\caption{Large-scale case study results using the \(\text{L}_\text{MCC}\) method with a 6-hour time limit}
\label{tab:scenario-results}
\scriptsize
\resizebox{\textwidth}{!}{%
\begin{tabular}{cc cc c c ccc ccc c}
\toprule
\multicolumn{2}{c}{Charging Location} &
\multicolumn{2}{c}{Policy} &
$|\mathcal{J}|$ &
$J^{\text{open}}$ &
\multicolumn{3}{c}{Total Chargers $\left( \sum_{j \in \mathcal{J}} S_{jk} \right)$} &
\multicolumn{3}{c}{$\overline{W}_{jk}$ (min)} &
$\text{Gap}^m$ (\%) \\
\cmidrule(lr){1-2} \cmidrule(lr){3-4} \cmidrule(lr){7-9} \cmidrule(lr){10-12}
TAZs & Depots & Single & Multi & & &
$k_1$ & $k_2$ & $k_3$ &
$k_1$ & $k_2$ & $k_3$ & \\
\midrule
\checkmark & \checkmark & \(\times\) & \checkmark & 2014 & 436 & 218 & 109 & 302 & 123.57 & 31.53 & 16.59 & 78.5 \\
\bottomrule
\end{tabular}%
}
\scriptsize
\vspace{2cm}
$k_1$: Slow chargers, $k_2$: Moderate chargers, $k_3$: Fast chargers.
\end{table*}

From \Cref{tab:scenario-results} we observe the results of our large-scale SCLA implementation, which successfully solved an instance with 449,367 households and 2,014 potential station locations (including both TAZs and existing depots). Using a 6-hour time limit and the McCormick relaxation (\(\text{L}_{\text{MCC}}\)) formulation, the method identified 436 best station locations out of 2,014 candidates, requiring only 21.6\% of the available sites. Across these opened stations, a total of 629 chargers were deployed---218 slow chargers, 109 moderate chargers, and 302 fast chargers---indicating a preference for higher-speed charging options in the best configuration.

Slow chargers averaged 1.36 chargers per station, with a combined waiting and service time (system time) of 123.58 minutes. Moderate chargers averaged 1.00 charger per station and achieved a shorter system time of 31.53 minutes. Fast chargers also averaged 1.00 charger per station, providing the lowest system time of just 16.60 minutes. From a computational perspective, the maximum Gurobi solve time per subproblem was 3019.04 seconds, and the primal heuristic phase was completed in 179.94 seconds. Although the final optimality gap of 78.5\% is substantial, it is a reasonable outcome given the unprecedented scale of this problem, representing the largest SCLA instance with queue congestion to date.

\begin{figure}
    \centering
    \begin{subfigure}[t]{0.48\textwidth}
        \centering
        \includegraphics[width=\textwidth]{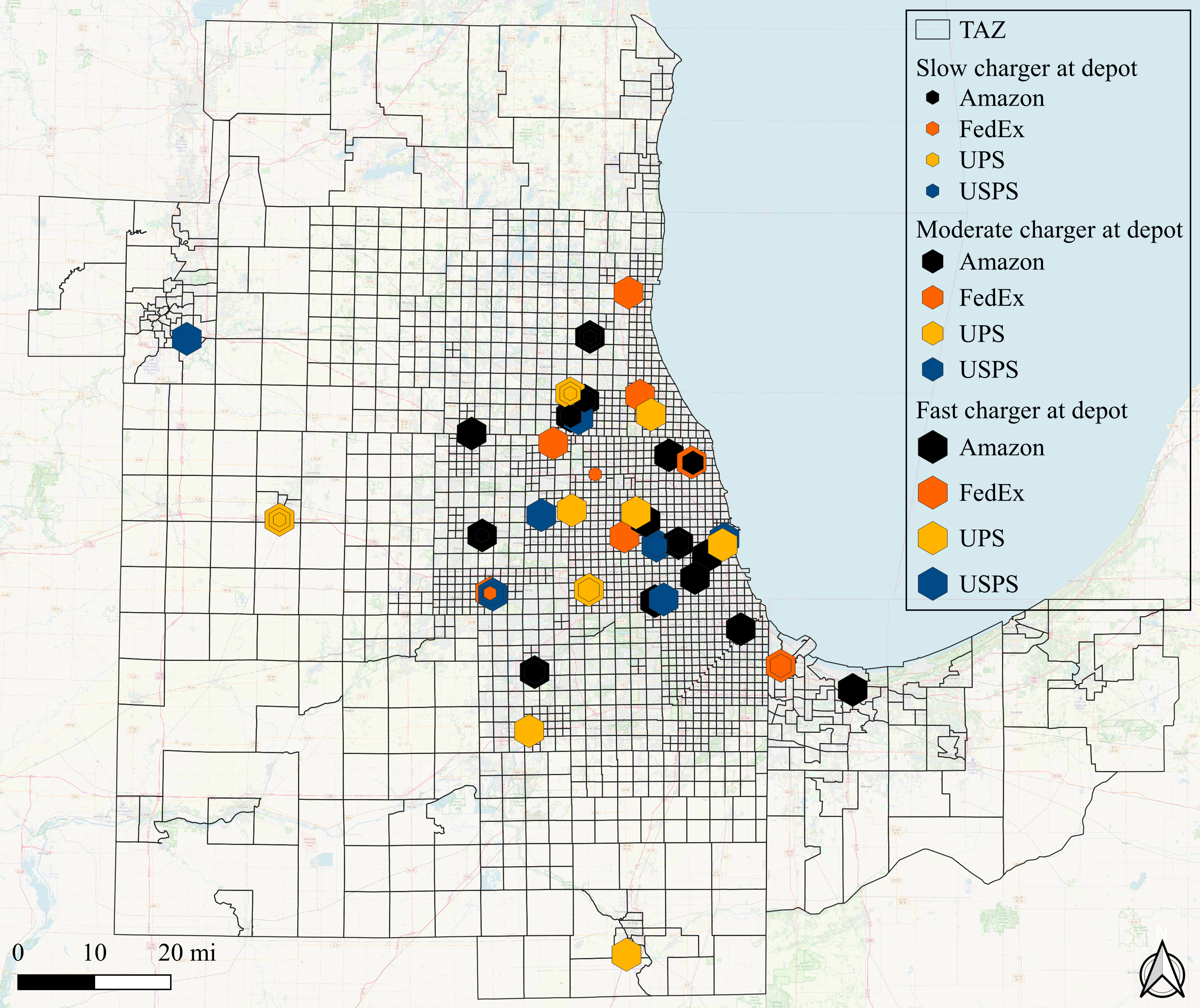}
        \caption{Depots.}\label{fig:open_depots}
    \end{subfigure}
    \quad
    \begin{subfigure}[t]{0.48\textwidth}
        \centering
        \includegraphics[width=\textwidth]{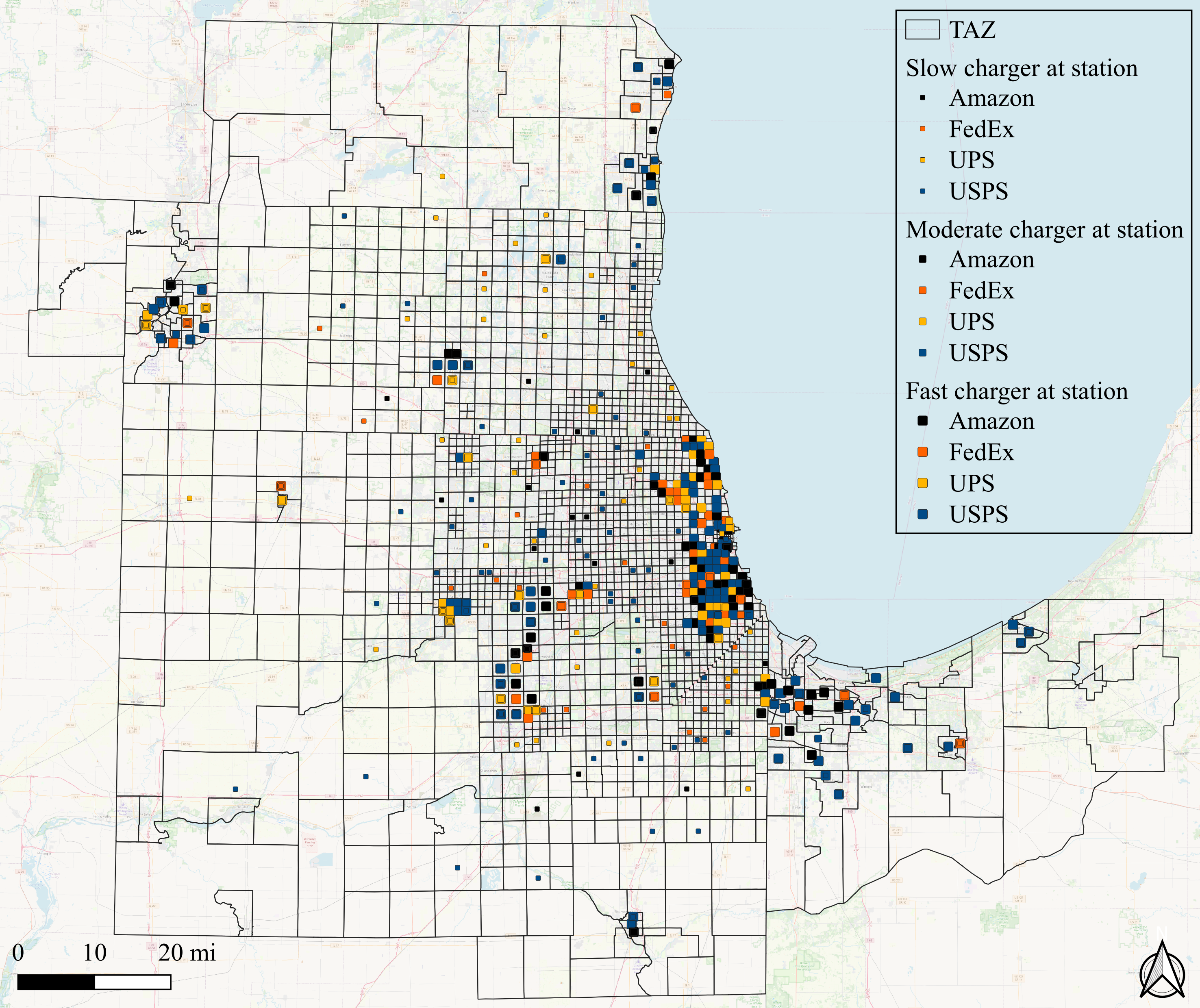}
        \caption{TAZs.}\label{fig:open_TAZs}
    \end{subfigure}

    \caption{Spatial distribution of open depot and TAZ stations}
    \label{fig:spatial_solution}
\end{figure} 

\Cref{fig:spatial_solution} shows the stations selected in the best solution identified. \Cref{fig:spatial_solution} (b) focuses on TAZ-based stations (squares) that are opened in the solution. Each station symbol also indicates which charger types (slow, moderate, or fast) were installed at that location.

Several trends emerge from these figures. First, most depot stations opened in or around the city center (\Cref{fig:spatial_solution} (a)), suggesting higher demands and thus a stronger incentive to co-locate chargers near major demand points. The TAZ-based stations are more widely dispersed across the metropolitan region (\Cref{fig:spatial_solution} (b)). These TAZ stations host a mix of fast and moderate chargers, while depots include more slow chargers. This pattern appears consistent with the model’s goal of balancing high-traffic locations (which benefit from faster chargers to limit queues) with cost-effective service in lower-density areas (where slower chargers may suffice).

The solution indicates that combining depot and TAZ options provides good spatial coverage while still adhering to budget and queue constraints. Many of the TAZ sites offering fast chargers are located near major corridors or in higher-demand neighborhoods, minimizing detour time and wait times for vehicles. Meanwhile, more modest charging installations at peripheral TAZ stations reduce congestion at busier stations without excessively increasing costs. These location and charger-type decisions reflect the model’s attempt to serve 449,367 households efficiently under multiagency operations where any household can use any open station, reinforcing that fast chargers are allocated to busier or central sites and slow or moderate chargers serve moderately traveled or outlying zones.

\subsection{Sensitivity Analysis}

\begin{figure*}[ht]
    \centering
    \begin{tikzpicture}

    
    \begin{axis}[
       name=plot1,
       at={(0,0)},
       anchor=north west,
       width=7.5cm,
       height=6.5cm,
       xlabel={Increase in $C^{\phi}$ (\%)},
       ylabel={Increase in $\mathbb{C}$ (\%)},
       symbolic x coords={0,10,20,30,40,50,80,100,150,200},
       x tick label style={font=\scriptsize},
       xtick=data,
       ymin=0,
       ymax=140,
       legend pos=north west
    ]
       \addplot[
           color=blue!70!black,
           mark=square*,
           mark size=2,
           error bars/.cd,
           y dir=both,
           y explicit,
           error bar style={line width=0.8pt}
       ] coordinates {
           (0,0.00) +- (0,0.00)
           (10,6.56) +- (0,0.14)
           (20,13.12) +- (0,0.28)
           (30,19.68) +- (0,0.41)
           (40,26.25) +- (0,0.55)
           (50,32.81) +- (0,0.69)
           (80,52.49) +- (0,1.10)
           (100,65.62) +- (0,1.38)
           (150,98.42) +- (0,2.06)
           (200,131.23) +- (0,2.75)
       };
       \addlegendentry{$|\mathcal{I}| = 50$}
    
       \addplot[
           color=red!70!black,
           mark=triangle*,
           mark size=2,
           error bars/.cd,
           y dir=both,
           y explicit,
           error bar style={line width=0.8pt}
       ] coordinates {
           (0,0.00) +- (0,0.00)
           (10,5.98) +- (0,0.22)
           (20,11.96) +- (0,0.43)
           (30,17.95) +- (0,0.64)
           (40,23.93) +- (0,0.86)
           (50,29.91) +- (0,1.07)
           (80,47.86) +- (0,1.72)
           (100,59.83) +- (0,2.15)
           (150,89.74) +- (0,3.22)
           (200,119.66) +- (0,4.29)
       };
       \addlegendentry{$|\mathcal{I}| = 100$}
    
       \addplot[
           color=black,
           mark=diamond*,
           mark size=2,
           error bars/.cd,
           y dir=both,
           y explicit,
           error bar style={line width=0.8pt}
       ] coordinates {
           (0,0.00) +- (0,0.00)
           (10,5.23) +- (0,0.16)
           (20,10.41) +- (0,0.27)
           (30,15.63) +- (0,0.41)
           (40,20.85) +- (0,0.55)
           (50,26.07) +- (0,0.70)
           (80,41.71) +- (0,1.12)
           (100,52.20) +- (0,1.46)
           (150,78.38) +- (0,1.99)
           (200,104.36) +- (0,2.87)
       };
       \addlegendentry{$|\mathcal{I}| = 200$}

       \node[anchor=south west,font=\small] at (rel axis cs:0,1.02) {(a)};
    \end{axis}
    \node[below=1cm] at (plot1.south) {(a)};

    \begin{axis}[
        name=plot2,
        at={(8.5cm,0)},
        anchor=north west,
        width=7.5cm,
        height=6.5cm,
        xlabel={Adjustment for $C^\tau$ (\%)},
        ylabel={Increase in $\mathbb{C}$ (\%)},
        symbolic x coords={0,2,4,6,8,10},
        xtick=data,
        ymin=0,
        ymax=28,
        legend pos=north west
    ]
        \addplot[
            color=blue!70!black,
            mark=square*,
            mark size=2,
            error bars/.cd,
            y dir=both,
            y explicit,
            error bar style={line width=0.8pt}
        ] coordinates {
            (0,0.00) +- (0,0.00)
            (2,2.82) +- (0,0.32)
            (4,8.38) +- (0,0.90)
            (6,12.67) +- (0,1.61)
            (8,15.90) +- (0,1.78)
            (10,18.83) +- (0,2.05)
        };
        \addlegendentry{$|\mathcal{I}| = 50$}
    
        \addplot[
            color=red!70!black,
            mark=triangle*,
            mark size=2,
            error bars/.cd,
            y dir=both,
            y explicit,
            error bar style={line width=0.8pt}
        ] coordinates {
            (0,0.00) +- (0,0.00)
            (2,4.26) +- (0,0.36)
            (4,10.06) +- (0,0.94)
            (6,14.90) +- (0,1.20)
            (8,18.82) +- (0,1.35)
            (10,21.93) +- (0,1.53)
        };
        \addlegendentry{$|\mathcal{I}| = 100$}
    
        \addplot[
            color=black,
            mark=diamond*,
            mark size=2,
            error bars/.cd,
            y dir=both,
            y explicit,
            error bar style={line width=0.8pt}
        ] coordinates {
            (0,0.00) +- (0,0.00)
            (2,5.04) +- (0,0.60)
            (4,11.78) +- (0,0.94)
            (6,16.09) +- (0,1.20)
            (8,19.60) +- (0,1.62)
            (10,22.82) +- (0,2.01)
        };
        \addlegendentry{$|\mathcal{I}| = 200$}

        \node[anchor=south west,font=\small] at (rel axis cs:0,1.02) {(b)};
    \end{axis}
    \node[below=1cm] at (plot2.south) {(b)};

    \begin{axis}[
       name=plot3,
       at={(0,-7cm)},
       anchor=north west,
       width=7.5cm,
       height=6.5cm,
       xlabel={Decrease in $C^{\xi}$ (\%)},
       ylabel={Decrease in $\mathbb{C}$ (\%)},
       symbolic x coords={0,10,20,30,40,50,80},
       xtick=data,
       ymin=0,
       ymax=17,
       legend pos=north west
    ]
       \addplot[
           color=blue!70!black,
           mark=square*,
           mark size=2,
           error bars/.cd,
           y dir=both,
           y explicit,
           error bar style={line width=0.8pt}
       ] coordinates {
           (0,0.00) +- (0,0.00)
           (10,1.92) +- (0,0.04)
           (20,3.83) +- (0,0.08)
           (30,5.75) +- (0,0.12)
           (40,7.66) +- (0,0.16)
           (50,9.58) +- (0,0.20)
           (80,15.36) +- (0,0.32)
       };
       \addlegendentry{$|\mathcal{I}| = 50$}
    
       \addplot[
           color=red!70!black,
           mark=triangle*,
           mark size=2,
           error bars/.cd,
           y dir=both,
           y explicit,
           error bar style={line width=0.8pt}
       ] coordinates {
           (0,0.00) +- (0,0.00)
           (10,1.75) +- (0,0.06)
           (20,3.49) +- (0,0.13)
           (30,5.24) +- (0,0.19)
           (40,6.99) +- (0,0.25)
           (50,8.74) +- (0,0.31)
           (80,14.47) +- (0,0.35)
       };
       \addlegendentry{$|\mathcal{I}| = 100$}
    
       \addplot[
           color=black,
           mark=diamond*,
           mark size=2,
           error bars/.cd,
           y dir=both,
           y explicit,
           error bar style={line width=0.8pt}
       ] coordinates {
           (0,0.00) +- (0,0.00)
           (10,1.55) +- (0,0.08)
           (20,3.07) +- (0,0.12)
           (30,4.70) +- (0,0.15)
           (40,6.26) +- (0,0.38)
           (50,8.09) +- (0,0.30)
           (80,14.06) +- (0,0.26)
       };
       \addlegendentry{$|\mathcal{I}| = 200$}

       \node[anchor=south west,font=\small] at (rel axis cs:0,1.02) {(c)};
    \end{axis}
    \node[below=1cm] at (plot3.south) {(c)};

    \begin{axis}[
       name=plot4,
       at={(8.5cm,-7cm)},
       anchor=north west,
       width=7.5cm,
       height=6.5cm,
       xlabel={Increase in $P_k$ (\%)},
       ylabel={Decrease in $\mathbb{C}$ (\%)},
       symbolic x coords={0,20,40,60,80,100},
       xtick=data,
       ymin=0,
       ymax=4.2,
       legend pos=north west
    ]
       \addplot[
           color=blue!70!black,
           mark=square*,
           mark size=2,
           error bars/.cd,
           y dir=both,
           y explicit,
           error bar style={line width=0.8pt}
       ] coordinates {
           (0,0.00) +- (0,0.00)
           (20,0.47) +- (0,0.05)
           (40,0.81) +- (0,0.09)
           (60,1.06) +- (0,0.12)
           (80,1.26) +- (0,0.14)
           (100,1.41) +- (0,0.16)
       };
       \addlegendentry{$|\mathcal{I}| = 50$}
    
       \addplot[
           color=red!70!black,
           mark=triangle*,
           mark size=2,
           error bars/.cd,
           y dir=both,
           y explicit,
           error bar style={line width=0.8pt}
       ] coordinates {
           (0,0.00) +- (0,0.00)
           (20,0.72) +- (0,0.06)
           (40,1.22) +- (0,0.10)
           (60,1.61) +- (0,0.13)
           (80,1.90) +- (0,0.16)
           (100,2.14) +- (0,0.18)
       };
       \addlegendentry{$|\mathcal{I}| = 100$}
    
       \addplot[
           color=black,
           mark=diamond*,
           mark size=2,
           error bars/.cd,
           y dir=both,
           y explicit,
           error bar style={line width=0.8pt}
       ] coordinates {
           (0,0.00) +- (0,0.00)
           (20,1.08) +- (0,0.09)
           (40,1.84) +- (0,0.15)
           (60,2.37) +- (0,0.24)
           (80,2.84) +- (0,0.24)
           (100,3.17) +- (0,0.26)
       };
       \addlegendentry{$|\mathcal{I}| = 200$}

       \node[anchor=south west,font=\small] at (rel axis cs:0,1.02) {(d)};
    \end{axis}
    \node[below=1cm] at (plot4.south) {(d)};
        
    \end{tikzpicture}
    
    \caption{Sensitivity analysis of changes in $\mathbb{C}$ under various parameter adjustments.}
    \label{fig:sen_analysis}
\end{figure*}
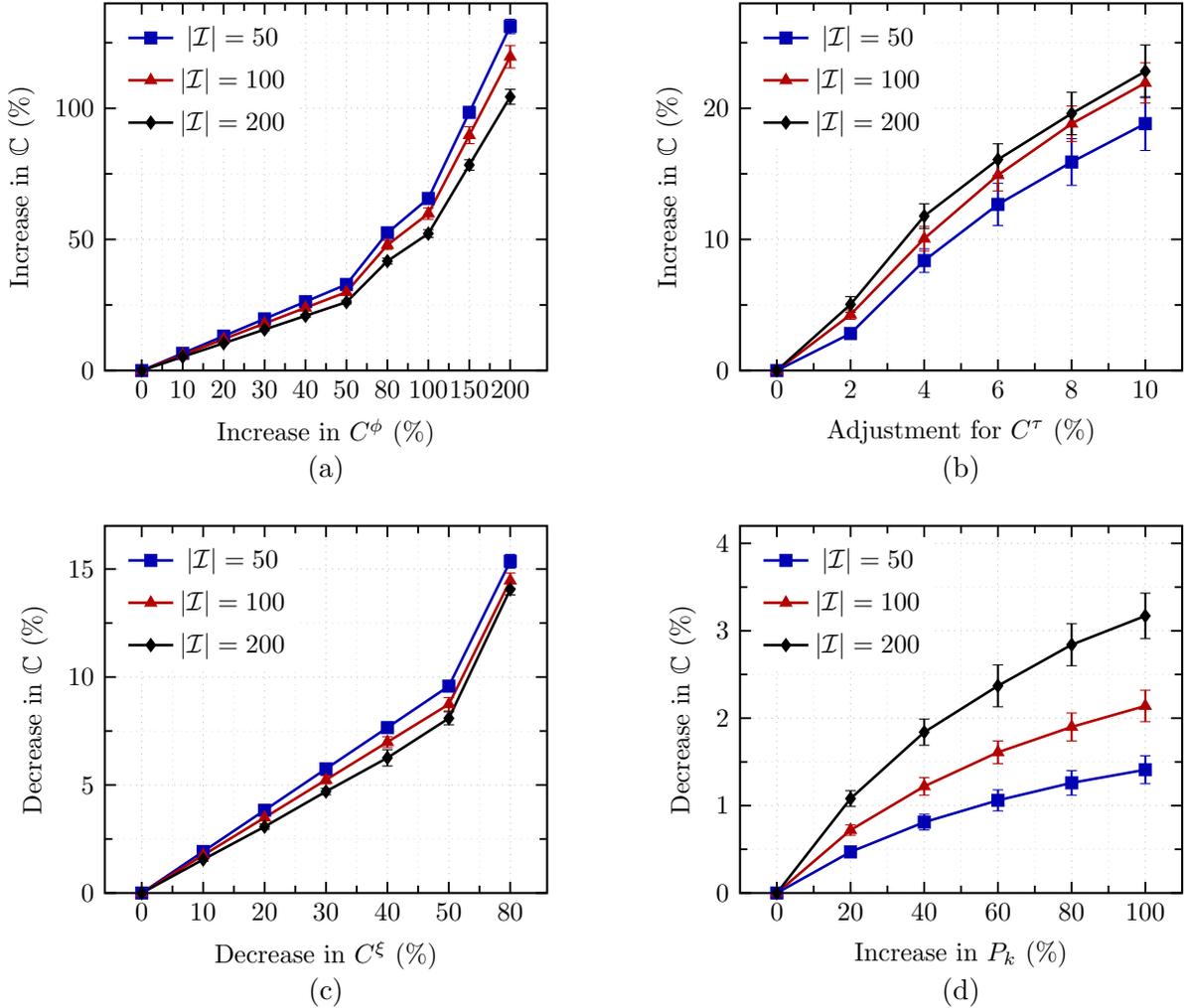

We perform sensitivity analyses to understand how changes in key model parameters impact total costs $(\mathbb{C})$ across different household sizes $(|\mathcal{I}| = 50, 100, 200)$. We assume a multiagency policy, meaning any open station whether at a TAZ or an existing depot can serve any household. For each parameter configuration and household size we generate five random instances; the average results are presented in \Cref{fig:sen_analysis}. The analysis focuses on four main parameters: station costs $(C^{\phi})$, waiting time costs $(C^{\tau})$, charger costs $(C^{\xi})$, and charger powers $(P_k)$ $\forall k \in \mathcal{K}$.  \Cref{fig:sen_analysis}(a) shows that station costs $(C^{\phi})$ have the most pronounced influence on $\mathbb{C}$, especially for scenarios with fewer households. For instance, a 100\% increase in $C^{\phi}$ leads to cost escalations of approximately 52.20\% to 65.62\%, with $|\mathcal{I}|=50$ displaying the greatest sensitivity. This indicates that infrastructure-related expenses disproportionately affect smaller sets of 
households. In contrast, \Cref{fig:sen_analysis}(b) demonstrates that waiting time costs $(C^{\tau})$ have a nonlinear and increasingly significant impact as the number of households grows. A 10\% increase in $C^{\tau}$ raises total costs by approximately 18.83\% to 22.82\%. Larger $|\mathcal{I}|$ values amplify this effect since more households imply higher EV arrival rates and more frequent charging events, thereby magnifying the consequences of driver idle time. Reductions in charger costs $(C^{\xi})$ show a more uniform effect across different household sizes (\Cref{fig:sen_analysis}(c)). An 80\% decrease in $C^{\xi}$ leads to roughly 14--15\% total cost savings, irrespective of $|\mathcal{I}|$. This uniformity suggests that cost-effective charging equipment investments benefit last-mile delivery operations across all scales. Improvements in charger powers $(P_k)$ yield more modest, yet scale-dependent benefits (\Cref{fig:sen_analysis}(d)). Doubling $P_k$ results in cost reductions of about 1.41\%, 2.14\%, and 3.17\% for $|\mathcal{I}|=50, 100,$ and $200$, respectively. Thus, while enhanced charging speeds consistently reduce total costs, their impact grows as the number of households (and, consequently, station utilization) increases. 

As a whole, these results highlight the interdependence between infrastructure costs, operational factors, and charging technology in shaping the economics of last-mile EV delivery. Although station costs dominate when household counts are low, waiting time and charger performance considerations become increasingly critical as the delivery network expands. Therefore, we believe that strategic decision-making should balance fixed infrastructure investments with enhancements in charging operations and technology, ensuring better performance across varying population sizes.

\section{Conclusion}
\label{sec:Conclusion}

Our paper presents a fresh approach to charging infrastructure planning by modeling the EV charging station location problem as a stochastic location model with congestion and immobile servers (SLCIS), which we refer to as the stochastic charger location and allocation problem (SCLA). We formulate the problem as a mixed-integer quadratic program and present two exact formulations, although the model building becomes computationally prohibitive for large scales using the exact methods. To address this issue, we have developed a Lagrangian-based dual decomposition framework, decomposing the problem into station-level subproblems, and we present two formulations for their solution. We use high-performance parallel computing across multiple partitions and solve the subproblems in parallel via a semi-relaxation and cutting-plane approach. We employed a three-step rounding heuristic to get integral solutions and used a subgradient method to improve lower bounds. Our extensive computational study on Chicago metropolitan data demonstrates that this framework consistently produces high-quality solutions for extremely large instances, whereas exact formulations cannot even be built. In practice, multiagency collaboration with combined depot and TAZ installations can yield substantial cost savings and reduced congestion. In conclusion, our paper yields valuable insights for strategic EV charging infrastructure planning. Our large-scale SCLA solution methodology proves highly effective in solving complex charging network design problems. It offers a tool for optimizers, urban planners, and policymakers navigating the transition to electric vehicle infrastructure. 

Future research directions may include incorporating time-varying demand patterns, exploring diverse geographic distributions of potential locations, and assessing the impact of emerging EV charger technologies on optimal network design. One could also refine our subgradient updates with more advanced techniques, such as proximal or bundle methods, to further accelerate convergence and tighten dual bounds. As we advance toward sustainable urban mobility, our paper provides a tool for informed decision-making in EV charging infrastructure development in the United States.

\section{Acknowledgment}

This material is based upon work supported by the U.S. Department of Energy, Office of Science, under contract number DE-AC02-06CH11357. 
This report and the work described were sponsored by the U.S. Department of Energy (DOE) Vehicle Technologies Office (VTO) under the Pathways to Affordable, Convenient, and Efficient (ACE) Regional Mobility, an initiative of the Energy Efficient Mobility Systems (EEMS) Program. Erin Boyd, a DOE Office of Critical Minerals and Energy Innovation (CMEI) manager, played an important role in establishing the project concept, advancing implementation, and providing guidance. The authors remain responsible for all findings and opinions presented in the paper. The findings are not suggestions for agencies to implement given the assumptions made in this study.

Computations for this research were performed on the Pennsylvania State University’s Institute for Computational and Data Sciences’ Roar supercomputer.

\clearpage







\clearpage

\section*{Appendix}\label{sec:appendix}

\section*{Relation between the $M/M/s_{jk}$ and $M/G/s_{jk}$ Systems}
\label{app:relation_mgs_mms}

We consider the Martin/Allen--Cunneen approximation \citep{allen2014probability} for the $M/G/s_{jk}$ system 
\citep[Eq.~(6.90)]{bolch2006queueing} as used in \citep{cokyasar2023additive}. 
In our notation, the expected time in the system is
\begin{equation}
\label{eq:MGk-AC-time-in-system}
\mathbb{W}(\rho_{jk}, s_{jk})_{(M/G/s)}
~\approx~
\frac{\mathbb{P}(\rho_{jk}, s_{jk})}{\mu_k s_{jk} (1-\rho_{jk})}
\cdot
\frac{1+c_{k}^2}{2 }
+ \frac{1}{\mu_k}
\end{equation}
where $c_{k}^2$ is the squared coefficient of variation of the charging-time distribution for charger type~$k$. 
Here, $\mu_k$ denotes a fixed service rate parameter for charger type~$k$ and is therefore omitted from the argument list of $\mathbb{W}$. For the $M/M/s_{jk}$ system, the Erlang--C formula provides the exact closed form expression
\begin{equation}
\label{eq:MMk-ErlangC-time-in-system}
\mathbb{W}(\rho_{jk}, s_{jk})_{(M/M/s)}
~=~
\frac{\mathbb{P}(\rho_{jk}, s_{jk})}{\mu_k s_{jk} (1-\rho_{jk})}
+ \frac{1}{\mu_k}
\end{equation}

Comparing \eqref{eq:MGk-AC-time-in-system} and \eqref{eq:MMk-ErlangC-time-in-system} shows that the
$M/G/s_{jk}$ expression differs from $M/M/s_{jk}$ only by a multiplicative factor on the queuing term:
\begin{equation}
\label{eq:MGk-as-factor-of-MMk}
\mathbb{W}(\rho_{jk}, s_{jk})_{(M/G/s)}
~\approx~
\frac{1+c_{k}^2}{2}
\Big(\mathbb{W}(\rho_{jk}, s_{jk})_{(M/M/s)} - \tfrac{1}{\mu_k}\Big)
+ \frac{1}{\mu_k}
\end{equation}
When service times are exponential ($c_{k}^2 = 1$), 
\eqref{eq:MGk-AC-time-in-system} reduces exactly to the Erlang--C formula 
\eqref{eq:MMk-ErlangC-time-in-system}. For any constant value of $c_{k}^2$, our exact method applies straightforwardly with minimal modification. When $c_{k}^2$ is treated as a decision variable, the exact method needs a more complicated treatment by decomposing the waiting time approximation into two functions and proving the convexity over one of them. We omit providing further details regarding this treatment as they can be found in \citep{cokyasar2023additive}.

\section*{Cut Instantiation under the $M/G/s_{jk}$ Approximation}
\label{app:mgs_cuts}

We recall that in the cutting-plane formulation of \Cref{sec:cutting_plane}, the nonlinear waiting-time function $\mathbb{W}(\rho_{jk}, s_{jk})$ 
was linearized by introducing supporting hyperplanes defined at candidate utilization points 
$\Tilde{\rho}_{jk} \in (0,1)$. 
For each $(j,k)$, let $A_{jk}$ and $B_{jk}$ denote the intercept and slope of the affine approximation 
to $\mathbb{W}^\nu(\rho_{jk}) = \mathbb{P}(\rho_{jk}, s_{jk}) / (1 - \rho_{jk})$ 
evaluated at $\Tilde{\rho}_{jk}$. 
The integer variable $z_{cjk}$ indicates whether $c$ chargers of type $k$ are installed 
at station $j$, and $x_{ijk}$ denotes the assignment of household $i$ to charger type $k$ 
at station $j$. 

From \eqref{eq:MGk-as-factor-of-MMk}, the expected time in system under the 
$M/G/s_{jk}$ model can be expressed as a scaled variant of the $M/M/s_{jk}$ formula:
\[
\mathbb{W}(\rho_{jk}, s_{jk})_{(M/G/s)}
~\approx~
\frac{1+c_{k}^2}{2}
\Big(\mathbb{W}(\rho_{jk}, s_{jk})_{(M/M/s)} - \tfrac{1}{\mu_k}\Big)
+ \frac{1}{\mu_k}
\]
Since $\mathbb{W}(\rho_{jk}, s_{jk})_{(M/M/s)}$ is strictly increasing and convex in $\rho_{jk}$ 
for fixed $s_{jk}$ (Property 4.1), 
$\mathbb{W}(\rho_{jk}, s_{jk})_{(M/G/s)}$ preserves these properties. 
Hence, the cutting-plane formulation extends directly to the $M/G/s_{jk}$ case 
with a multiplicative factor to the waiting time component. Hence, the resulting valid cut is given by
\begin{equation}
\label{W_cut_MGk}
\begin{aligned}
    W_{jk} 
    \geq &
    \ \frac{1+c_{k}^2}{2}
    \left(
    \frac{A_{jk}z_{c jk}}{\mu_k c}
    + 
    \frac{B_{jk} \sum_{i\in \mathcal{I}_j} \lambda_i (z_{c jk} + x_{ijk} - 1)}{\mu_k^2 c^2}
    \right)
    + 
    \frac{z_{c jk}}{\mu_k},
    \\
    &\forall c \in \mathcal{S} \setminus \{0\}, \;
    j \in \mathcal{J}, \;
    k \in \mathcal{K}
\end{aligned}
\end{equation}

When $c_{k}^2 = 1$, the constraint \eqref{W_cut_MGk} 
reduces exactly to the original $M/M/s_{jk}$ cut \eqref{W_cut_exact}. For non-exponential charging-time distributions ($c_{k}^2 \neq 1$), 
the factor $(1+c_{k}^2)/2$ provides a valid lower bound for the waiting time during the branch-and-bound process.

\section*{Impact of Under-Dispersion on Congestion Estimates}
\label{app:gi_gs_analysis}

Delivery demand counts at the household or cluster level may exhibit \textit{under-dispersion} relative to a Poisson process, particularly in dense urban 
areas where multiple delivery points (e.g., apartment complexes or multi-tenant buildings) 
receive at most one delivery per day from a carrier. 
To understand how such under-dispersion affects our household-level Poisson assumption (\Cref{sec:modeling_household_demand}), 
we compare the $M/M/s$ system with the more general $GI/G/s$ system using 
the Allen--Cunneen approximation formula~\citep{allen2014probability,bolch2006queueing}, 
which extends the $M/M/s$ model to account for renewal interarrival processes and 
service-time variability.

Let $c_{a,jk}^2$ and $c_{k}^2$ denote the squared coefficients of variation of the interarrival and service times, respectively, for charger type $k$ at station $j$. Then the utilization of that station–charger pair is
\[
\rho_{jk} = \frac{\bar{\lambda}_{jk}}{s_{jk}\mu_k}
\]
where $\bar{\lambda}_{jk} = \sum_{i\in \mathcal{I}_j} \lambda_i x_{ijk}$ and $\lambda_i = \gamma_i \pi_i$ denotes the charging rate for the vehicle delivering to household $i$ (see \Cref{sec:modeling_household_demand}).

Under the Allen--Cunneen approximation \citep{allen2014probability}, the expected total time in system at $(j,k)$ is
\begin{equation}\label{eq:ac_general}
\mathbb{W}(GI/G/s)
~\approx~
\frac{\rho_{jk}/\mu_k}{1-\rho_{jk}}
\cdot
\frac{c_{a,jk}^2 + c_{k}^2}{2}
+ \frac{1}{\mu_k}.
\end{equation}

In the $M/M/s$ case ($c_{a,jk}^2 = c_{k}^2 = 1$), we obtain
\begin{equation}\label{eq:mm_special}
\mathbb{W}(M/M/s)
~\approx~
\frac{\rho_{jk}/\mu_k}{1-\rho_{jk}}
+ \frac{1}{\mu_k}.
\end{equation}

Comparing \eqref{eq:ac_general} and \eqref{eq:mm_special} shows that the congestion (queuing) component of $\mathbb{W}$ in the $GI/G/s$ case is scaled by a factor
\[
\frac{c_{a,jk}^2 + c_{k}^2}{2}
\]
relative to the $M/M/s$ case. Since our modeling approach assumes exponential charging service ($c_{k}^2 = 1$), we have
\[
\mathbb{W}(GI/G/s)
~\approx~
\frac{c_{a,jk}^2 + 1}{2}\;
\frac{\rho_{jk}/\mu_k}{1-\rho_{jk}}
+ \frac{1}{\mu_k}.
\]
Hence, when demand counts are under-dispersed ($c_{a,jk}^2 < 1$),
\[
\mathbb{W}(GI/G/s) \le \mathbb{W}(M/M/s),
\]
showing that our $M/M/s_{jk}$ formulation yields conservative congestion estimates. 

Furthermore, as illustrated by our dataset for the Chicago metropolitan area (see \Cref{fig:household_distribution}), 
households are highly nonuniformly distributed—dense, clustered zones in the downtown Chicago Loop area coexist with more dispersed suburban regions like Naperville. Such heterogeneity implies that some local demand streams may indeed be under-dispersed, while others may be approximately Poisson. Therefore, over a strategic planning horizon, modeling household demand as independent Poisson processes provides a simple yet conservative approximation that captures the aggregate variability observed across heterogeneous urban regions.


\begin{algorithm}[H]
\small
\caption{Construction of $\mathcal{J}_i$ and $\mathcal{I}_j$}
\label{alg:construct_Ji_Ij}
\begin{algorithmic}[1]
\REQUIRE 
  $\mathcal{I}$: set of households; 
  $\mathcal{J}$: set of candidate stations; 
  $k_c$: number of nearest stations to consider; 
  $\mathbf{x}(i)$: \((\text{latitude}, \text{longitude})\) for each household $i \in \mathcal{I}$; 
  $\mathbf{x}(j)$: \((\text{latitude}, \text{longitude})\) for each station $j \in \mathcal{J}$.
\ENSURE 
  $\mathcal{J}_i$: up to $k_c$ closest stations for each $i$; 
  $\mathcal{I}_j$: set of households for which $j$ is among their $k_c$ closest.
\STATE Build a 2D $k$-d tree from the station coordinates $\mathbf{x}(j)$, $j \in \mathcal{J}$.
\FOR{$i \in \mathcal{I}$}
    \STATE Query the $k_c$ nearest stations to $\mathbf{x}(i)$ using the $k$-d tree.
    \STATE $\mathcal{J}_i(i) \gets$ those $k_c$ stations
\ENDFOR
\FOR{$j \in \mathcal{J}$}
    \STATE $\mathcal{I}_j(j) \gets \{\,i \in \mathcal{I}\,\mid\,j \in \mathcal{J}_i(i)\}$
\ENDFOR
\RETURN $\mathcal{J}_i, \mathcal{I}_j$
\end{algorithmic}
\end{algorithm}


\newcommand{\Assign}{\leftarrow}

\begin{algorithm}[H]
\small
\caption{Station Opening}
\label{alg:adaptive_station_opening}
\begin{algorithmic}[1]
\REQUIRE $ \mathcal{J} $: set of potential station locations, $ \mathcal{I} $: set of all households, $ y_j^t $ for all $ j \in \mathcal{J} $: station open states at iteration $ t $, $ \mathcal{I}_j $ for all $ j \in \mathcal{J} $: set of households coverable by station $ j $

\ENSURE $ y_j' $: updated station open states, $ \mathcal{U}^c $: set of uncovered households after updating station states

\STATE Initialize $ y_j' \Assign y_j^t $ for all $ j \in \mathcal{J} $ \hfill $ \triangleright $ Initialize station states
\STATE $ O \Assign \{ j \in \mathcal{J} \mid y_j^t = 1 \} $ \hfill $ \triangleright $ Set of open stations
\STATE $ C \Assign \bigcup_{j \in O} \mathcal{I}_j $ \hfill $ \triangleright $ Set of covered households
\STATE $ \mathcal{U}^c \Assign \mathcal{I} \setminus C $ \hfill $ \triangleright $ Set of uncovered households

\IF{$ \mathcal{U}^c = \emptyset $}
    \RETURN $ y_j' $, $ \emptyset $ \hfill $ \triangleright $ All households are covered; no new openings required
\ENDIF

\WHILE{$ \mathcal{U}^c \neq \emptyset $ \AND $ \mathcal{J} \setminus O \neq \emptyset $}
    \STATE For each $ j \in \mathcal{J} \setminus O $, compute $ w_j \Assign |\mathcal{I}_j \cap \mathcal{U}^c| $ \hfill $ \triangleright $ Compute weights
    \STATE $ \mathcal{J}_{\text{cand}} \Assign \{ j \in \mathcal{J} \setminus O \mid w_j > 0 \} $ \hfill $ \triangleright $ Candidate stations
    \IF{$ \mathcal{J}_{\text{cand}} = \emptyset $}
        \STATE \textbf{break} \hfill $ \triangleright $ No more stations can cover remaining households
    \ENDIF
    \STATE $ W \Assign \sum_{j \in \mathcal{J}_{\text{cand}}} w_j $ \hfill $ \triangleright $ Total weight for probability calculation
    \STATE For each $ j \in \mathcal{J}_{\text{cand}} $, compute $ p_j \Assign w_j / W $ \hfill $ \triangleright $ Calculate selection probabilities
    \STATE Randomly select station $ j^* $ from $ \mathcal{J}_{\text{cand}} $ based on probabilities $ p_j $ \hfill $ \triangleright $ Select station probabilistically
    \STATE Set $ y_{j^*}' \Assign 1 $ \hfill $ \triangleright $ Open station $ j^* $
    \STATE Update $ O \Assign O \cup \{ j^* \} $ \hfill $ \triangleright $ Update set of open stations
    \STATE Update $ \mathcal{U}^c \Assign \mathcal{U}^c \setminus \mathcal{I}_{j^*} $ \hfill $ \triangleright $ Remove covered households
\ENDWHILE

\IF{$ \mathcal{U}^c = \emptyset $}
    \RETURN $ y_j' $, $ \emptyset $ \hfill $ \triangleright $ Coverage achieved for all households
\ENDIF

\end{algorithmic}
\end{algorithm}


\newcommand{\ceildiv}[2]{\left\lceil \dfrac{#1}{#2} \right\rceil}

\begin{algorithm}[H]
\small
\caption{Household Assignment Iteration}
\label{alg:household_iteration}
\begin{algorithmic}[1]
\REQUIRE Sets $\mathcal{I}, \mathcal{J}, \mathcal{K}$; parameters $x_{ijk}^t, y_j^t, \lambda_i, \mu_k, \epsilon, \bar{S}_{jk}, EW, \mathcal{J}_i$; optional flag \texttt{probabilistic} $\in \{\texttt{true}, \texttt{false}\}$
\ENSURE Updated variables $ x_{ijk}' $, $ s_{jk}' $, $ W_{jk}' $, $ y_j' $; set of unassigned households $\mathcal{U}^a$

\STATE Initialize $x_{ijk}' \gets 0$ for all $i,j,k$
\STATE $s_{jk}' \gets s_{jk}^t$, $W_{jk}' \gets W_{jk}^t$

\STATE Compute $D_{jk} \gets \sum_{i' \in \mathcal{I}} \lambda_{i'} x_{i'jk}^t$ for all $j,k$

\STATE Initialize $\mathcal{U}^a \gets \emptyset$
\STATE Sort $\mathcal{I}$ in descending order of $\lambda_i$

\FOR{each $i \in \mathcal{I}$}
    \STATE $V_i \gets \{ j \in \mathcal{J}_i : y_j^t = 1 \}$
    \IF{$V_i = \emptyset$}
        \STATE $\mathcal{U}^a \gets \mathcal{U}^a \cup \{ i \}$
        \STATE \textbf{continue}
    \ENDIF

    \STATE $\mathcal{C}_i \gets \{ (j,k) : j \in V_i, k \in \mathcal{K} \}$
    \IF{\texttt{probabilistic} = \texttt{true}}
        \STATE $\mathcal{C}_i' \gets \{(j,k)\in\mathcal{C}_i : x_{ijk}^t > 0\}$
        \STATE Normalize $p_{ijk} = \dfrac{x_{ijk}^t}{\sum_{(j',k')\in\mathcal{C}_i'} x_{ij'k'}^t}$ for $(j,k)\in\mathcal{C}_i'$
        \STATE Set $x_{ijk}' \gets p_{ijk}$ for all $(j,k)\in\mathcal{C}_i'$
        \STATE \textbf{continue} \hfill $ \triangleright $ Skip deterministic feasibility checks
    \ENDIF
    \STATE Sort $\mathcal{C}_i$ by (1) descending $x_{ijk}^t$, (2) ascending $C^{\xi}_k$, (3) ascending $W_{jk}$ \hfill $ \triangleright $ Break ties at random

    \STATE \textit{assigned} $\gets$ \textbf{false}

    \FOR{each $(j,k) \in \mathcal{C}_i$}
        \STATE \textit{assigned} $\gets$ \textsc{CheckAndAssign}($i,j,k,\lambda_i,\mu_k,\epsilon,\bar{S}_{jk},EW,D_{jk},s_{jk}',W_{jk}'$)
        \IF{\textit{assigned} = \textbf{true}}
            \STATE \textbf{break} \hfill $ \triangleright $ Move to next household
        \ENDIF
    \ENDFOR

    \IF{\textit{assigned} = \textbf{false}}
        \STATE $\mathcal{U}^a \gets \mathcal{U}^a \cup \{ i \}$
    \ENDIF
\ENDFOR

\FOR{each $j \in \mathcal{J}$}
    \STATE $y_j' \gets 1$ if $\exists\, i,k: x_{ijk}'=1$, else $y_j' \gets 0$
\ENDFOR

\RETURN $x_{ijk}', s_{jk}', W_{jk}', y_j', \mathcal{U}^a$
\end{algorithmic}
\end{algorithm}

\begin{algorithm}[H]
\small
\caption{Check and Assign Feasibility}
\label{alg:check_and_assign}
\begin{algorithmic}[1]
\REQUIRE Household $i$, station $j$, charger type $k$, parameters $\lambda_i, \mu_k, \epsilon, \bar{S}_{jk}, EW$, and current values $D_{jk}, s_{jk}', W_{jk}'$; optional weight $p_{ijk}$ (default $1$)
\ENSURE Returns \textbf{true} if assignment is successful, else \textbf{false}. Updates $x_{ijk}', s_{jk}', W_{jk}'$ on success.

\STATE $D_{jk}^\text{c} \gets D_{jk} + p_{ijk}\lambda_i$ \hfill $ \triangleright $ $p_{ijk}=1$ in deterministic mode
\STATE $s_{jk}^\text{req} \gets \lceil \frac{D_{jk}^\text{c}}{\mu_k (1-\epsilon)} \rceil$

\IF{$s_{jk}^\text{req} > \bar{S}_{jk}$}
    \STATE \textbf{return false} \hfill $ \triangleright $ Capacity constraint violated
\ENDIF

\STATE $s_{jk}^\text{c} \gets \max(s_{jk}', s_{jk}^\text{req})$

\WHILE{$s_{jk}^\text{c} \leq \bar{S}_{jk}$}
    \STATE $\rho_{jk}^\text{c} \gets \frac{D_{jk}^\text{c}}{s_{jk}^\text{c}\mu_k}$
    \IF{$\rho_{jk}^\text{c} \geq 1$}
        \STATE \textbf{break} \hfill $ \triangleright $ Cannot assign with current $s_{jk}^\text{c}$
    \ENDIF

    \STATE $W_{jk}^\text{c} \gets \mathbb{W}(\rho_{jk}^\text{c}, s_{jk}^\text{c}, \mu_k)$
    \IF{$W_{jk}^\text{c} \leq EW + \frac{1}{\mu_k}$}
        \STATE $D_{jk} \gets D_{jk}^\text{c}$
        \STATE $s_{jk}' \gets s_{jk}^\text{c}$
        \STATE $x_{ijk}' \gets p_{ijk}$
        \STATE $W_{jk}' \gets W_{jk}^\text{c}$
        \STATE \textbf{return true}
    \ELSE
        \STATE $s_{jk}^\text{c} \gets s_{jk}^\text{c} + 1$
    \ENDIF
\ENDWHILE

\STATE \textbf{return false} \hfill $ \triangleright $ No feasible assignment found
\end{algorithmic}
\end{algorithm}


\begin{algorithm}[H]
\small
\caption{Calculate Additional Stations for Overload}
\label{alg:additional_stations_overload}
\begin{algorithmic}[1]
\REQUIRE Unassigned households $ \mathcal{U}^a $; current depot states $ y_j' $; depot accessibility $ \mathcal{J}_i $; household assignments $ \mathcal{I}_j $
\ENSURE Updated depot states $ y_j' $

\STATE $ \mathcal{C} \Assign \bigcup_{i \in \mathcal{U}^a} \mathcal{J}_i $ \hfill $ \triangleright $ Candidate depots accessible to unassigned households
\STATE $ \mathcal{C} \Assign \mathcal{C} \setminus \{ j \in \mathcal{J} \mid y_j' = 1 \} $ \hfill $ \triangleright $ Exclude already open depots
\IF{$ \mathcal{C} = \emptyset $}
    \RETURN $ y_j' $ \hfill $ \triangleright $ No additional depots can be opened
\ENDIF
\STATE For each $ j \in \mathcal{C} $, compute $ w_j \Assign |\mathcal{I}_j \cap \mathcal{U}^a| $
\STATE $ j^* \Assign \arg\max_{j \in \mathcal{C}} w_j $ \hfill $ \triangleright $ Depot covering most unassigned households
\STATE $ y_{j^*}' \Assign 1 $ \hfill $ \triangleright $ Open depot $ j^* $
\RETURN $ y_j' $
\end{algorithmic}
\end{algorithm}

\begin{algorithm}[H]
\small
\caption{Primal Feasible Solution Heuristic}
\label{alg:primal_feasible}
\begin{algorithmic}[1]
\REQUIRE Sets $ \mathcal{I}, \mathcal{J}, \mathcal{K} $; variables $ x_{ijk}^t $, $ y_j^t $; parameters $ \lambda_i $, $ \mu_k $, $ \epsilon $, $ \bar{S}_{jk} $, $ EW $, $ \mathcal{J}_i $, $ \mathcal{I}_j$; optional flag \texttt{probabilistic} $\in \{\texttt{true},\texttt{false}\}$
\ENSURE Feasible assignment $ x_{ijk}' $, updated depot states $ y_j' $, charger allocations $ s_{jk}' $, waiting times $ W_{jk}' $

\STATE Call \Cref{alg:adaptive_station_opening} with $ \mathcal{J}, \mathcal{I}, y_j^t, \mathcal{I}_j $ to obtain $ y_j' $ and $ \mathcal{U}^c $
\IF{$ \mathcal{U}^c = \emptyset $} 
    \STATE Initialize $ \mathcal{U}^a \gets \emptyset $
    \WHILE{$ \mathcal{U}^a \neq \emptyset $} 
        \STATE $(x_{ijk}', s_{jk}', W_{jk}', y_j', \mathcal{U}^a) \gets$
        \STATE \quad \textsc{HouseholdAssignmentIteration}($\mathcal{I}, \mathcal{J}, \mathcal{K}, x_{ijk}^t, y_j', \lambda_i, \mu_k, \epsilon, \bar{S}_{jk}, EW, \mathcal{J}_i$, \texttt{probabilistic})
        \IF{$ \mathcal{U}^a \neq \emptyset $}
            \STATE Call \Cref{alg:additional_stations_overload} with $ \mathcal{U}^a $, $ y_j' $, $ \mathcal{J}_i $ to update $ y_j' $
        \ENDIF
    \ENDWHILE
    \RETURN $ x_{ijk}' $, $ y_j' $, $ s_{jk}' $, $ W_{jk}' $
\ELSE
    \RETURN Failure
\ENDIF
\end{algorithmic}
\end{algorithm}

\clearpage
\bibliographystyle{elsarticle-harv} 
\bibliography{TR_PartB_revision1/references}

@article{kelley1960cutting,
  doi = {10.1137/0108053},
  title={The cutting-plane method for solving convex programs},
  author={Kelley, Jr, James E},
  journal={Journal of the Society for Industrial and Applied Mathematics},
  volume={8},
  number={4},
  pages={703--712},
  year={1960},
  publisher={SIAM}
}

@book{stewart2009,
 doi = {10.1515/9781400832811},
 author = {William J. Stewart},
 publisher = {Princeton University Press},
 title = {Probability, {Markov} Chains, Queues, and Simulation: The Mathematical Basis of Performance Modeling},
 year = {2009}
}

@article{berman2007multiple,
  title={The multiple server location problem},
  author={Berman, Oded and Drezner, Zvi},
  journal={Journal of the Operational Research Society},
  volume={58},
  number={1},
  pages={91--99},
  year={2007},
  publisher={Taylor \& Francis},
  doi = {10.1057/palgrave.jors.2602126}
}

@incollection{bermanandkrass,
  doi = {10.1007/978-3-030-32177-2_17},
  booktitle={Location Science},
  author = {Berman, Oded and Krass, Dmitry},
  title = {Stochastic Location Models with Congestion},
  editor={Laporte, Gilbert and Nickel, Stefan and Saldanha da Gama, Francisco},
  year={2019},
  edition = {$2^{nd}$},
  chapter = {17},
  pages = {477--531},
  publisher={Springer},
  address = {Switzerland},
  isbn = {978-3-030-32176-5}
}

@article{aboolian2008location,
  doi = {10.1080/07408170701411385},
  title={Location and allocation of service units on a congested network},
  author={Aboolian, Robert and Berman, Oded and Drezner, Zvi},
  journal={IIE Transactions},
  volume={40},
  number={4},
  pages={422--433},
  year={2008},
  publisher={Taylor \& Francis}
}

@article{elhedhli2006service,
  title={Service system design with immobile servers, stochastic demand, and congestion},
  author={Elhedhli, Samir},
  journal={Manufacturing \& Service Operations Management},
  volume={8},
  number={1},
  pages={92--97},
  year={2006},
  publisher={INFORMS},
  doi = {10.1287/msom.1050.0094}
}

@article{amiri1997solution,
  title={Solution procedures for the service system design problem},
  author={Amiri, Ali},
  journal={Computers \& Operations Research},
  volume={24},
  number={1},
  pages={49--60},
  year={1997},
  publisher={Elsevier},
  doi = {10.1016/S0305-0548(96)00022-6}
}

@article{wang2002algorithms,
  doi = {10.1023/A:1020961732667},
  title={Algorithms for a facility location problem with stochastic customer demand and immobile servers},
  author={Wang, Qian and Batta, Rajan and Rump, Christopher M},
  journal={Annals of Operations Research},
  volume={111},
  pages={17--34},
  year={2002},
  publisher={Springer}
}

@article{castillo2009social,
  title={Social optimal location of facilities with fixed servers, stochastic demand, and congestion},
  author={Castillo, Ignacio and Ingolfsson, Armann and Sim, Thaddeus},
  journal={Production and Operations Management},
  volume={18},
  number={6},
  pages={721--736},
  year={2009},
  publisher={SAGE Publications Sage CA: Los Angeles, CA},
  doi={10.1111/j.1937-5956.2009.01034.x}
}

@article{halfin1981heavy,
  title={Heavy-traffic limits for queues with many exponential servers},
  author={Halfin, Shlomo and Whitt, Ward},
  journal={Operations Research},
  volume={29},
  number={3},
  pages={567--588},
  year={1981},
  publisher={INFORMS},
  doi = {10.1287/opre.29.3.567}
}

@article{borst2004dimensioning,
  title={Dimensioning large call centers},
  author={Borst, Sem and Mandelbaum, Avi and Reiman, Martin I},
  journal={Operations Research},
  volume={52},
  number={1},
  pages={17--34},
  year={2004},
  publisher={INFORMS},
  doi = {10.1287/opre.1030.0081}
}

@article{vidyarthi2014efficient,
  title={Efficient solution of a class of location-allocation problems with stochastic demand and congestion},
  author={Vidyarthi, Navneet and Jayaswal, Sachin},
  journal={Computers \& Operations Research},
  volume={48},
  pages={20--30},
  year={2014},
  publisher={Elsevier},
  doi = {10.1016/j.cor.2014.02.014}
}

@article{ahmadi2018location,
  doi = {10.48550/arXiv.1809.00080},
  title={Location and capacity planning of facilities with general service-time distributions using conic optimization},
  author={Ahmadi-Javid, Amir and Berman, Oded and Hoseinpour, Pooya},
  journal={arXiv preprint arXiv:1809.00080},
  year={2018}
}

@article{ahmadi2020linear,
  title={Linear formulations and valid inequalities for a classic location problem with congestion: a robust optimization application},
  author={Ahmadi-Javid, Amir and Ramshe, Nasrin},
  journal={Optimization Letters},
  volume={14},
  number={5},
  pages={1265--1285},
  year={2020},
  publisher={Springer},
  doi={10.1007/s11590-019-01419-8}
}

@article{ahmadi2022convexification,
  title={Convexification of queueing formulas by mixed-integer second-order cone programming: An application to a discrete location problem with congestion},
  author={Ahmadi-Javid, Amir and Hoseinpour, Pooya},
  journal={INFORMS Journal on Computing},
  volume={34},
  number={5},
  pages={2621--2633},
  year={2022},
  publisher={INFORMS},
  doi = {10.1287/ijoc.2021.1125}
}

@ARTICLE{Syam_2008,title={A multiple server location-allocation model for service system design},year={2008},author={Siddhartha S. Syam},mag_id={1981972625},journal={Computers \& Operations Research}, doi = {10.1016/j.cor.2006.10.019}}

@ARTICLE{Elhedhli_2018,title={Service system design with immobile servers, stochastic demand and concave-cost capacity selection},year={2018},author={Samir Elhedhli and Samir Elhedhli and Yan Wang and Yan Wang and Ahmed Saif and Ahmed Saif},doi={10.1016/j.cor.2018.01.019},pmid={null},pmcid={null},mag_id={2787928828},journal={Computers \& Operations Research}}

@ARTICLE{Goez_2017,title={Second-Order Cone Optimization Formulations for Service System Design Problems with Congestion},year={2017},author={Julio C. Góez and Julio C. Góez and Miguel F. Anjos and Miguel F. Anjos},doi={10.1007/978-3-030-12119-8_5},pmid={null},pmcid={null},mag_id={2912123179},journal={null}}

@article{aboolian2022efficient,
  url = {https://www.csupom.com/uploads/1/1/4/8/114895679/n20p1formatted.pdf},
  title={An Efficient Approach for Service System Design with Immobile Servers, Stochastic Demand, Congestion, and Consumer Choice},
  author={Aboolian, Robert and Elhedhli, Samir and Karimi, Majid},
  journal={Journal of Supply Chain and Operations Management},
  volume={20},
  number={1},
  pages={1},
  year={2022}
}

@ARTICLE{Vidyarthi_2015,title={The impact of directed choice on the design of preventive healthcare facility network under congestion.},year={2015},author={Navneet Vidyarthi and Navneet Vidyarthi and Onur Kuzgunkaya and Onur Kuzgunkaya},doi={10.1007/s10729-014-9274-2},pmid={24879402},pmcid={null},mag_id={2092597604},journal={Health Care Management Science}}

@ARTICLE{Etebari_2019,title={A column generation algorithm for the choice-based congested location-pricing problem},year={2019},author={Farhad Etebari},doi={10.1016/j.cie.2019.03.023},pmid={null},pmcid={null},mag_id={2920956709},journal={Computers \& Industrial Engineering}}

@ARTICLE{Aboolian_2009,title={The multiple server center location problem},year={2009},author={Robert Aboolian and Robert Aboolian and Oded Berman and Oded Berman and Zvi Drezner and Zvi Drezner},doi={10.1007/s10479-008-0341-2},pmid={null},pmcid={null},mag_id={2058073949},journal={Annals of Operations Research}}

@article{cokyasar2023additive,
  title={Additive manufacturing capacity allocation problem over a network},
  author={Cokyasar, Taner and Jin, Mingzhou},
  journal={IISE Transactions},
  volume={55},
  number={8},
  pages={807--820},
  year={2023},
  publisher={Taylor \& Francis},
  doi = {10.1080/24725854.2022.2120222}
}

@article{barahona2000volume,
  doi = {10.1007/s101070050002},
  title={The volume algorithm: Producing primal solutions with a subgradient method},
  author={Barahona, Francisco and Anbil, Ranga},
  journal={Mathematical Programming},
  volume={87},
  pages={385--399},
  year={2000},
  publisher={Springer}
}

@article{guignard2003lagrangean,
  doi = {10.1007/bf02579036},
  title={Lagrangean relaxation},
  author={Guignard, Monique},
  journal={Top},
  volume={11},
  pages={151--200},
  year={2003},
  publisher={Springer}
}

@book{glover1998tabu,
  doi = {10.1007/978-1-4615-6089-0},
  title={Tabu Search},
  author={Glover, Fred and Laguna, Manuel},
  year={1998},
  publisher={Springer}
}

@article{mccormick1976computability,
  doi = {10.1007/BF01580665},
  title={Computability of global solutions to factorable nonconvex programs: {Part I}—Convex underestimating problems},
  author={McCormick, Garth P},
  journal={Mathematical Programming},
  volume={10},
  number={1},
  pages={147--175},
  year={1976},
  publisher={Springer}
}

@misc{IEA2024,
  author       = {{International Energy Agency}},
  title        = {Global {EV} Outlook 2024},
  year         = {2024},
  publisher    = {IEA},
  address      = {Paris},
  url          = {https://www.iea.org/reports/global-ev-outlook-2024},
  note         = {Licence: CC BY 4.0}
}

@article{grassmann1983convexity,
  title={The convexity of the mean queue size of the {M/M/c} queue with respect to the traffic intensity},
  author={Grassmann, W},
  journal={Journal of Applied Probability},
  volume={20},
  number={4},
  pages={916--919},
  year={1983},
  publisher={Cambridge University Press},
  doi = {10.2307/3213605}
}

@article{lee1983note,
  title={A note on the convexity of performance measures of {M/M/c} queueing systems},
  author={Lee, Hau Leung and Cohen, Morris A},
  journal={Journal of Applied Probability},
  volume={20},
  number={4},
  pages={920--923},
  year={1983},
  publisher={Cambridge University Press},
  doi = {10.2307/3213606}
}

@article{barahona2005near,
  doi = {10.1016/j.disopt.2003.03.001},
  title={Near-optimal solutions to large-scale facility location problems},
  author={Barahona, Francisco and Chudak, Fabi{\'a}n A},
  journal={Discrete Optimization},
  volume={2},
  number={1},
  pages={35--50},
  year={2005},
  publisher={Elsevier}
}

@article{auld2016polaris,
  title={{POLARIS}: Agent-based modeling framework development and implementation for integrated travel demand and network and operations simulations},
  author={Auld, Joshua and Hope, Michael and Ley, Hubert and Sokolov, Vadim and Xu, Bo and Zhang, Kuilin},
  journal={Transportation Research Part C: Emerging Technologies},
  volume={64},
  pages={101--116},
  year={2016},
  publisher={Elsevier},
  doi = {10.1016/j.trc.2015.07.017}
}

@techreport{wood20232030,
  doi = {10.15483/1996812},
  title={The 2030 national charging network: Estimating {US} light-duty demand for electric vehicle charging infrastructure},
  author={Wood, Eric and Borlaug, Brennan and Moniot, Matt and Lee, Dong-Yeon DY and Ge, Yanbo and Yang, Fan and Liu, Zhaocai},
  year={2023},
  institution={National Renewable Energy Laboratory, Golden, CO (United States)}
}

@article{DolsakPrakash2021,
  author = {Dolsak, Nives and Prakash, Aseem},
  title = {The Lack Of {EV} Charging Stations Could Limit {EV} Growth},
  journal = {Forbes},
  year = {2021},
  month = {May},
  day = {5},
  url = {https://www.forbes.com/sites/prakashdolsak/2021/05/05/the-lack-of-ev-charging-stations-could-limit-ev-growth/},
  note = {Accessed on July 31, 2024}
}

@article{davatgari2024electric,
  title={Electric vehicle supply equipment location and capacity allocation for fixed-route networks},
  author={Davatgari, Amir and Cokyasar, Taner and Subramanyam, Anirudh and Larson, Jeffrey and Mohammadian, Abolfazl Kouros},
  journal={European Journal of Operational Research},
  volume={317},
  number={3},
  pages={953--966},
  year={2024},
  publisher={Elsevier},
  doi = {10.1016/j.ejor.2024.04.022}
}

@article{hale2003location,
  title={Location science research: A review},
  author={Hale, Trevor S and Moberg, Christopher R},
  journal={Annals of Operations Research},
  volume={123},
  pages={21--35},
  year={2003},
  publisher={Springer},
  doi = {10.1023/A:1026110926707}
}

@mastersthesis{davatgari2021location,
  url = {https://docs.lib.purdue.edu/dissertations/AAI30504864/},
  title={Location planning for electric charging stations and wireless facilities in the era of autonomous vehicle operations},
  author={Davatgari, Amir},
  year={2021},
  school={Purdue University}
}

@misc{ElectrifyAmerica2024,
  title        = {Pricing and Plans for {EV} Charging},
  howpublished = {\url{https://www.electrifyamerica.com/pricing/}},
  note         = {Accessed: 2024-08-13},
  year         = 2024,
  author       = {{Electrify America}}
}

@misc{LightningEMotors2022,
  title        = {Lightning Electric Class 6 LowCab Forward},
  howpublished = {\url{https://lightningemotors.com/}},
  note         = {Accessed on: 03 August 2024},
  year         = 2022,
  author       = {{Lightning eMotors}}
}

@techreport{SmithCastellano2015,
  author      = {Michael Smith and Joseph Castellano},
  title       = {Costs Associated with Non-Residential Electric Vehicle Supply Equipment},
  institution = {U.S. Department of Energy},
  year        = 2015,
  pages       = {1--43},
  month       = nov,
  url         = {https://afdc.energy.gov/files/u/publication/evse_cost_report_2015.pdf},
  note        = {Accessed on: 13 August 2024}
}

@misc{Williams2020,
  author      = {K. Williams},
  title       = {How Many Hours Can a Truck Driver Drive?},
  year        = 2020,
  howpublished= {\url{https://www.cdljobs.com/news-notes/news/how-many-hours-can-a-truck-driver-drive}},
  note        = {Accessed on: 13 August 2024}
}

@misc{RivianEDV2024,
  author       = {{Rivian EDV}},
  title        = {Rivian {EDV}},
  year         = 2024,
  howpublished = {\url{https://en.wikipedia.org/wiki/Rivian_EDV}},
  note         = {Specific section on battery capacity referenced, accessed on: 13 August 2024}
}

@misc{UPSDriverSalary2024,
  author       = {{Upper Inc.}},
  title        = {{UPS} Driver Salary},
  year         = 2024,
  howpublished = {\url{https://www.upperinc.com/blog/ups-driver-salary}},
  note         = {Specific section on average hourly wage after four years referenced, accessed on: 13 August 2024}
}

@misc{FedExRevenue2024,
  author       = {{FedEx}},
  title        = {{FedEx} Company Structure},
  year         = 2024,
  howpublished = {\url{https://www.fedex.com/en-us/about/company-structure.html}},
  note         = {Data on annual revenue, accessed on: 13 August 2024}
}

@misc{FedExOperationalDays2024,
  author       = {{FedEx}},
  title        = {{FedEx} Company Structure},
  year         = 2024,
  howpublished = {\url{https://www.fedex.com/en-us/about/company-structure.html}},
  note         = {Data on operational days per year, accessed on: 13 August 2024}
}

@misc{FedExFleetSize2024,
  author       = {{FedEx}},
  title        = {{FedEx} Company Structure},
  year         = 2024,
  howpublished = {\url{https://www.fedex.com/en-us/about/company-structure.html}},
  note         = {Data on fleet size, accessed on: 13 August 2024}
}

@misc{FedExDailyPackagesDelivered2024,
  author       = {{FedEx}},
  title        = {{FedEx} Company Structure},
  year         = 2024,
  howpublished = {\url{https://www.fedex.com/en-us/about/company-structure.html}},
  note         = {Data on daily packages delivered, accessed on: 13 August 2024}
}

@techreport{bennett2022estimating,
  doi = {10.2172/1894645},
  title={Estimating the Breakeven Cost of Delivered Electricity to Charge Class 8 Electric Tractors},
  author={Bennett, Jesse and Mishra, Partha and Miller, Eric and Borlaug, Brennan and Meintz, Andrew and Birky, Alicia},
  year={2022},
  institution={National Renewable Energy Lab.(NREL), Golden, CO (United States)}
}

@article{liu2017locating,
  title={Locating multiple types of charging facilities for battery electric vehicles},
  author={Liu, Haoxiang and Wang, David ZW},
  journal={Transportation Research Part B: Methodological},
  volume={103},
  pages={30--55},
  year={2017},
  publisher={Elsevier},
  doi = {10.1016/j.trb.2017.01.005}
}

@article{yilmaz2012review,
  title={Review of battery charger topologies, charging power levels, and infrastructure for plug-in electric and hybrid vehicles},
  author={Yilmaz, Murat and Krein, Philip T},
  journal={IEEE Transactions on Power Electronics},
  volume={28},
  number={5},
  pages={2151--2169},
  year={2012},
  publisher={IEEE},
  doi = {10.1109/TPEL.2012.2212917}
}

@Misc{Gurobi2024,
  author = {{Gurobi Optimization, LLC}},
  title = {Gurobi Optimizer Reference Manual},
  year = {2024},
  howpublished = {Online},
  url = {https://www.gurobi.com}
}

@article{dalcin2008mpi,
  doi = {10.1016/j.jpdc.2007.09.005},
  title={{MPI} for {Python}: Performance improvements and {MPI-2} extensions},
  author={Dalc{\'\i}n, Lisandro and Paz, Rodrigo and Storti, Mario and D’El{\'\i}a, Jorge},
  journal={Journal of Parallel and Distributed Computing},
  volume={68},
  number={5},
  pages={655--662},
  year={2008},
  publisher={Elsevier}
}

@article{kchaou2021charging,
  title={Charging station location problem: A comprehensive review on models and solution approaches},
  author={Kchaou-Boujelben, Mouna},
  journal={Transportation Research Part C: Emerging Technologies},
  volume={132},
  pages={103376},
  year={2021},
  publisher={Elsevier},
  doi = {10.1016/j.trc.2021.103376}
}

@article{luo2020electric,
  title={Electric vehicle charging station location towards sustainable cities},
  author={Luo, Xiangyu and Qiu, Rui},
  journal={International Journal of Environmental Research and Public Health},
  volume={17},
  number={8},
  pages={2785},
  year={2020},
  publisher={MDPI},
  doi = {10.3390/ijerph17082785}
}

@article{chen2020optimal,
  title={Optimal charging facility location and capacity for electric vehicles considering route choice and charging time equilibrium},
  author={Chen, Rui and Qian, Xinwu and Miao, Lixin and Ukkusuri, Satish V},
  journal={Computers \& Operations Research},
  volume={113},
  pages={104776},
  year={2020},
  publisher={Elsevier},
  doi = {10.1016/j.cor.2019.104776}
}

@article{cui2019electric,
  title={Electric vehicle charging station placement method for urban areas},
  author={Cui, Qiushi and Weng, Yang and Tan, Chin-Woo},
  journal={IEEE Transactions on Smart Grid},
  volume={10},
  number={6},
  pages={6552--6565},
  year={2019},
  publisher={IEEE},
  doi = {10.1109/TSG.2019.2907262}
}

@article{erdougan2022establishing,
  title={Establishing a statewide electric vehicle charging station network in {Maryland}: A corridor-based station location problem},
  author={Erdo{\u{g}}an, Sevgi and {\c{C}}apar, {\.I}smail and {\c{C}}apar, {\.I}brahim and Nejad, Mohammad Motalleb},
  journal={Socio-Economic Planning Sciences},
  volume={79},
  pages={101127},
  year={2022},
  publisher={Elsevier},
  doi = {10.1016/j.seps.2021.101127}
}

@article{davidov2017stochastic,
  title={Stochastic expansion planning of the electric-drive vehicle charging infrastructure},
  author={Davidov, Sreten and Panto{\v{s}}, Milo{\v{s}}},
  journal={Energy},
  volume={141},
  pages={189--201},
  year={2017},
  publisher={Elsevier},
  doi = {10.1016/j.energy.2017.09.065}
}

@article{xi2013simulation,
  title={Simulation-optimization model for location of a public electric vehicle charging infrastructure},
  author={Xi, Xiaomin and Sioshansi, Ramteen and Marano, Vincenzo},
  journal={Transportation Research Part D: Transport and Environment},
  volume={22},
  pages={60--69},
  year={2013},
  publisher={Elsevier},
  doi = {10.1016/j.trd.2013.02.014}
}

@article{jordan2022electric,
  title={Electric vehicle charging stations emplacement using genetic algorithms and agent-based simulation},
  author={Jord{\'a}n, Jaume and Palanca, Javier and Mart{\'\i}, Pasqual and Julian, Vicente},
  journal={Expert Systems with Applications},
  volume={197},
  pages={116739},
  year={2022},
  publisher={Elsevier},
  doi = {10.1016/j.eswa.2022.116739}
}

@misc{fareye2024,
  title={Electric Vehicles for Last Mile Delivery: Future, Benefits \& Challenges},
  author={{FareEye}},
  year={2024},
  howpublished={\url{https://fareye.com/resources/blogs/electric-vehicles-for-last-mile-delivery}},
  note={Accessed: January 2024}
}

@misc{fedex2023,
  title={{FedEx} Deploys Electric Vehicles to Advance Sustainability Goal of Zero-Emissions Last-Mile Delivery in {India}},
  author={{FedEx Express}},
  year={2023},
  month={October},
  howpublished={\url{https://newsroom.fedex.com/newsroom/middle-east-indian-subcontinent-and-africa/fedex-deploys-electric-vehicles-to-advance-sustainability-goal-of-zero-emissions-last-mile-delivery-in-india}},
  note={Press Release}
}

@article{boyd2003subgradient,
  url = {https://web.stanford.edu/class/ee392o/subgrad_method.pdf},
  title={Subgradient methods},
  author={Boyd, Stephen and Xiao, Lin and Mutapcic, Almir},
  journal={Lecture Notes of EE392o, Stanford University, Autumn Quarter},
  year={2003}
}

@misc{nrel,
  doi = {10.2172/1604308},
title = "{R\&D} Insights for Extreme Fast Charging of Medium- and Heavy-Duty Vehicles: Insights from the {NREL} Commercial Vehicles and Extreme Fast Charging Research Needs Workshop, {August} 27-28, 2019",
author = "Kevin Walkowicz and Andrew Meintz and John Farrell",
year = "2020",
type = "Other",
}

@misc{EDF2024ChargingNeeds,
  author       = {{Environmental Defense Fund}},
  title        = {New study charts medium- and heavy-duty charging needs for 18 states and Washington, D.C.},
  howpublished = {\url{https://blogs.edf.org/energyexchange/2024/10/08/new-study-charts-medium-and-heavy-duty-charging-needs-for-18-states-and-washington-d-c/}},
  year         = {2024},
  month        = oct,
  note         = {Accessed: 2025-05-12}
}

@article{bertsimas1991stochastic,
  doi = {10.1287/opre.39.4.601},
  title={A stochastic and dynamic vehicle routing problem in the {Euclidean} plane},
  author={Bertsimas, Dimitris J and Van Ryzin, Garrett},
  journal={Operations Research},
  volume={39},
  number={4},
  pages={601--615},
  year={1991},
  publisher={INFORMS}
}

@article{brown2005statistical,
  doi = {10.1198/016214504000001808},
  title={Statistical analysis of a telephone call center: A queueing-science perspective},
  author={Brown, Lawrence and Gans, Noah and Mandelbaum, Avishai and Sakov, Anat and Shen, Haipeng and Zeltyn, Sergey and Zhao, Linda},
  journal={Journal of the American Statistical Association},
  volume={100},
  number={469},
  pages={36--50},
  year={2005},
  publisher={Taylor \& Francis}
}

@article{xie2018long,
  doi = {10.1016/j.tre.2017.11.014},
  title={Long-term strategic planning of inter-city fast charging infrastructure for battery electric vehicles},
  author={Xie, Fei and Liu, Changzheng and Li, Shengyin and Lin, Zhenhong and Huang, Yongxi},
  journal={Transportation Research Part E: Logistics and Transportation Review},
  volume={109},
  pages={261--276},
  year={2018},
  publisher={Elsevier}
}

@article{xie2021integrated,
  doi = {10.1016/j.apenergy.2021.117142},
  title={Integrated {U.S.} nationwide corridor charging infrastructure planning for mass electrification of inter-city trips},
  author={Xie, Fei and Lin, Zhenhong},
  journal={Applied Energy},
  volume={298},
  pages={117142},
  year={2021},
  publisher={Elsevier}
}

@article{kinay2023charging,
  doi ={10.1287/trsc.2021.0494},
  title={Charging station location and sizing for electric vehicles under congestion},
  author={K{\i}nay, {\"O}mer Burak and Gzara, Fatma and Alumur, Sibel A},
  journal={Transportation Science},
  volume={57},
  number={6},
  pages={1433--1451},
  year={2023},
  publisher={INFORMS}
}

@article{burke1956output,
  doi = {10.1287/opre.4.6.699},
  title={The output of a queuing system},
  author={Burke, Paul J},
  journal={Operations Research},
  volume={4},
  number={6},
  pages={699--704},
  year={1956},
  publisher={INFORMS}
}

@article{jackson1957networks,
  doi = {10.1287/opre.5.4.518},
  title={Networks of waiting lines},
  author={Jackson, James R},
  journal={Operations Research},
  volume={5},
  number={4},
  pages={518--521},
  year={1957},
  publisher={INFORMS}
}

@book{bose2013introduction,
  doi = {10.1007/978-1-4615-0001-8},
  title={{An Introduction to Queueing Systems}},
  author={{Bose, Sanjay K}},
  year={2013},
  publisher={Springer Science \& Business Media}
}

@article{albin1984approximating,
  title={{Approximating a point process by a renewal process, II: Superposition arrival processes to queues}},
  author={{Albin, Susan L}},
  journal={Operations Research},
  volume={32},
  number={5},
  pages={1133--1162},
  year={1984},
  publisher={INFORMS}
}

@book{medhi2002stochastic,
  title={{Stochastic models in queueing theory}},
  author={{Medhi, Jyotiprasad}},
  year={2002},
  publisher={Elsevier}
}

@book{whitt2002stochastic,
  title={{Stochastic-process limits: an introduction to stochastic-process limits and their application to queues}},
  author={{Whitt, Ward}},
  year={2002},
  publisher={Springer}
}

@book{bolch2006queueing,
  title={Queueing networks and Markov chains: modeling and performance evaluation with computer science applications},
  author={Bolch, Gunter and Greiner, Stefan and De Meer, Hermann and Trivedi, Kishor S},
  year={2006},
  edition={Second},
  publisher={John Wiley \& Sons}
}

@book{allen2014probability,
  title={Probability, statistics, and queueing theory},
  author={Allen, Arnold O},
  year={2014},
  publisher={Academic press}
}

@article{jung2014stochastic,
  title={Stochastic dynamic itinerary interception refueling location problem with queue delay for electric taxi charging stations},
  author={Jung, Jaeyoung and Chow, Joseph YJ and Jayakrishnan, R and Park, Ji Young},
  journal={Transportation Research Part C: Emerging Technologies},
  volume={40},
  pages={123--142},
  year={2014},
  publisher={Elsevier}
}

@article{cosmetatos1976some,
  title={{Some approximate equilibrium results for the multi-server queue (M/G/r)}},
  author={{Cosmetatos, George P}},
  journal={{Journal of the Operational Research Society}},
  volume={27},
  number={3},
  pages={615--620},
  year={1976},
  publisher={Taylor \& Francis}
}

\vfill
\framebox{\parbox{.90\linewidth}{\scriptsize The submitted manuscript has been created by
        UChicago Argonne, LLC, Operator of Argonne National Laboratory (``Argonne'').
        Argonne, a U.S.\ Department of Energy Office of Science laboratory, is operated
        under Contract No.\ DE-AC02-06CH11357.  The U.S.\ Government retains for itself,
        and others acting on its behalf, a paid-up nonexclusive, irrevocable worldwide
        license in said article to reproduce, prepare derivative works, distribute
        copies to the public, and perform publicly and display publicly, by or on
        behalf of the Government.  The Department of Energy will provide public access
        to these results of federally sponsored research in accordance with the DOE
        Public Access Plan \url{http://energy.gov/downloads/doe-public-access-plan}.}}
\end{document}